\documentclass[12pt]{article}
\usepackage{amsmath,amssymb,amsfonts,amsthm,bbm}
\usepackage{epic,eepic,epsfig,longtable}
\usepackage{multirow,verbatim}
\usepackage{array}
\usepackage{graphicx}
\usepackage{floatrow}
\usepackage{subcaption}
\usepackage{paralist}
\usepackage{latexsym}
\usepackage{comment}
\usepackage{epsfig}
\usepackage{setspace}
\usepackage{CJK}
\usepackage{color}

\usepackage{geometry}
\usepackage{algorithm}
\usepackage{algpseudocode}

\usepackage{rotating}

\usepackage[authoryear, round]{natbib}
\def\spacingset#1{\renewcommand{\baselinestretch}%
{#1}\small\normalsize} \spacingset{1}

\makeatletter
\def\singlespace{\def\baselinestretch{1}\@normalsize}

\numberwithin{equation}{section}

\renewcommand{\hat}{\widehat}

\renewcommand{\hat}{\widehat}

\newcommand{\bfm}[1]{\ensuremath{\mathbf{#1}}}

   \def\bA{\bfm A}  
   \def\bB{\bfm B}

    \def\FF{\mathbb{F}}

   \def\bI{\bfm I}  
   \def\bJ{\bfm J}  
     
   \def\bL{\bfm L}  
     
   \def\bN{\bfm N}

\def\bq{\bfm q}   \def\bQ{\bfm Q}

   \def\bW{\bfm W}  
   \def\bX{\bfm X}

\newcommand{\bfsym}[1]{\ensuremath{\boldsymbol{#1}}}

 \def\bgamma{\bfsym \gamma}             \def\bGamma{\bfsym \Gamma}

 \def\bmu{\bfsym {\mu}}                 
 
            \def\bTheta {\bfsym {\Theta}}
           
 \def\bsigma{\bfsym \sigma}             \def\bSigma{\bfsym \Sigma}

 \def\btau{\bfsym {\tau}}

\def\bvartheta{\bfsym{\vartheta}}	\def\bvarsigma{\boldsymbol{\varsigma}}

\DeclareMathOperator{\E}{E}

\def\newpage{\vfill\eject}

\def\today{\ifcase\month\or
  January\or February\or March\or April\or May\or June\or
  July\or August\or September\or October\or November\or December\fi
  \space\number\day, \number\year}

\newdimen\biblioindent    \biblioindent=30pt

 at 8truept

\newcommand{\beq}{\begin{equation}}
  \newcommand{\eeq}{\end{equation}}
\newcommand{\beqn}{\begin{eqnarray}}
  \newcommand{\eeqn}{\end{eqnarray}}
\newcommand{\beqnn}{\begin{eqnarray*}}
  \newcommand{\eeqnn}{\end{eqnarray*}}

\allowdisplaybreaks
\usepackage{verbatim}
\renewcommand{\baselinestretch}{1.66}
\def\tilde{\widetilde}

\def\FF{\mathcal{F}}
\def\[{\left [}  \def\]{\right ]} \def\({\left (}  \def\){\right )}
 \def\endpf{$\blacksquare$}
\def\hat{\widehat}

 \def \E {\mathrm{E}}    
 \def \1 {\mathbf{1}} 
\def \Y {Z}

\newtheorem{thm}{Theorem}

\newtheorem{lemma}{Lemma}

\newtheorem{remark}{Remark}
\theoremstyle{proposition}
\newtheorem{proposition}{Proposition}

\newtheorem{assumption}{Assumption}

\usepackage{lipsum}

\title{Adaptive Robust Large Volatility Matrix Estimation Based on High-Frequency Financial Data\footnote{
Minseok Shin is a Ph.D. student, College of Business, KAIST, Seoul 02455, South Korea.
Donggyu Kim is an Ewon Assistant Professor,  College of Business, KAIST, Seoul 02455, South Korea.
His research was supported by KAIST Basic Research Funds by Faculty (A0601003029) and the National Research Foundation of Korea (NRF) (2021R1C1C1003216).
Jianqing Fan is Frederick L. Moore'18 Professor of Finance, Department of Operations Research and Financial Engineering, Princeton University, Princeton, NJ 08544,  USA.  
His research was supported by NSFC grant No.71991471 and 71991470.
}
  }

\author{Minseok Shin$^1$,  Donggyu Kim$^1$\footnote{corresponding author. Tel: +82-02-958-3448. E-mail addresses: minseokshin@kaist.ac.kr (M. Shin),  donggyukim@kaist.ac.kr (D. Kim), jqfan@princeton.edu (J. Fan). }, and  Jianqing Fan$^{2,3, \dag}$ \\
$^1$KAIST, $^2$Capital  University of Economics and Business, and $^3$ Princeton University}

\begin{document}
\maketitle
\begin{spacing}{1.65}

\begin{abstract}
Several novel statistical methods have been developed to estimate large integrated volatility matrices based on high-frequency financial data.
To investigate their asymptotic behaviors, they require a sub-Gaussian or finite high-order moment assumption for observed log-returns, which cannot account for the heavy-tail phenomenon of stock-returns.
Recently, a robust estimator was developed to handle heavy-tailed distributions with some bounded fourth-moment assumption.
However, we often observe that log-returns have heavier tail distribution than the finite fourth-moment and that the degrees of heaviness of tails are heterogeneous across asset and over time.
In this paper, to deal with the heterogeneous heavy-tailed distributions,
we develop an adaptive robust integrated volatility estimator that employs pre-averaging and truncation schemes based on jump-diffusion processes.
We call this an adaptive robust pre-averaging realized volatility (ARP) estimator.
We show that the ARP estimator has a sub-Weibull tail concentration with only finite 2$\alpha$-th moments for any $\alpha>1$.
In addition, we establish matching upper and lower bounds to show that the ARP estimation procedure is optimal.
To estimate large integrated volatility matrices using the approximate factor model, the ARP estimator is further regularized using the principal orthogonal complement thresholding (POET) method.
The numerical study is conducted to check the finite sample performance of the ARP estimator.
\end{abstract}


\section{Introduction}
In modern financial studies and practices, volatility estimation is fundamental in risk management, performance evaluation, and portfolio allocation.
Due to the wide availability of high-frequency financial data, many well-performing volatility estimation methods have been developed to estimate integrated volatilities.
Examples include two-time scale realized volatility (TSRV) \citep{zhang2005tale}, multi-scale realized volatility (MSRV) \citep{zhang2006efficient, zhang2011estimating}, wavelet estimator \citep{fan2007multi},
pre-averaging realized volatility (PRV) \citep{christensen2010pre, jacod2009microstructure}, kernel realized volatility (KRV) \citep{barndorff2008designing, barndorff2011multivariate},
quasi-maximum likelihood estimator (QMLE) \citep{ait2010high, xiu2010quasi}, and the local method of moments \citep{bibinger2014estimating}.
One of the stylized facts of financial data is the existence of price jumps, and empirical studies have shown that the decomposition of daily variation into its continuous and jump components can better explain the volatility dynamics  \citep{ait2012testing,   andersen2007roughing, barndorff2006econometrics, corsi2010threshold, song2020volatility}.
For example,  \citet{fan2007multi} and \citet{zhang2016jump} employed the wavelet method to identify the jumps based on noisy high-frequency data.
\citet{mancini2004estimation} studied a threshold method for jump detection and presented the order of an optimal threshold, and \citet{davies2018data} further examined a data-driven threshold method.
 These estimation methods perform well for a small number of assets.
However, we often encounter a large number of assets in practices such as portfolio allocation, which results in the curse of dimensionality.
To overcome the curse of dimensionality and obtain an efficient and effective large volatility estimator, we often impose the approximate factor structure on the volatility matrix \citep{fan2018robust, fan2013large, fan2018eigenvector, kim2019factor}.
For example, to account for common market factors such as sector, firm size, and book-to-market ratios, the factor-based high-dimensional It\^o process is widely employed, and the idiosyncratic volatility is assumed to be sparse \citep{ait2017using,fan2016incorporating, fan2016overview,   kim2018Large, kong2018systematic}.
The  principal orthogonal complement thresholding (POET) method \citep{fan2013large} is often employed to estimate these low-rank plus sparse matrices.

\begin{figure}[!ht]
\centering
\includegraphics[width = 0.6 \textwidth]{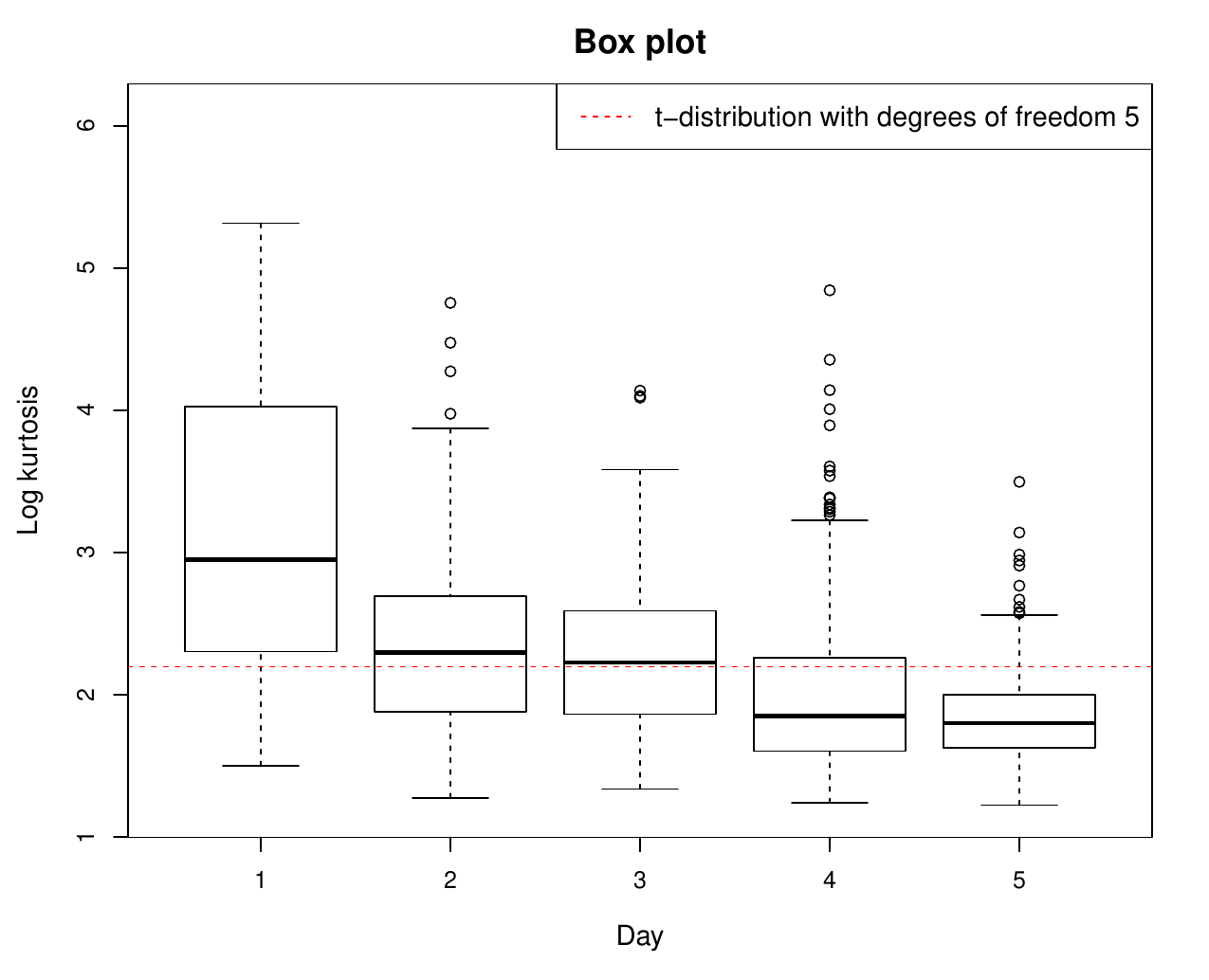}
\caption{The box plots for the daily distributions of log kurtoses calculated from 1-min log returns based on the most liquid 200 stocks on the S\&P 500 index.  The 5 selected box plots correspond to 5 days with the largest  IQR, the 75th, 50th, 25th, and 0th (smallest) percentile of the IQR among 501 trading days in 2015--2016, respectively.  The red dash marks the kurtosis for the $t_5$-distribution.  There are many assets whose 1-minute log-returns have empirical kurtosis larger than $t_5$-distribution and the numbers depend on the day of the market conditions.}\label{Fig-1}
\end{figure}

The performance of the factor-based large volatility matrix estimator critically depends on the accuracy of each element of the integrated volatility estimator.
Specifically, sub-Weibull tail concentration for the input volatility matrix estimator is required to investigate its asymptotic behaviors \citep{fan2013large, kim2019factor, kong2018systematic}.
However, one stylized fact of stock-return data is heavy-tailedness, which violates the sub-Gaussian assumption of the stock-return data.
Recently, with a bounded fourth-moment assumption on the microstructural noise, \citet{fan2018robust} developed a robust estimation method, which can attain sub-Gaussian tail concentration with the optimal convergence rate.
See also \citet{catoni2012challenging, minsker2018sub}.
However, we often observe that the bounded fourth-moment condition is violated. 
For example, Figure \ref{Fig-1} shows the five selected box plots of daily log kurtoses for the 1-minute log-returns of the 200 most liquid assets on the S\&P 500 index, calculated using the previous tick scheme.
These five selected plots correspond to the days with the largest interquartile range (IQR), the 75th, 50th, 25th, and 0th (smallest) percentile of the IQR among 501 trading days in 2015--2016.
From Figure \ref{Fig-1}, we find that the log-return data are heavy-tailed and also have heterogeneous degrees of heaviness of tails across the different assets and days.
See also \citet{cont2001empirical, mao2018stochastic, massacci2017tail}.
These facts generate the demand for developing an adaptive robust estimation method that can handle heterogeneous heavy-tailedness.

In this paper, we develop an adaptive robust integrated volatility estimator based on jump-diffusion processes contaminated by microstructural noises.
We first use the pre-averaging scheme  \citep{jacod2009microstructure} to adjust the unbalanced order relationship between the microstructural noises and true log-returns.
We then employ the truncation method \citep{minsker2018sub} using the daily moment conditions of assets.
Specifically, we truncate pre-averaged variables according to their heavy-tailedness, which allows for adaptive learning merits to be enjoyed.
Also, the truncation method sufficiently mitigates the effect of the jump signal on the pre-averaged variables.
We call the proposed estimator the adaptive robust pre-averaging realized volatility (ARP) estimator.
We show that the ARP estimator has sub-Weibull tail concentration, with finite 2$\alpha$-th moment assumption for any $\alpha>1$.
Also, by establishing matching upper and lower bounds for a pre-averaging estimator class, we show that the proposed estimator can achieve the optimal rate.
To make inferences about large integrated volatility matrices using the ARP estimator, we employ the factor-based high-dimensional It\^o process, which admits a low-rank plus sparse structure for the volatility matrix.
The ARP estimator is further regularized by the POET method, and we also investigate the theoretical properties of the POET estimator from the ARP estimator.
 Finally, we propose an estimation method for the tail index $\alpha$, resulting a full data-driven procedure, and investigate its properties.

The rest of the paper is organized as follows. Section \ref{SEC-2} provides the model and data structure for estimating each daily integrated volatility, and Section \ref{SEC-3} introduces an adaptive robust integrated volatility estimation method.
Section \ref{SEC-4} establishes the concentration property and optimality of the proposed estimator.
In Section \ref{SEC-5}, with the approximate factor model, we show how to estimate the large integrated volatility matrix using the proposed estimation procedure.  There, we  also propose the tail index estimator and show its concentration property. 
Section \ref{SEC-6} presents a simulation study to illustrate the finite sample performances of the proposed estimator and applies the estimation method to the  200 most-liquid assets on the S\&P 500.
The conclusion is presented in Section \ref{SEC-7}, and proofs are collected in the  Appendix.

\section{The model setup} \label{SEC-2}
We first introduce some notations.
For any given  $p_1$ by $p_2$ matrix $\bA = \left(A_{ij}\right)_{1 \leq i \leq p_1, 1 \leq j \leq p_2}$, let
 \begin{equation*}
  	\|\bA\|_1 = \max\limits_{1 \leq j \leq p_2}\sum\limits_{i = 1}^{p_1}|A_{ij}|,\hspace{0.5cm} \|\bA\|_\infty = \max\limits_{1 \leq i \leq p_1}\sum\limits_{j = 1}^{p_2}|A_{ij}| , \hspace{0.5cm}
    \| \bA \| _{\max} = \max_{i,j} | A_{ij}|.
 \end{equation*}
The matrix spectral norm $\|\bA\|_2$ is  the square root of the largest eigenvalue of $\bA\bA^\top$,
and the Frobenius norm of $\bA$ is denoted by $\|\bA\|_F = \sqrt{ \mathrm{tr}(\bA^{\top} \bA) }$.
We will use $C$ to denote a generic positive constant whose value is free of $n$ and $p$ and may vary between appearances.

Let $\bX(t)= (X_1 (t), \ldots, X_p(t))^\top$ be the vector of true log-prices for $p$ assets at time $t$.
To model the high-frequency financial data, we often employ the jump-diffusion process as follows:
\begin{eqnarray} \label{eq-2.1}
 d \bX(t) &=& d \bX^{c}(t)+\bJ(t) \cr
 		&=& \bmu (t)dt +  \bsigma ^\top(t)d\bW_t+\bJ(t),
\end{eqnarray}
where $\bX^{c}(t)= (X_1^{c}(t), \ldots, X_p^{c}(t))^\top$  with $\bX^{c}(0)=\bX(0)$  is the vector of true continuous log-prices at time $t$, $\bmu(t) = (\mu_1(t), \ldots, \mu_p(t))^\top$ is a drift vector, $\bsigma(t)$ is a $q$ by $p$ matrix, $\bW_t$ is a $q$ dimensional independent Brownian motion, and the stochastic processes $\bmu(t)$, $\bX(t)$, $\bX^{c}(t)$, and $\bsigma(t)$ are defined on a filtered probability space $(\Omega, \FF, \{\FF_t, t \in [0, 1]\},  P)$ with filtration $\FF_t$ satisfying the usual conditions.
 For the jump part, $\bJ(t) = (J_1(t), \ldots, J_p(t))^\top$ denotes the jump process and each $J_i(t)$ is decomposed as 
 \begin{equation*}
J_i(t) = J_{i1}(t) + J_{i2}(t),
\end{equation*}
where 
 \begin{equation*}
 J_{i1}(t) =  \int_{0}^{t} \int_{|x| > 1} x N_i(dx,ds)
\end{equation*}
is the sum of the big jumps, $N_i$ is a Poisson random measure defined on $\mathbb{R} \times [0, 1]$ with intensity measure $I_i(dx)ds$,
 \begin{equation*}
 J_{i2}(t) =  \int_{0}^{t} \int_{|x| \leq 1} x[N_i(dx,ds) - I_i(dx)ds ]
\end{equation*}
is the compensated sum of small jumps, and $N_i(dx,ds) - I_i(dx)ds$ is the compensated Poisson random measure. 
We assume that $J_{i1}(t)$ is a compound Poisson process with the bounded intensity and $J_{i2}(t)$ is an independent L\'evy process.
For the small jumps, we define the Blumenthal-Gettor index $\pi$ as
\begin{equation*}
\pi =\inf \{\delta \geq 0,  \quad \int_{|x| \leq 1} |x|^{\delta} I_i(dx) < \infty \quad \text{for all }i=1, \ldots, p \}.
\end{equation*}
We note that the Blumenthal-Gettor index $\pi$ measures the activity of small jumps.
We describe the condition for $\pi$ in Assumption \ref{Assumption1}.
This kind of jump processes is widely used to explain price jumps \citep{ait2020high, andersen2021volatility, cont2007nonparametric, li2017adaptive, mancini2009non}. 
The instantaneous volatility matrix of the continuous log-price $\bX^c(t)$ is
\begin{equation*}
\bgamma(t) = \left(\gamma_{ij}(t)\right)_{1\leq i,j\leq p} =  \bsigma^\top(t)\bsigma(t),
\end{equation*}
and their  quadratic variation  is
\begin{eqnarray*}
 \left[\bX^c, \bX^c \right]_t &=& \int_0^t\bgamma(s) ds = \left( \int_0^t \gamma  _{ij}(s) ds\right)_{1 \leq i,j \leq p}  \cr
	&=&  \int_{0}^t \bsigma^\top(s)\bsigma(s)ds .
\end{eqnarray*}
The parameter of interest is the integrated volatility matrix of $\bX^c(t)$,
\begin{equation} \label{eq-2.2}
	\bGamma = \[ \bX^c, \bX^c \] _1 = \int_{0}^1 \bgamma(s) ds.
\end{equation}

Unfortunately, we cannot directly observe the true log-prices $\bX(t)$.
In fact, observed high-frequency data are contaminated by microstructural noises.
Furthermore, high-frequency data encounter a non-synchronization problem, in which transactions for multiple assets often arrive asynchronously.
In this regard, we assume that the observed log-price $Y_i(t_{i,k})$ obeys the following model:
\begin{equation}  \label{eq-2.3}
Y_i (t_{i,k}) = X_i(t_{i,k}) + \epsilon_{i}(t_{i,k}) \quad  \text{ for } i = 1, \ldots, p, k = 0, \ldots, n_i ,
\end{equation}
where $t_{i,k}$ is the $k$-th observation time point of the $i$-th asset, and for each given $i=1, \ldots, p,$ $\epsilon_{i}(t_{i,k}), k=0, \ldots, n_i,$ are i.i.d. noises with a mean of zero.
We assume that for $i,j=1, \ldots, p,$ $\E\left[ \epsilon _{i}\left( t \right) \epsilon _{j}\left( t \right) \right]=\eta _{ij}$ and $\epsilon_{i}(t)$ is independent of  $\epsilon_{j}(t')$ for $t\neq t'$.

To handle the microstructural noise issue, several estimation methods have been developed \citep{ ait2010high,  barndorff2008designing, barndorff2011multivariate, bibinger2014estimating, christensen2010pre,   fan2007multi,  jacod2009microstructure, xiu2010quasi,zhang2005tale,  zhang2006efficient, zhang2011estimating}.
They work well for a finite number of assets and are widely adopted to develop large volatility matrix estimation procedures \citep{kim2018adaptive, kim2016asymptotic,  wang2010vast}.
However, the observed log-prices are heavy-tailed, so these methods cannot lead to the estimators with the sub-Weibull concentration bound, which is essential for making inferences about the asymptotic behaviors of the large volatility matrix estimation procedures.
To tackle the heavy-tail issue, \citet{fan2018robust} proposed the robust pre-averaging realized volatility estimation procedure, which can achieve the sub-Gaussian tail concentration with only the finite fourth-moment condition on the microstructural noise.
However, as shown in Figure \ref{Fig-1}, the degrees of heaviness of tails of assets' log-returns is heterogeneous across assets and over time.
Furthermore, jumps in the true log-price process can also cause heavy-tailed distributions.
To account for these features,  we accommodate heterogeneous degrees of tail distributions based on the jump-diffusion process contaminated by microstructural noises.
We allow  each asset to have its own order of the highest finite absolute moment (see Assumption \ref{Assumption1} in Section \ref{SEC-4} for details).

\section{Adaptive robust pre-averaging realized volatility estimator} \label{SEC-3}
In this section, we introduce an adaptive robust integrated volatility estimation procedure to handle non-synchronization, price jump, and microstructural noise.
To handle the non-synchronization problem, we consider the generalized sampling time proposed by \citet{ait2010high}.
We note that the generalized sampling time scheme includes synchronization schemes such as previous tick \citep{zhang2011estimating, wang2010vast} and refresh time \citep{barndorff2011multivariate, fan2012vast}.
See also  \citet{bibinger2014estimating, chen2019five, fan2019structured, hayashi2005covariance, hayashi2011nonsynchronous,  malliavin2009fourier, park2016estimating}.
We define the generalized sampling time as follows.

\paragraph{Definition 1.} \citep{ait2010high}. A sequence of time points $\btau=\left\{\tau_0, \ldots, \tau_n \right\}$ is said to be the generalized sampling time if
\begin{enumerate}
\item [(1)] $0=\tau_0<\tau_1<\cdots<\tau_{n-1}<\tau_{n}=1$;
\item [(2)] there exists at least one observation for each asset between consecutive $\tau_j$'s; and
\item [(3)] the time intervals, $\left\{\Delta_{j}=\tau_j-\tau_{j-1}; j =1,\ldots, n \right\}$, satisfy
$\sup_j \Delta_{j}\xrightarrow{p} 0$.
\end{enumerate}
Examples of generalized sampling scheme include the previous tick approach discussed in \cite{zhang2011estimating} and the refresh time scheme in \cite{barndorff2011multivariate}.
For the $i$-th asset, we select arbitrary observation, $Y_i (\tau_{i,k})$, between $\tau_{k-1}$ and $\tau_{k}$.
In other words, we choose any $\tau_{i,k}\in ( \tau _{k-1},\tau _{k}] \cap \{ t_{i,l},l= 0,1,\ldots ,n_{i}\}$ for $i=1, \ldots, p$.

Based on synchronized time, $\btau$, we adopt the pre-averaging method to manage  the microstructural noise \citep{jacod2009microstructure}.
For the observed log-returns, $Y_i (\tau_{i,k+1}) - Y_i (\tau_{i,k}), i=1,\ldots, p, k=1, \ldots, n-1$, the variance of the microstructural noise $2 \eta_{ii}$ dominates the continuous log-return volatility $\int_{\tau_{i,k}}^{\tau_{i,k+1}} \gamma_{ii}(t)dt$.
Therefore, it is difficult to estimate the integrated volatility without smoothing for denoising.
To adjust the order relationship between the noises and continuous log-returns, we use the following pre-averaged data to suppress the noises \citep{christensen2010pre, jacod2009microstructure}:
\begin{equation}\label{eq-3.1}
 \Y_i(\tau_{k}) = \sum _{l=0}^{K_n-1} g\left(\frac{l}{K_n}\right)\left  \{ Y_i (\tau_{i,k+l+1}) - Y_i (\tau_{i,k+l}) \right \}  \quad \text{ for } i=1, \ldots, p , k=1, \ldots, n-K_n,
\end{equation}
where the  weight function $g(\cdot)$ is continuous and piecewise continuously differentiable with a piecewise Lipschitz derivative $g^\prime$ and satisfies $g(0) = g(1) = 0$ and $\int_0^1\{g(t)\}^2dt > 0$.
In this paper, we choose bandwidth parameter $K_n$ as $C_{K}n^{1/2}$ for some constants $C_{K}$, which yields the optimal rate $n^{-1/4}$ for the resulting estimator.
Then, the continuous log-returns and noises in $\Y_i(\tau_{k})$'s are of the same order of magnitude  \citep{fan2018robust}.
However,  the pre-averaged random variables still have heterogeneous heavy-tails across assets.
Furthermore, there exist jump variations in the pre-averaged data, which create outliers.
To handle these problems, we robustly estimate the volatility matrix by applying an adaptive truncation method according to the tails of the continuous part of the data.

We define the quadratic pre-averaged random variables
\begin{equation}  \label{eq-3.2}
Q_{ij}(\tau_k) = \frac{n-K_n}{ \phi K_n}\Y_i(\tau_k)\Y_j(\tau_k) \quad \text{ for } i,j=1,\ldots, p , k=1, \ldots, n-K_n,
\end{equation}
where $\phi = \frac{1}{K_n} \sum_{\ell = 0 }^{K_n-1} \left \{ g \( \frac{\ell}{K_n}\) \right \}^2 $, and  let
\begin{equation} \label{eq-3.3}
	\alpha_{ij}=2\wedge \dfrac{2\alpha_{i}\alpha_{j}}{\alpha_{i}+\alpha_{j}},
\end{equation}
where   $\alpha_{i}$ is the order of the highest finite moment for the continuous part of  $Q_{ii}(\tau_k)$ (see Assumption \ref{Assumption1} in Section \ref{SEC-4}). 
 We note that the harmonic mean provides that the continuous part of  $Q_{ij}(\tau_k)$ has the finite $\dfrac {2\alpha_{i}\alpha_{j}}{\alpha_{i}+\alpha_{j}}$-th moment (see Proposition \ref{Proposition3}(a) in the Appendix).
Furthermore, since we can obtain the optimal convergence rate under the finite second moment condition for the continuous part of  $Q_{ij}(\tau_k)$ \citep{fan2018robust} even when it has higher moments, we set $\alpha_{ij} \leq 2$.
Therefore, we define $\alpha_{ij}$ as  \eqref{eq-3.3}.

Then, to handle the heterogeneous heavy-tails, we propose the following adaptive truncation method:
\begin{equation} \label{eq-3.4}
	\hat{T}^{\alpha }_{ij, \theta} =  \dfrac {1}{(n-K_{n})\theta_{ij}}\sum_{k=1}^{n-K_n}\psi_{\alpha_{ij}}\left\{\theta _{ij}Q_{ij}\left(\tau_{k}\right)\right\},
\end{equation}
where $\theta _{ij}$ is a truncation parameter specified in Theorem \ref{Theorem1} and implemented by \eqref{eq-5.17}, and $\psi _{\alpha }\left( x\right)$ is a bounded non-decreasing function defined for $\alpha  \in ( 1,  2]$ as follows:
\begin{equation*}
	\psi_{\alpha} \left( x\right)=
 \begin{cases}
-\log (1- t_{\alpha} +c_{\alpha}t_{\alpha}^{\alpha})   & \text{ if } x\geq t_{\alpha}   \\
-\log (1- x +c_{\alpha}x^{\alpha}) & \text{ if } 0\leq x\leq t_{\alpha} \\
\log (1+ x + c_{\alpha}|x|^{\alpha}) & \text{ if } -t_{\alpha}\leq x\leq 0 \\
\log (1 -t_{\alpha} + c_{\alpha}t_{\alpha}^{\alpha}) &  \text{ if }x\leq -t_{\alpha},\\
\end{cases}
\end{equation*}
where  $c_{\alpha}=\max \left \{(\alpha -1)/\alpha , \sqrt {(2-\alpha)/\alpha }\right \}$, and $t_{\alpha}=\left( 1/\alpha c_{\alpha}\right) ^{1/(\alpha -1)}.$
 We note that to obtain the optimal sub-Weibull tail concentration for the ARP estimator, the inequality $-\log (1- x +c_{\alpha}|x|^{\alpha}) \leq \psi _{\alpha}\left( x\right) \leq \log (1+ x + c_{\alpha}|x|^{\alpha})$ is required. 
Since the term $c_{\alpha}$ appears in the upper bound, we need to choose $c_{\alpha}$ as small as possible  under the  following condition  
\begin{equation*}
1+ x + c_{\alpha}|x|^{\alpha} > 0  \,\,\, \text{and} \,\,\,  -\log (1- x +c_{\alpha}|x|^{\alpha}) \leq \log (1+ x + c_{\alpha}|x|^{\alpha})
\end{equation*}
for all $x \in \mathbb{R}$.
To satisfy the above inequality and obtain the sharp bound of the  concentration inequality, we choose $c_{\alpha}=\max \left \{(\alpha -1)/\alpha , \sqrt {(2-\alpha)/\alpha }\right \}$. 
On the other hand, since $f(x)=1+ x + c_{\alpha}|x|^{\alpha}$ has the minimum at $x= -\left( 1/\alpha c_{\alpha}\right) ^{1/(\alpha -1)}$, to make $\psi _{\alpha}\left( x\right)$ non-decreasing, we choose  $t_{\alpha}=\left( 1/\alpha c_{\alpha}\right) ^{1/(\alpha -1)}$. 
 This proposed truncation is also able to mitigate the impact of jumps since the impact of each jump for the truncated quadratic pre-averaged random variables is bounded by the truncation level $CK_n/(n \theta_{ij})$.
In this point of view, the proposed approach is different from the literature that handles the jumps by finding the exact jump locations.
For example, \citet{cont2007nonparametric} and \citet{mancini2009non} first detected the jump locations and estimated the integrated volatility only using the time intervals with no jumps.
Their theorems rely on the condition that the increments of the continuous part is bounded.
However, in this paper, we assume that the microstructural noise has the heavy tails, which makes it hard to distinguish the jumps and outliers from the heavy-tailed microstructural noise. 
Thus, to handle this issue, we first determine the truncation parameter according to the moment condition of the microstructural noises and then control the impact of jumps using the selected truncation parameter. 
We note that other truncation methods can also achieve a similar goal (see \citet{fan2016shrinkage}).
It will be shown that the proposed adaptive robust estimator $\hat{T}^{\alpha}_{ij, \theta}$ possesses the sub-Weibull concentration bounds (see Theorem \ref{Theorem1}).
Furthermore, by employing the tick-by-tick overlapping scheme, we can mitigate the effect of the irregular observation time point, which can cause the inconsistency for the pre-averaging estimator without the tick-by-tick overlapping \citep{mykland2019algebra}.

The adaptive robust estimator $\hat{T}^{\alpha}_{ij, \theta}$ is, however, not a consistent estimator of the true integrated volatility  $\Gamma_{ij}$, because the noises still remain in each $Q_{ij}(\tau_k)$.  Indeed, it will be shown that  $\hat{T}^{\alpha}_{ij, \theta}$ converges to
\begin{equation} \label{eq-3.5}
 T_{ij} = \Gamma_{ij}   + \rho_{ij},
\end{equation}
where
\begin{equation*}
 \rho_{ij}=\frac{\sum_{k=1}^{n}\mathbf{1}(\tau_{i,k}=\tau_{j,k})}{ \phi K_n}\zeta   \eta_{ij}, \quad \zeta = \sum\limits_{l = 0}^{K_n-1}\left \{ g\left(\frac{l}{K_n}\right) - g\left(\frac{l+1}{K_n}\right)\right\} ^2 = O\(\frac{1}{K_n}\),
\end{equation*}
with the covariance of noise $\eta_{ij}$ defined in \eqref{eq-2.3}, and $\mathbf{1}(\cdot)$ is the indicator function.
Hence, to estimate the integrated volatility $\Gamma_{ij}$, we adjust $\hat{T}^{\alpha}_{ij, \theta}$ by subtracting an estimator of $\rho_{ij}$.
For this purpose, let us first define an adaptive robust estimator, $\hat{\rho}^{\alpha }_{ij, \theta}$, as
\begin{equation}  \label{eq-3.6}
	\hat{\rho}^{\alpha }_{ij, \theta} =  \dfrac {\zeta}{\phi K_n \theta_{\rho,ij}}\sum_{k=1}^{n-1}\psi_{\alpha_{ij}}\left\{\theta _{\rho,ij}Q_{\rho,ij}\left(\tau_{k}\right)\right\},
\end{equation}
where \begin{equation}  \label{eq-3.7}
Q_{\rho, ij}(\tau_k) =\frac{1}{2 }  \{ Y_i(\tau_{i,k+1})- Y_i(\tau_{i,k})\} \{ Y_j(\tau_{j,k+1})- Y_j(\tau_{j,k})\}
\end{equation}
for $i,j=1, \ldots, p , k=1, \ldots, n-1,$ and $\theta _{\rho,ij}$ is truncation parameter that will be specified in Theorem \ref{Theorem3}.
We now define the integrated volatility estimator as follows:
\begin{equation}  \label{eq-3.8}
\hat{\Gamma}_{ij} ^{\alpha}= \hat{T}^{\alpha}_{ij, \theta} - \hat{\rho}^{\alpha }_{ij, \theta}.
\end{equation}
We call this the adaptive robust pre-averaging realized volatility (ARP) estimator.   
This provides a preliminary consistent estimate of $\hat{\Gamma}_{ij}$, which will be further regularized.

\section{Theoretical properties of the ARP estimator} \label{SEC-4}
In this section, we show the concentration property and optimality of the ARP estimator  by establishing matching upper and lower bounds for both $\hat{T}^{\alpha}_{ij, \theta}$ and $\hat{\rho}^{\alpha }_{ij, \theta}$.
We define
\begin{eqnarray*}
	&&Y_i^{c} (\tau_{i,k}) = X_i^{c}(\tau_{i,k}) + \epsilon_{i}(\tau_{i,k}), \quad   \Y_i^{c}(\tau_{k}) = \sum _{l=0}^{K_n-1} g\left(\frac{l}{K_n}\right)\left  \{ Y_i^{c} (\tau_{i,k+l+1}) - Y_i^{c} (\tau_{i,k+l}) \right \}, \cr
 &&Q_{ij}^{c}(\tau_k) = \frac{n-K_n}{ \phi K_n}\Y_i^{c}(\tau_k)\Y_j^{c}(\tau_k),  \quad \text{and}\cr
 && Q_{\rho, ij}^{c}(\tau_k) =\frac{1}{2 }  \{ Y_i^{c}(\tau_{i,k+1})- Y_i^{c}(\tau_{i,k})\} \{ Y_j^{c}(\tau_{j,k+1})- Y_j^{c}(\tau_{j,k})\},
\end{eqnarray*}
where $X_i^{c}(t)$ is the true continuous log-price process defined in \eqref{eq-2.1}, and the superscript $c$ represents the continuous part of the true log-price.
Now, to investigate asymptotic properties of $\hat{T}^{\alpha}_{ij, \theta}$, we make the following assumptions.

\begin{assumption} \label{Assumption1}
~
\begin{enumerate}
\item [(a)] There exist positive constants $\nu_\mu$ and $\nu_\gamma$  such  that
\begin{equation*}
	\max_{1\leq i \leq p}\sup_{0\leq t \leq 1}  | \mu_{i}(t)| \leq \nu_\mu \text{ a.s.}, \text{ and } \; \max_{1\leq i \leq p}\sup_{0\leq t \leq 1}    \gamma_{ii} (t) \leq \nu_\gamma \text{ a.s.}
\end{equation*}
\item [(b)] The generalized sampling time $\left\{\tau_j\right\}$ is independent of the true log-price related  processes defined in \eqref{eq-2.1}  and the noise $\epsilon_i(t_{i,k})$.
The time intervals, $\left\{\Delta_{j}=\tau_j-\tau_{j-1}, \text{ } 1\leq j \leq n \right\}$, satisfy
\begin{equation*}
	\max_{1\leq j \leq n} \Delta_{j} \leq Cn^{-1}  \text{ a.s.}
\end{equation*}
\item [(c)]There exist positive constants, $\nu_Q$ and $\alpha_1, \ldots, \alpha_p>1$, such that for all $1\leq k \leq n-K_{n}$,
\begin{equation*}
\max_{1\leq i \leq p} \E \left\{  \left |Q^{c}_{ii}(\tau_{k})    \right  | ^{\alpha_{i}} \right\} \leq \nu_Q.
\end{equation*}

\item [(d)] For $i,j=1, \ldots, p,$
$\E\left[ \epsilon _{i}\left( t_{i,k}  \right) \epsilon _{j}\left( t_{j,k^{\prime}} \right) \right]=\eta _{ij}$ for $t_{i,k} = t_{j, k^{\prime}}$,
$\epsilon_{i}(t_{i,k})$ is independent of  $\epsilon_{j}(t_{j,k^{\prime}})$ for $t_{i,k} \neq t_{j,k^{\prime}}$, and $\epsilon_{i}(\cdot)$ is independent of the true log-price related processes.

\item [(e)] The Blumenthal-Gettor index $\pi$ satisfies  $\pi  \in [ 0, 2)$, and we have
\begin{equation*} 
I_i(|x|> n^{(1-\alpha_{i})/(4\alpha_{i}-2\alpha_{i} \pi)}) < \infty \,\, \text{ for }i=1, \ldots, p.
\end{equation*}
\end{enumerate}
\end{assumption}
\begin{remark}
For Assumption \ref{Assumption1}(a), the boundedness condition of the instantaneous volatility process $\gamma_{ii}(t)$ can be relaxed to the locally boundedness condition when we investigate the asymptotic behaviors of volatility estimators, such as their convergence rate (see \citet{ait2017using}).
Specifically, Lemma 4.4.9 in \citet{jacod2012discretization} indicates that if the asymptotic result, such as convergence in probability or stable convergence in law, is satisfied under the boundedness condition, it is also satisfied under the locally boundedness condition.
From this point of view, because we consider a finite time period, it is sufficient to investigate the asymptotic properties under the boundedness condition.
Thus, Assumption \ref{Assumption1}(a) is not restrictive.
\end{remark}

\begin{remark}
	First of all, the finite 2$\alpha_i$-th moment condition for the microstructural noises with  Assumption \ref{Assumption1}(a),(b),(d) implies  Assumption \ref{Assumption1}(c).   
    It is the finite moment condition, which entails that the quadratic pre-averaged variable, ${Q}^{c}_{ij}(\tau_k)$, for the continuous part satisfies
 \begin{equation}\label{eq-4.1}
\E \left\{  \left |{Q}^{c}_{ij}(\tau_k)    \right  | ^{\alpha_{ij}} \Big | \FF_{\tau_k} \right\} \leq C \text{ a.s.} 
\end{equation} 
for all $1\leq i,j \leq p$ and $1\leq k \leq n-K_{n}$,  where $\alpha_{ij}$ is defined in \eqref{eq-3.3} (see Proposition \ref{Proposition3}(a) in the Appendix).
To account for the heterogeneous heavy-tailedness, we allow the tail index $\alpha_{i}$ to vary from 1 to infinity.
If $\alpha_{i}=2$ for all $i =1,\ldots,p$, it is the similar setting as that of  \citet{fan2018robust}, and the ARP estimator has universal truncation, which we call the universal robust pre-averaging realized volatility (URP) estimator.
The main difficulty in extending the theory for the URP estimator to that for the ARP estimator is that  $\alpha_{i}$'s may not be integers.
Thus, to obtain the optimal convergence rate, we need to choose appropriate truncation function and truncation parameters and carefully derive the upper bounds.
We note that when $\alpha_i < 2$, the URP estimator needs to choose the smallest tail index for all tails.
 This choice cannot provide the optimal convergence rate for the heterogeneous heavy-tailedness. 
On the other hand, to implement the ARP procedure, we need to estimate the tail indices, which may cause some estimation errors.
To overcome this problem, we proposed the novel tail index estimation procedure and show its concentration properties in Section \ref{SEC-5.2}. 
To investigate the heterogeneous heavy-tail, we compare the ARP and URP estimators in the numerical studies.

\end{remark}

\begin{remark}
  Assumption \ref{Assumption1}(e) allows the infinity number of small jumps. 
  However, we additionally need the finite number of jumps whose sizes are bigger than $n^{(1-\alpha_{i})/(4\alpha_{i}-2\alpha_{i} \pi)}$ due to the outliers coming from the heavy-tailed observations.
  Specifically, since the microstructural noises have the heavy tails,  $Q^{c}_{ii}(\tau_{k})$  can have large values that are comparable to the jumps. 
  To handle this issue, we need  the condition $I_i(|x|> n^{(1-\alpha_{i})/(4\alpha_{i}-2\alpha_{i} \pi)}) < \infty$.
  We note that if $\pi=0$, Assumption \ref{Assumption1}(e) is not required. 
  Furthermore, to investigate asymptotic theorems for the tail index estimation in Section \ref{SEC-5.2}, we need  the condition $\pi=0$ (see Assumption \ref{Assumption4}(b)). 
\end{remark}

\begin{remark}
To obtain the sub-Weibull tail concentration, we technically need some bounded condition for the target parameter, such as Assumption \ref{Assumption1}(a). 
However, if we change the parameter of interest to the sum of conditional expected values as follows:
\begin{eqnarray*}  
	\bGamma^{new} &=& \sum_{k=0}^{n-1} \E \left \{ \int_{\tau_k} ^{\tau_{k+1}} \bgamma(t) dt  \middle | \FF_{\tau_k} \right \}   
\end{eqnarray*}
and impose the bounded condition for the conditional expected values, namely,
$$
\max_{1\leq i  \leq p} \max_{0 \leq k \leq n-1} \E \[  \left\{(\tau_{k+1}- \tau_k)^{-1} \int_{\tau_k} ^{\tau_{k+1}} \gamma_{ii} (t)dt  \right \} ^{\alpha_{i}}  \middle | \FF_{\tau_k} \] \leq \nu_\gamma   \text{ a.s.} ,
 $$
we can obtain the same result for the new target parameter $	\bGamma^{new}$.
 Under this condition, the random fluctuation of the instantaneous volatility process can be the source of the heavy-tailedness. 
  Then, the proposed procedure can estimate the new target parameter $\bGamma^{new} $ with the optimal rate $n^{(1-\alpha_{ij})/ 2\alpha_{ij}}$. 
\end{remark}

The theorem below shows that $ \hat{T}^{\alpha}_{ij,\theta}$ has the sub-Weibull tail concentration with a convergence rate of $n^{(1-\alpha_{ij})/2\alpha_{ij}}$.

\begin{thm}\label{Theorem1}
(Upper bound)
Under the models \eqref{eq-2.1} and \eqref{eq-2.3} and Assumption \ref{Assumption1}, let $  \delta ^{-1} \in [ n^{c}, e^{n^{1/2}} ]$ for some positive constant $c >0$.
Take
\begin{equation*}
\theta_{ij} =    \left(\frac {K_{n}\log \left(3K_{n}^{2}\delta ^{-1}\right)}{\left( \alpha_{ij} -1\right) c_{\alpha_{ij} }S_{ij}(n-K_{n})} \right)^{1/\alpha_{ij}},
\end{equation*}
where $S_{ij} = \dfrac {1}{n-K_{n}}\sum ^{n-K_{n}}_{k=1}U_{ij}\left( \tau_{k}\right)$, and $U_{ij}\left( \tau_{k}\right)$'s are some positive constants defined in Proposition \ref{Proposition3}(a).
Then, we have, for a sufficiently large $n$,
\begin{equation}\label{eq-4.2}
		 \Pr  \left \{   | \hat{T}^{\alpha}_{ij,\theta} - T_{ij} | \leq  C  \left( n^{-1/2}\log \delta ^{-1}\right) ^{(\alpha_{ij} -1)/\alpha_{ij} }
		  \right \}  \geq  1-\delta.
\end{equation}
\end{thm}

Theorem \ref{Theorem1} indicates that  $\hat{T}^{\alpha}_{ij,\theta}$ has the sub-Weibull concentration bound with the convergence rate of $n^{(1-\alpha_{ij})/2\alpha_{ij}}$.
Specifically, as long as the number of assets, $p$, satisfies $  p^{b+2} \in [ n^{c}, e^{n^{1/2}} ]$ for some positive constant $c >0$, we have
\begin{equation*}
		 \Pr  \left \{  \max_{1 \leq i,j \leq p}  | \hat{T}^{\alpha}_{ij,\theta} - T_{ij} | \geq  C_b  \( n^{-1/2}\log p \) ^{(\alpha -1)/\alpha }
		   \right \}  \leq  p^{-b}
\end{equation*}
for any constant $b>0$ and $\alpha=\min_{1\leq i \leq p} \alpha_{i}$, where $C_b$ is some constant depending on $b$, which is the essential condition for investigating inferences of large integrated volatility matrix (see Proposition \ref{Proposition1}).
An interesting finding is that there is a trade-off between the convergence rate $n^{(1-\alpha_{ij})/2\alpha_{ij}}$ and the tail indices $\alpha_{i}$ and $\alpha_{j}$.
This raises the question of whether the upper bound in \eqref{eq-4.2} is optimal.

Let $\hat{T}_{ij}\left( Q_{ij}\left( \tau_{k}\right) , \delta \right) = \hat{T}_{ij}\left( Q_{ij}\left( \tau_{1}\right), \ldots, Q_{ij}\left( \tau_{n-K_{n}}\right) , \delta \right)$ be any pre-averaging estimator for $T_{ij}$ defined in \eqref{eq-3.5}, which takes the values of pre-averaged variables $Q_{ij}\left( \tau_{k}\right), k=1, \ldots, n-K_n$, defined in \eqref{eq-3.2}.
The following theorem establishes the lower bound for the maximum concentration probability among the class of pre-averaging estimators $\hat{T}_{ij}\left( Q_{ij}\left( \tau_{k}\right) , \delta \right)$ which satisfy $\max_{1\leq i \leq p} \E \left\{  \left |Q^{c}_{ii}(\tau_{k})    \right  | ^{\alpha_{i}} \right\} \leq C$ for all $1\leq k \leq n-K_{n}$.

\begin{thm}\label{Theorem2}
(Lower bound)
Under the assumptions in Theorem \ref{Theorem1}, let $\alpha_{ij} \in \left(1,2\right)$ for some $1\leq i,j \leq p$.
Then, we have, for a sufficiently large $n$, 
\begin{equation}\label{eq-4.3}
		\min _{\hat{T}_{ij}\left( Q_{ij}\left( \tau_{k}\right) ,\delta \right)} \max_{\bQ^c \in \Xi (\alpha_1,\ldots, \alpha_p)} \Pr  \left \{   | \hat{T}_{ij}\left( Q_{ij}\left( \tau_{k}\right) ,\delta \right) - T_{ij} | \geq  C  \left( n^{-1/2}\log \delta ^{-1} \right) ^{(\alpha_{ij} -1)/\alpha_{ij} }
		 \right \}  \geq  \delta,
\end{equation}
where 
\begin{equation*}
\Xi (\alpha_1,\ldots, \alpha_p) =\{\bQ^c = (Q^{c}_{ii}(\tau_{k}))_{i=1,\ldots,p, k=1,\ldots, n-K_n} : \max_{i, k} \E \left\{  \left |Q^{c}_{ii}(\tau_{k})    \right  | ^{\alpha_{i}} \right\} \leq C \}.
\end{equation*}
\end{thm}

Theorem \ref{Theorem2} shows that  the lower bound is $n^{(1-\alpha_{ij})/2\alpha_{ij}}$, which matches the upper bound  in Theorem \ref{Theorem1}.
Thus, the proposed estimator obtains the optimal convergence rate of $n^{(1-\alpha_{ij})/2\alpha_{ij}}$.

\begin{remark}
To handle the microstructural noise, we use the sub-sampling scheme, and the number of non-overlapping quadratic pre-averaged variables $Q_{ij}(\tau_k)$ is $Cn^{1/2}$, which is known as the optimal choice.
That is, we are only able to  use   $n^{1/2}$ observations to estimate $T_{ij}$ due to the microstructural noise, which is the cost of managing the microstructural noise.
Thus, the optimal convergence rate is expected to be the square root of the rates of the estimators that are not affected by the microstructural noise.
From this point of view, the convergence rate $n^{(1-\alpha_{ij})/2\alpha_{ij}}$ is consistent with the results in \citet{devroye2016sub} and \citet{sun2019adaptive}.
\end{remark}

Recall that the ARP estimator has the following bias adjustment:
\begin{equation}\label{eq-4.4}
\hat{\Gamma}^{\alpha }_{ij} = \hat{T}^{\alpha}_{ij, \theta} - \hat{\rho}^{\alpha }_{ij, \theta}.
\end{equation}
Thus, to establish the concentration inequality for  the ARP estimator $\hat{\Gamma}^{\alpha}_{ij}$,  we need to investigate $\hat{\rho}^{\alpha}_{ij,\theta}$.
To do this, we use the  quadratic log-return random variables  $Q_{\rho,ij}(\tau_k)$ defined in \eqref{eq-3.7} and require the following moment condition. 
\begin{assumption} \label{Assumption2}
There exists a positive constant $\nu_{\rho,Q}$ such that
\begin{equation*}
\max_{1\leq i \leq p} \E \left\{  \left |Q^{c}_{\rho,ii}(\tau_{k})    \right  | ^{\alpha_{i}} \right\} \leq \nu_{\rho,Q}
\end{equation*}
for all $1\leq k \leq n-1$.
\end{assumption}

\begin{remark}
Assumption \ref{Assumption2} indicates that the continuous part of the observed log-return, $Y_i^{c} (\tau_{i,k+1})-Y_i^{c} (\tau_{i,k})$, has finite $2\alpha_{i}$-th moment.
We note that Assumption \ref{Assumption1}(c) is satisfied under Assumption \ref{Assumption1}(a),(b),(d) and Assumption \ref{Assumption2} (see Proposition \ref{Proposition4} in the Appendix).
\end{remark}
 
With this $\alpha_{ij}$-th moment condition, we establish the concentration inequalities for the ARP estimator $\hat{\Gamma}^{\alpha }_{ij}$ in the following theorem.

\begin{thm}\label{Theorem3}
(Upper bound)
Under the assumptions in Theorem \ref{Theorem1} and Assumption \ref{Assumption2}, take
\begin{equation*}
\theta_{\rho,ij} =    \left(\frac {\log \left(6\delta ^{-1}\right)}{\left( \alpha_{ij} -1\right) c_{\alpha_{ij} }S_{\rho,ij}(n-1)} \right)^{1/\alpha_{ij}},
\end{equation*}
where
$	S_{\rho,ij} = \dfrac {1}{n-1}\sum ^{n-1}_{k=1}U_{\rho,ij}\left( \tau_{k}\right)$, and $U_{\rho,ij}\left( \tau_{k}\right)$'s are some positive constants defined in Proposition \ref{Proposition3}(b). 
Then, for a sufficiently large $n$, we have
\begin{equation}\label{eq-4.5}
		 \Pr  \left \{   | \hat{\rho}^{\alpha}_{ij,\theta} - \rho_{ij} | \leq  C  \left( n^{-1}\log \delta ^{-1}\right) ^{(\alpha_{ij} -1)/\alpha_{ij} }
		 \right \}  \geq  1-\delta
\end{equation}
and
\begin{equation}\label{eq-4.6}
		 \Pr  \left \{   | \hat{\Gamma}^{\alpha}_{ij} - \Gamma_{ij} | \leq  C  \left( n^{-1/2}\log \delta ^{-1}\right) ^{(\alpha_{ij} -1)/\alpha_{ij} }
		 \right \}  \geq  1-2\delta.
\end{equation}
\end{thm}

Theorem \ref{Theorem3} shows that $\hat{\rho}^{\alpha}_{ij,\theta}$  has a sub-Weibull tail concentration bound with a convergence rate of $n^{(1-\alpha_{ij})/\alpha_{ij}}$, which is negligible compared to the upper bound in Theorem \ref{Theorem1}.
Thus, the ARP estimator has a sub-Weibull tail concentration with an optimal convergence rate of $n^{(1-\alpha_{ij})/2\alpha_{ij}}$  as shown in Theorems \ref{Theorem1}--\ref{Theorem2}.
Although the upper bound for $\hat{\rho}^{\alpha}_{ij,\theta}$ is dominated by the upper bound in Theorem \ref{Theorem1}, it is worth checking whether $\hat{\rho}^{\alpha}_{ij,\theta}$ is an optimal estimator.
Let $\hat{\rho}_{ij}\left( Q_{\rho,ij}\left(\tau_{k}\right) ,\delta \right) = \hat{\rho}_{ij}\left( Q_{\rho,ij}\left(\tau_{1}\right), \ldots, Q_{\rho,ij}\left(\tau_{n-1}\right) ,\delta \right)$ be any estimator for $\rho_{ij}$, possibly depending on $\delta$.
The following theorem provides a lower bound for the maximum concentration probability among the class of estimators $\hat{\rho}_{ij}\left( Q_{\rho,ij}\left(\tau_{k}\right) ,\delta \right)$ which satisfy $\max_{1\leq i \leq p} \E \left\{  \left |{Q}^{c}_{\rho,ii}(\tau_k)    \right  | ^{\alpha_{i}}   \right\}   \leq C$ for all $1\leq k \leq n-1$.

\begin{thm}\label{Theorem4}
(Lower bound)
Under the assumptions in Theorem \ref{Theorem3}, let $\alpha_{ij} \in \left(1,2\right)$ for some $1\leq i,j \leq p$.
Then, we have, for a sufficiently large $n$, 
\begin{equation}\label{eq-4.7}
		\min _{\hat{\rho}_{ij}\left( Q_{\rho, ij}\left( \tau_{k}\right) ,\delta \right)} \max_{\bQ_{\rho} ^c \in \Xi_{\rho}  (\alpha_1,\ldots, \alpha_p)}  \Pr  \left \{   | \hat{\rho}_{ij}\left( Q_{\rho,ij}\left( \tau_{k}\right) ,\delta \right) - \rho_{ij} | \geq  C  \left( \log \delta ^{-1}  /n\right) ^{(\alpha_{ij} -1)/\alpha_{ij} }
		 \right \}  \geq  \delta,
\end{equation}
where \begin{equation*}
\Xi_{\rho} (\alpha_1,\ldots, \alpha_p) =\{\bQ_{\rho}^c = (Q^{c}_{\rho, ii}(\tau_{k}))_{i=1,\ldots,p, k=1,\ldots, n-1} : \max_{i, k} \E \left\{  \left |Q^{c}_{\rho, ii}(\tau_{k})    \right  | ^{\alpha_{i}} \right\} \leq C \}.
\end{equation*}
\end{thm}

The upper bound in \eqref{eq-4.5} and  lower bound in \eqref{eq-4.7} match, which implies  that $\hat{\rho}^{\alpha}_{ij,\theta}$ achieves the optimal rate.
In sum, the proposed estimators for $T_{ij}$ and $\rho_{ij}$ are both optimal in terms of convergence rate, which implies that the ARP estimator is also optimal in the class of pre-averaging approaches.

\section{Application to large volatility matrix estimation}\label{SEC-5}

In this section, we discuss how to estimate large integrated volatility matrices based on the approximate factor model using  the ARP estimator.
Specifically, we assume that the integrated volatility matrix has the following low-rank plus sparse structure:
\begin{equation*}
 \bGamma  =   \bTheta  + \bSigma = \sum_{k=1}^r \bar{\lambda}_k \bar{\bq}_k \bar{\bq} _k ^\top + \bSigma,
\end{equation*}
where $\bTheta$ is a low-rank volatility matrix with the $i$-th largest eigenvalue $\bar{\lambda}_i>0$ and its associated eigenvector $\bar{\bq}_i$.
The low-rank volatility matrix $\bTheta$ accounts for the factor effect on the volatility matrix, and we assume that its rank $r$ is bounded.
The sparse volatility matrix $\bSigma$ stands for idiosyncratic risk and satisfies the following sparse condition:
\begin{equation}\label{eq-5.1}
\max_{ 1 \leq i \leq p} \sum_{j=1}^p | \Sigma_{ij}| ^{q}   (\Sigma_{ii} \Sigma_{jj} ) ^{(1-q)/2}  \leq M_{\sigma} s_p \text{ a.s.},
\end{equation}
where $M_\sigma$ is a positive random variable with $\E \( M_{\sigma}^2 \) < \infty$, $q \in [0,1)$, and $s_p$  is a deterministic function of $p$ that grows slowly in $p$.
When $\Sigma_{ii}$ is bounded from below, and $q=0$, $s_p$ measures the maximum number of nonvanishing elements in each row of matrix $\bSigma$.
This low-rank plus sparse structure is widely adopted for studying large matrix inferences \citep{ait2017using, bai2003inferential, bai2002determining,  fan2018robust,  fan2018eigenvector, jung2022next, kim2018Large,  stock2002forecasting}.

\subsection{Principal orthogonal complement thresholding}\label{SEC-5.1}
To harness the low-rank plus sparse structure, we employ the POET  method \citep{fan2013large} as follows.
 We first decompose an input volatility matrix using the ARP estimators in \eqref{eq-3.8} as follows:
  $$\hat{\bGamma} =\( \hat{\Gamma}^{\alpha }_{ij} \)_{i,j=1,\ldots, p}= \sum_{k=1}^p \hat{\lambda}_k \hat{ \bq}  _{k}  \hat{\bq}_k ^\top, $$
where $\hat{\lambda}_i$ is the $i$-th largest eigenvalue of $\hat{\bGamma}$, and $\hat{\bq}_i$ is its corresponding eigenvector.
 Then, using the first $r$ principal components, we estimate the low-rank volatility matrix $\bTheta$ by
  $$
  \hat{\bTheta}= \sum_{k=1}^r \hat{\lambda}_k \hat{\bq} _{k} \hat{\bq}_{k} ^\top.
  $$
  To estimate the sparse  volatility matrix $\bSigma$, we first calculate the input idiosyncratic volatility matrix  estimator  $\tilde{\bSigma} =(\tilde{\Sigma}_{ij})_{1\leq i,j\leq p} = \hat{\bGamma} - \hat{\bTheta}$ and employ the adapted thresholding method as follows:
\begin{equation*}
	\hat{\Sigma}_{ij} =
\begin{cases}
 \tilde{\Sigma}_{ij} \vee 0  , & \text{ if } i= j\\
s_{ij} ( \tilde{\Sigma}_{ij}) \1 ( |\tilde{\Sigma}_{ij}| \geq \varpi_{ij} ) ,  & \text{ if } i \neq j
\end{cases}
\quad \text{ and } \quad \hat{\bSigma} =(\hat{\Sigma}_{ij})_{1\leq i,j\leq p},
\end{equation*}
where the thresholding function $s_{ij} (\cdot) $ satisfies that $|s_{ij} (x) -x |\leq \varpi_{ij}$, and the adaptive thresholding level $\varpi_{ij} = \varpi_n\,\sqrt{ ( \tilde{\Sigma}_{ii} \vee 0 ) ( \tilde{\Sigma}_{jj} \vee 0 )}$, which corresponds to the correlation thresholding at level $\varpi_{n}$.
Examples of the thresholding function $s_{ij}(x)$ include the soft thresholding function $s_{ij}(x)=x-  \mbox{sign}(x) \varpi_{ij}$ and the hard thresholding function $s_{ij}(x)=x$.
The tuning parameter $\varpi_n$ will be specified in Proposition \ref{Proposition1}.
In the empirical study, we use the hard thresholding method.
We note that the large volatility matrix estimation method proposed by \citet{kong2018systematic} first estimates the instantaneous factor and idiosyncratic volatility matrices by applying principal component analysis (PCA) to the instantaneous volatility matrices.
Then, by aggregating the instantaneous factor and idiosyncratic volatility matrix estimators, they estimate the integrated factor and idiosyncratic volatility matrix, respectively. 
Similarly, we can use the proposed robust estimation procedure to estimate the instantaneous input volatility matrix. 
Then, we may obtain  similar results.

  With the low-rank volatility matrix estimator $\hat{\bTheta}=(\hat{\Theta}_{ij})_{1\leq i,j \leq p}$ and the sparse  volatility matrix estimator  $\hat{\bSigma}=(\hat{\Sigma}_{ij})_{1\leq i,j \leq p} $, we estimate the integrated volatility matrix $\bGamma$ by
\begin{equation*}
	\tilde{\bGamma}_{POET} = \hat{\bTheta} + \hat{\bSigma} .
\end{equation*}
 To investigate asymptotic behaviors of high-dimensional statistical inference methods such as the POET estimator, the sub-Weibull concentration inequality is  required and is satisfied by the ARP estimator, as shown in Theorem \ref{Theorem3}.
 Thus, the POET estimator based on the ARP estimators can enjoy  asymptotic properties similar to those established in \citet{fan2018robust}.
To study its asymptotic behaviors,  we  need the following   technical conditions, imposed by \citet{fan2018robust}, but the sub-Weibull concentration rate is different because we consider heterogeneous heavy-tailedness.

\begin{assumption}\label{Assumption3}
~
\begin{itemize}
\item [(a)] Let $D_\lambda  = \min \{\bar{ \lambda}_i -  \bar{\lambda}_{i+1} : 1 \leq i \leq r\}$, $  (\lambda _1 + p M_{\sigma} )  / D_{\lambda} \leq C_1$ a.s., and $D_{\lambda} \geq C_2  p $ a.s. for some generic constants $C_1$ and $C_2$,
where  $\bar{\lambda}_{r+1}=0$, $M_{\sigma}$ is defined in \eqref{eq-5.1}, and $\bar{\lambda}_i$ and $\lambda_i$ are the $i$-th eigenvalues of $\bTheta$ and $\bGamma$, respectively.

\item[(b)] For some fixed constant $C_3$,  we have
$$
 \frac{p}{r} \max_{ 1\leq  i \leq p} \sum_{j=1}^r \bar{q}_{ij} ^2  \leq C_3 \text{ a.s.},
  $$
where $\bar{\bq}_j = (\bar{q}_{1j}, \ldots, \bar{q}_{pj})^\top$ is the $j$-th eigenvector of $\bTheta$.

\item[(c)]The smallest eigenvalue  of $\bSigma$ stays away from zero almost surely.

\item[(d)]
$ s_p/ \sqrt{p}  +\left( n^{-1/2}\log p\right) ^{(\alpha_{\min} -1)/\alpha_{\min}} =o(1)$, where $\alpha_{\min}=\min_{1\leq i \leq p} \alpha_{i}$.
\end{itemize}
\end{assumption}

Under Assumption \ref{Assumption3},  we can establish the following proposition, similar to the proof of Theorem 3 in  \citet{fan2018robust}.
Below, we assume a generic input $\hat{\bGamma}$ that satisfies \eqref{eq-5.2}.
In particular, the ARP estimator satisfies the condition, as shown in Theorem \ref{Theorem3}.

\begin{proposition}\label{Proposition1}
Under the model  \eqref{eq-2.1}, let $\alpha_{\min}=\min_{1\leq i \leq p} \alpha_{i}$ and assume that the concentration inequality,
\begin{equation}\label{eq-5.2}
	\Pr \left \{ \max_{1\leq i, j \leq p} | \hat{\Gamma}_{ij} - \Gamma_{ij} | \geq  C  \left( n^{-1/2}\log p\right) ^{(\alpha_{\min} -1)/\alpha_{\min} }  \right \} \leq  p^{-1},
\end{equation}
Assumption \ref{Assumption3}, and the sparse condition \eqref{eq-5.1} are met.
Take $ \varpi_n=C_{\varpi} \beta_n$ for some large fixed constant $C_{\varpi}$, where $\beta_n =   M_\sigma    s_p/ p +  \left( n^{-1/2}\log p\right) ^{(\alpha_{\min} -1)/\alpha_{\min}}$.
Then, we have, for a sufficiently large n, with probability greater than $1-2p^{-1}$,
\begin{eqnarray}
	&& \| \hat{\bSigma}  - \bSigma  \|_{2} \leq C  M_{\sigma} s_p   \beta_n^{1-q}, \label{eq-5.3}\\
	&& \| \hat{\bSigma}  - \bSigma  \|_{\max} \leq C   \beta_n , \label{eq-5.4}\\
	&&\| \tilde{\bGamma}_{POET} - \bGamma \| _{\bGamma} \leq  C \[  p^{1/2} \left( n^{-1/2}\log p\right) ^{(2\alpha_{\min} -2)/\alpha_{\min} }   +   M_{\sigma} s_p  \beta_n^{1-q} \] ,  \label{eq-5.5} \quad \text{and} \\
	&& \| \tilde{\bGamma}_{POET} - \bGamma \| _{\max} \leq  C   \beta_n,  \label{eq-5.6}
\end{eqnarray}
where the relative Frobenius norm $\| \bA \| _{\bGamma} ^2 = p^{-1} \|\bGamma ^{-1/2}  \bA  \bGamma ^{-1/2}  \| _F^2.$
Furthermore, suppose that $M_{\sigma} s_p   \beta_n^{1-q} = o(1)$.
Then, with probability approaching 1, the minimum eigenvalue of $\hat{\bSigma}$ is bounded away from 0, $\tilde{\bGamma}_{POET}$ is non-singular,  
\begin{eqnarray}
	&& \| \hat{\bSigma}^{-1}  - \bSigma^{-1} \|_{2} \leq  C  M_{\sigma} s_p   \beta_n^{1-q}, \quad \text{and} \label{eq-5.7}\\
	&& \| \tilde{\bGamma}_{POET}^{-1}  - \bGamma^{-1} \|_{2} \leq  C  M_{\sigma} s_p   \beta_n^{1-q}. \label{eq-5.8}
\end{eqnarray}
\end{proposition}

\begin{remark}
Unlike Theorem 3 in \citet{fan2018robust}, Proposition \ref{Proposition1} imposes the sub-Weibull concentration condition \eqref{eq-5.2}, which is the optimal rate with only finite 2$\alpha_{\min}$-th moments, as shown in Theorems \ref{Theorem1}--\ref{Theorem4}.
Note that if $ 2p^{3} \in [ n^{c}, e^{n^{1/2}} ]$ for some positive constant $c >0$, Theorem \ref{Theorem3} shows that the ARP estimator satisfies \eqref{eq-5.2} for $\delta =1/\(2p^{3}\)$.
Also, the POET estimator is consistent in terms of the relative Frobenius norm as long as $p=o(n^{(2\alpha_{\min} -2)/\alpha_{\min}})$.
That is, the convergence rate is a function of the minimum tail index $\alpha_{\min}$.
\end{remark}

 In the numerical study, we often observe that $\tilde{\bGamma}_{POET}$  is singular when the sample size $n$ is small.
To overcome this issue, we add some small value to all of the diagonal entries of $\tilde{\bGamma}_{POET}$.
Specifically, we adjust $\tilde{\bGamma}_{POET}$ by
\begin{equation}\label{eq-5.9}
\hat{\bGamma}_{POET} = \tilde{\bGamma}_{POET} + c_{POET}n^{-1/2}\bI_{p}
\end{equation}
for some constant $c_{POET}>0$, where  $\bI_{p}$ is the $p$ dimensional identity matrix.
This type of adjustment is often used when estimating the large inverse matrices \citep{cai2011constrained, cai2016estimating}.
We note that the adjustment \eqref{eq-5.9} does not  affect the theoretical properties of the POET estimator.

\subsection{Tail index estimation}\label{SEC-5.2}
In this section, we propose an estimation procedure for the tail indices.
Specifically, we modify the Hill's estimator \citep{hill1975simple} as follows.
Let 
\begin{equation*}
H_i^{(j)} = j \log \(\dfrac{D_i^{(j)}}{D_i^{(j+1)}} \) \quad \text{for} \quad j=1, \ldots, n-2,
\end{equation*}
where $D_{i}^{(1)} \geq \ldots \geq D_{i}^{(n-1)}$ are the order statistics of  $\{ | Y_{i}\left(\tau_{i,2}\right)-Y_{i}\left(\tau_{i,1}\right)|, \ldots, | Y_{i}\left(\tau_{i,n}\right)-Y_{i}\left(\tau_{i,n-1}\right)| \}$.
Then, the tail index is estimated by
\begin{equation}\label{eq-5.10}
\tilde{\alpha}_{i} = \dfrac{u_n}{2 \sum_{j=1}^{u_n} H_i^{(j)} \1 \( H_i^{(j)} \leq \omega \)},
\end{equation}
where $u_n\to \infty$ is a given sequence such that $u_n/n \rightarrow 0$  and  $\omega$ is a truncation parameter which will be determined in Theorem \ref{Theorem5}. 
We note that  the indicator function is used to handle the jumps in $D_i^{(j)}$.

To investigate the theoretical properties of $\tilde{\alpha}_{i}$, we need the following assumptions.

\begin{assumption}\label{Assumption4}
~
\begin{itemize}
\item [(a)] For each $1 \leq i \leq p$, $| \epsilon_{i}\left(\tau_{i,2}\right)-\epsilon_{i}\left(\tau_{i,1}\right)|, \ldots, | \epsilon_{i}\left(\tau_{i,n}\right)-\epsilon_{i}\left(\tau_{i,n-1}\right)|$ have the cumulative distribution function $F_i(x)$ satisfying
\begin{equation*}
F_i(x) = 1-(c_i/x)^{2\alpha_i} \quad \text{for} \,\,\, x \geq C \(n/u_n\)^{1/(2\alpha_i)},
\end{equation*}
where $c_i$ is some positive constant.

\item[(b)] $\max\limits_{1 \leq i \leq p}{\alpha_i} \leq C$, $\max\limits_{1 \leq i \leq p}{c_i} \leq C$, and the Blumenthal-Gettor index $\pi=0$.

\item[(c)] $n^{\kappa}  \leq p$ for some positive constant  $\kappa$, $\(\log p\)^2 / u_n  \rightarrow 0$ as $n,p \rightarrow \infty$, and $u_n \leq n^{(\alpha_{\max} +1)/(2\alpha_{\max} +1)}$, where $\alpha_{\max}=\max\limits_{1 \leq i \leq p}\alpha_i$. 
\end{itemize}
\end{assumption}

\begin{remark}
Assumption \ref{Assumption4}(a) indicates that $| \epsilon_{i}\left(\tau_{i,k}\right)-\epsilon_{i}\left(\tau_{i,k-1}\right)|$, $k=2, \ldots, n$,  have the same tail distribution  as the Pareto distribution with a scale parameter $c_i$ and a shape parameter $2\alpha_i$, which has only finite $2\alpha_i$-th moment.
The Pareto distribution is widely used to model the heavy-tailed distribution \citep{chin2008heavy, coronel2005fitting, dagum2014income, nirei2016pareto}. 
Since we only require the condition for the tail, Assumption \ref{Assumption4}(a) is not restrictive. 
We note that for appropriate scale parameter $c_i$, the Pareto distribution with the shape parameter $2\alpha_i$ has asymptotically the same tails as the t-distribution with $2\alpha_i$ degrees of freedom.
Finally, we need the condition $u_n \leq n^{(\alpha_{\max} +1)/(2\alpha_{\max} +1)}$ to handle the jumps using the threshold method.
Since $(\alpha_{\max} +1)/(2\alpha_{\max} +1)>1/2$, we can choose $u_n= C_u n^{1/2}$, which results in the convergence rate $n^{-1/4}$ (see Theorem \ref{Theorem5}).
Thus, this condition is not restrictive.
\end{remark}

The theorem below shows that $\tilde{\alpha}_{i}$ has the sub-exponential tail concentration with a convergence rate of $u_n^{-1/2}$.

\begin{thm}\label{Theorem5}
Under the models \eqref{eq-2.1} and \eqref{eq-2.3} and Assumption \ref{Assumption4}, for any given positive constant $a$, choose $\omega = C_{\omega,a} \log p$ for some large constant $C_{\omega,a}>0$.
Then, for a sufficiently large $n$, we have, with probability at least $1-p^{-a}$, 
\begin{equation}\label{eq-5.11}
 \max_{1 \leq i \leq p} | \tilde{\alpha}_{i} - \alpha_i | \leq  C  u_n^{-1/2}\log p.
\end{equation}
\end{thm}

As discussed above, since $(\alpha_{\max} +1)/(2\alpha_{\max} +1)>1/2$, we can obtain the convergence rate \eqref{eq-5.11} of at least $n^{-1/4}$.
In this paper, we choose $u_n = C_u n^{1/2}$ for some constant $C_u$.
This choice is enough to obtain the theoretical results obtained in Theorem \ref{Theorem3}, which will be shown in Proposition \ref{Proposition2}.

Since $\tilde{\alpha}_{i}$ has the estimation error, $\tilde{\alpha}_{i}$ can be bigger than the true tail index $\alpha_i$.
In this case,  the truncation parameters $\theta_{ii}$ and $\theta_{\rho,ii}$ may go to zero and then the heavy-tailedness may not be handled well.
To tackle this obstacle, we adjust $\tilde{\alpha}_{i}$ as follows:
\begin{equation}\label{eq-5.12}
    \hat{\alpha}_i  = \tilde{\alpha}_{i}  - c_{\xi}n^{-\xi} \log p,
\end{equation}
where $c_{\xi}$ is some  positive constant and $0<\xi<1/4$.
Then, under the assumptions in Theorem \ref{Theorem5}, we have
\begin{equation}\label{eq-5.13}
  \Pr \left\{ \max_{1 \leq i \leq p} (\hat{\alpha}_i - \alpha_i ) <0 \quad \text{and} \quad  \max_{1 \leq i \leq p} | \hat{\alpha}_i - \alpha_i | \leq  C  n^{-\xi}\log p \right\} \geq 1-p^{-a}.
\end{equation}

The following proposition shows that with $\hat{\alpha}_i$, we can  obtain the same theoretical results as in Theorem \ref{Theorem3}.

\begin{proposition}\label{Proposition2}
Under the assumptions in Theorem \ref{Theorem1}, Theorem \ref{Theorem3}, and Theorem \ref{Theorem5}, let $\log p \leq Cn^{w}$ for some positive constant $w < \xi$. 
Choose the estimator of $\alpha_{ij}$ as
\begin{equation}\label{eq-5.14}
	\hat{\alpha}_{ij} = 2\wedge \dfrac {2\hat{\alpha}_{i}\hat{\alpha}_{j}}{\hat{\alpha}_{i}+\hat{\alpha}_{j}}.
\end{equation}
Then, for any given positive constant $a$, we have, with probability at least $1-p^{-a}$, 
\begin{equation}\label{eq-5.15}
 | \hat{\Gamma}^{\alpha}_{ij} - \Gamma_{ij} | \leq  C  \left( n^{-1/2}\log p \right) ^{(\alpha_{ij} -1)/\alpha_{ij}}.
\end{equation}
\end{proposition}

Proposition \ref{Proposition2} indicates that the ARP estimator with $\hat{\alpha}_{ij}$ also has the optimal convergence rate.
Thus, the POET estimator based on the ARP estimators also satisfies \eqref{eq-5.2}, and so it can enjoy the same theoretical results as in  Proposition \ref{Proposition1}.

\subsection{Discussion on the  tuning parameter selection}\label{SEC-5.3}

To implement the ARP estimation procedure, we need to choose tuning parameters.
In this section,  we discuss how to select the tuning parameters  for the numerical studies.
To obtain $\tilde{\alpha}_{i}$ in \eqref{eq-5.10}, we choose
\begin{equation*}
	u_n =c_u \lfloor n^{1/2} \rfloor \quad \text{and} \quad \omega = c_{\omega} \log p, 
\end{equation*}
where $c_u$ and $c_{\omega}$ are tuning parameters. 
Since the tail index should be bigger than 1, we adjust $\tilde{\alpha}_{i}$ as follows:
\begin{equation}\label{eq-5.16}
    \hat{\alpha}_i  = \max\{\tilde{\alpha}_{i}  - c_{\xi}n^{-\xi} \log p, 1.1 \},
\end{equation}
where $\xi$ and $c_{\xi}$ are tuning parameters. 
In the numerical study, we choose $c_u =2$, $c_{\omega}=1/3$, $\xi=0.2$, and $c_{\xi}=0.01$.
Then, we obtained  $\hat{\alpha}_{ij}$ based on \eqref{eq-5.14}.
For the estimation of $T_{ij}$ and $\rho_{ij}$, we note that the truncation function $\psi_{\alpha}(x)$ truncates the variables around zero, which produces some bias \citep{minsker2018sub}.
We also note that truncating the variables around some other constant does not affect the asymptotic theoretical results.
Thus, to improve the performance of the ARP estimator in the numerical perspective, we truncate the variables around their median value.
Specifically, let $MQ_{ij}$ be the median of $Q_{ij}\left( \tau_{k}\right)$, $k=1, \ldots, n-K_n$, and $MQ_{\rho,ij}$ be the median of $Q_{\rho,ij}\left( \tau_{k}\right)$, $k=1, \ldots, n-1$.
Then, we estimate $T_{ij}$ and $\rho_{ij}$ as follows: 
\begin{equation*}
	\hat{T}^{\alpha }_{ij, \theta} = MQ_{ij} +  \dfrac {1}{(n-K_{n})\theta_{ij}}\sum_{k=1}^{n-K_n}\psi_{\alpha_{ij}}\left\{\theta _{ij}\[Q_{ij}\left(\tau_{k}\right)-MQ_{ij}\]\right\}
\end{equation*}
and
\begin{equation*}
	\hat{\rho}^{\alpha }_{ij, \theta} =  \dfrac {(n-1)\zeta}{\phi K_n}\Big[MQ_{\rho,ij} + \dfrac{1}{(n-1)\theta_{\rho,ij}}\sum_{k=1}^{n-1}\psi_{\alpha_{ij}}\left\{\theta _{\rho,ij}\[Q_{\rho,ij}\left(\tau_{k}\right)-MQ_{\rho,ij}\]\right\}\Big].
\end{equation*} 
With this scheme, we choose the thresholding level as follows:
\begin{equation}\label{eq-5.17}
\theta_{ij} =    c\left(\frac {K_{n}\log p}{\left( \hat{\alpha}_{ij} -1\right) c_{\hat{\alpha}_{ij} }\hat{S}_{ij}(n-K_{n})} \right)^{1/\hat{\alpha}_{ij}}
\end{equation}
and
\begin{equation}\label{eq-5.18}
\theta_{\rho,ij} =    c\left(\frac {\log p}{\left( \hat{\alpha}_{ij} -1\right) c_{\hat{\alpha}_{ij} }\hat{S}_{\rho,ij}(n-1)} \right)^{1/\hat{\alpha}_{ij}},
\end{equation}
where $\hat{S}_{ij}=\dfrac{1}{n-K_n}\sum_{k=1}^{n-K_n}\left\{ \left| Q_{ij}\left( \tau_{k}\right) -MQ_{ij} \right| ^{\hat{\alpha}_{ij}}\right\}$, $\hat{S}_{\rho,ij}=\dfrac{1}{n-1}\sum_{k=1}^{n-1}\left\{ \left| Q_{\rho,ij}\left( \tau_{k}\right) - MQ_{\rho,ij} \right| ^{\hat{\alpha}_{ij}}\right\}$, and $c$ is a tuning parameter.
In the simulation study, we choose $c$ as 0.5.
For the empirical study, we choose $c$ that minimizes the corresponding mean squared prediction error (MSPE) for the in-sample period. 
Details can be found in Section \ref{SEC-6.2}.
In $\psi_{\alpha}(x)$, $c_{\alpha}$ is determined by $\hat{\alpha}_{ij}$, that is, $c_{\hat{\alpha}_{ij}}= \max \left \{(\hat{\alpha}_{ij} -1)/\hat{\alpha}_{ij} , \sqrt {(2-\hat{\alpha}_{ij})/\hat{\alpha}_{ij} }\right \}$.
In the pre-averaging stage, we choose $K_n= \lfloor n^{1/2} \rfloor$ and $g(x) = x  \wedge (1 - x)$.
Finally, we select $ c_{POET} = 3  \text{median}\left\{\lambda_1\( \tilde{\bGamma}_{POET}\), \ldots, \lambda_p\( \tilde{\bGamma}_{POET}\) \right\}$, where $\lambda_i\( \tilde{\bGamma}_{POET}\)$  is the $i$-th largest eigenvalue of $\tilde{\bGamma}_{POET}$.

\section{Numerical study}\label{SEC-6}
\subsection{A simulation study}\label{SEC-6.1}
To check the finite sample performance of the ARP estimator, we conducted a simulation study.
We considered the jump diffusion process and generated the data with frequency $1/n^{all}$.
We used the heterogeneous heavy-tail process (heavy-tail process 1), a homogeneous heavy-tail process (heavy-tail process 2), and a sub-Gaussian process. 
We generated the non-synchronized observation time points and employed the refresh time scheme.
Specifically, we considered  the following true log-price process:
\begin{equation*}
 d \bX(t)  =   \bmu (t)dt +  \bvartheta ^{\top} (t) d\bW_t^{*}  +  \bsigma ^\top(t)d\bW_t+\bL(t)d\bN(t),
\end{equation*}
where   $\bmu(t) = (0.02, \ldots, 0.02)^\top$, $\bW_t^{*}$ and  $\bW_t$ are $r$  and $p$ dimensional independent Brownian motions, respectively, $\bvartheta (t)$ and $ \bsigma (t)$ are $r$ by $r$ and $p$ by $p$ matrices, respectively,  $\bL(t)$ is the jump size, and $\bN(t)$ is the $p$ dimensional Poisson process with the intensity $\bI(t)$.
To generate two heavy-tail processes, we used a setting similar to those in \citet{wang2010vast} and \citet{fan2018robust}.
Specifically, let $\bsigma(t)$ be the Cholesky decomposition of $\bvarsigma (t)= (\varsigma_{ij} (t))_{1\leq i,j \leq p} $.
 The diagonal elements of $\bvarsigma(t)$ come from four different processes: geometric Ornstein-Uhlenbeck processes,  the sum of two CIR processes  \citep{cox1985theory, barndorff2002econometric}, the volatility process in Nelson's GARCH diffusion limit model \citep{wang2002asymptotic}, and the two-factor log-linear stochastic volatility process \citep{huang2005relative} with leverage effect.
Details can be found in \citet{wang2010vast}.
To control the tail behaviors of the instantaneous volatility matrix $\bvarsigma(t)$, we used the  $t$-distribution as follows:
\begin{equation*}
	\varsigma_{ii} (t_l)= \(1+  \left |t_{df_{i},l}\right | \) \, \tilde{\varsigma}_{ii}(t_l),
\end{equation*}
where for $l=1,\ldots, n^{all}$, $t_{df_{i},l}$ are the i.i.d. $t$-distributions with degrees of freedom $df_{i}$, $t_l=l/n^{all}$, and  $\tilde{\varsigma}_{ii}(t_l)$ were generated by the four processes listed above. 
 To account for the heterogeneous heavy-tailed distribution (heavy-tail process 1), $df_{i}$ were generated from the unif(3, 4), whereas, for the homogeneous heavy-tailed distribution (heavy-tail process 2), we set $df_{i}=5$.
To obtain the sparse instantaneous volatility matrix $\bvarsigma(t)$, we generated its off-diagonal elements as follows:
\begin{equation*}
  \varsigma_{ij}(t_l) = \left\{ \kappa(t_l)\right\}^{|i - j|} \sqrt{\varsigma_{ii}(t_l)\varsigma_{jj}(t_l)}, \quad 1 \leq i \neq j \leq p,
\end{equation*}
 where the process $\kappa(t)$ is given by
\begin{eqnarray*}
  &&\kappa(t) = \frac{e^{\frac{1}{2} u(t)} - 1}{e^{\frac{1}{2} u(t)} + 1}, \hspace{0.5cm} du(t) = 0.03 \{ 0.64 - u(t)\}dt + 0.118u(t)dW_{\kappa, t},\cr
  &&W_{\kappa, t} = \sqrt{0.96}W_{\kappa, t}^0 - 0.2 \sum\limits_{i = 1}^p W_{it}/\sqrt{p},
\end{eqnarray*}
and   $W_{\kappa, t}^0, \kappa=1,\ldots, p$, are  one dimensional Brownian motions independent of the Brownian motions $\bW_t ^\ast$ and $\bW_t$.
The low-rank instantaneous volatility matrix  $\bvartheta^\top(t)\bvartheta(t)$ is $\bB ^\top \{\bvartheta^f (t) \}^\top  \bvartheta^f (t) \bB$, where $\bB= (B_{ij}) _{1 \leq i \leq r, 1\leq j \leq p} \in \mathbb{R}^{r \times p}$, and $B_{ij}$ was generated from the unif(-0.7, 0.7), and $\bvartheta^f (t)$  was generated similar to $\bsigma(t)$.
Specifically, $\bvartheta^f (t)$ is the Cholesky decomposition of $\bvarsigma ^f (t)$, and the diagonal elements of $\bvarsigma ^f (t)$ at time $t_l$ were
\begin{equation*}
	\varsigma_{ii}^f (t_l)= \left \{ 1+   \left |t_{df_{i},l}^f \right | \right  \} \, \tilde{\varsigma}_{ii} ^f(t_l),
\end{equation*}
where $t_{df_{i},l}^f, l=1,\ldots, n^{all}$, are the i.i.d. $t$-distributions with degrees of freedom $df_{i}$, and $\tilde{\varsigma}_{ii}^f(t_l), l=1,\ldots, n^{all}$, were generated from the geometric Ornstein-Uhlenbeck processes.
The off-diagonal elements of $\bvarsigma ^f (t)$ were set as zero.
For the jump part, we chose $\bI(t) = (25, \ldots, 25)^\top$, and the jump size $L_i(t)$  was obtained from independent $t$-distribution with degrees of freedom $df_{i}$ and standard deviation $0.15 \sqrt{ \int _{0}^1  \gamma _{ii}(t)  dt }$.
We also generated a sub-Gaussian process similarly to the heavy-tail process, except that $t$-distribution terms were set as standard normal distribution terms.

To generate the observation time points, we first obtained $(n^{all}+1)$ sampling time points, $t_k=k/n^{all}$, $k=0, \ldots, n^{all}$. 
Based on these points, we generated non-synchronized data similar to the scheme in \citet{ait2010high}, as follows.
 First, $p$ random proportions $w_{i}, i=1, \ldots, p$, were independently generated from the unif(0.8, 1).
Second, we set each $t_k$ as the observation time point of the $i$-th asset if the independent Bernoulli random variable with parameter $w_i$ had a value of 1.
Third, the noise-contaminated high-frequency observations $Y_i(t_{i,k})$ were generated from the model \eqref{eq-2.3}.
Specifically, the noise $\epsilon_{i} (t_{i,k})$ was obtained from independent $t$-distribution with degrees of freedom $df_{i}$ and standard deviation $0.05 \sqrt{ \int _{0}^1  \gamma _{ii}(t)  dt }.$
We chose $p=200$ and $r=3$, and we varied $n^{all}$ from $1000$ to $4000$. We employed the refresh time scheme to obtain synchronized data.

To investigate the effect of the adaptiveness of the proposed ARP procedure, we introduce a universal robust pre-averaging realized volatility (URP) estimator, which uses the same estimation procedure as the ARP estimator with $\hat{\alpha}_{ij}=2$ for all $1\leq i,j \leq p$.
That is, the URP estimator truncates the pre-averaged variables with the universal tail index level.
Also, we employed the robust pre-averaging realized volatility (RPRV) estimator \citep{fan2018robust}, which can handle the bounded fourth moment condition as the URP estimator.
Specifically, we first obtained
\begin{equation*}
\hat{T}_{ij}^{\text{RPRV}}= \arg \min \limits_{x} \sum_{k=0}^{n-K_n} \ell_{h_{ij}}\left(\dfrac{n-K_n +1}{\phi K_n}Z_i(\tau_k) Z_j(\tau_k) \mathbf{1}\left\{ \left|Z_i(\tau_k) \right| \leq v_{i,n}\right\} \mathbf{1}\left\{ \left| Z_j(\tau_k) \right| \leq v_{j,n}\right\} - x \right),
\end{equation*}
where $h_{ij}$ is the truncation parameter, $K_n = \lfloor n^{1/2} \rfloor$, 
$\ell_{h}$ is the Huber loss
\begin{equation*}
\ell_ h  (x)= \begin{cases}
  2 h |x| -h^2 & \text{ if }  |x| \geq h \\ 
 x^2 & \text{ if }  |x| < h,
\end{cases}
\end{equation*}
and $v_{i,n} = c_{i,v}n^{-0.235}$ is a thresholding level for some constant $c_{i,v}$.
We chose $c_{i,v}$ as 7 times the sample standard deviation for the pre-averaged variables $n^{1/4}Z_{i}\left(\tau_k\right)$.
We note that the thresholding for $Z_i(\tau_k)$ is used to handle the jumps.
Then, we obtained 
\begin{equation*}
\hat{\eta}_{ii}^{\text{RPRV}}= \arg \min \limits_{x} \sum_{k=0}^{n-1} \ell_{h_{\eta,ii}}\left(\[Y_i(\tau_{k+1}) -Y_i(\tau_k) \]^2/2 - x \right),
\end{equation*}
where $h_{\eta,ii}$ is the truncation parameter.
With $\hat{T}_{ij}^{\text{RPRV}}$ and $\hat{\eta}_{ii}^{\text{RPRV}}$, we calculated the RPRV estimator as follows:
\begin{equation*}
\hat{\Gamma}_{ij}^{\text{RPRV}}= \hat{T}_{ij}^{\text{RPRV}} - \dfrac{n-K_n +1}{\phi K_n} \zeta \hat{\eta}_{ij}^{\text{RPRV}} \mathbf{1}(i=j).
\end{equation*}
In the numerical study, we chose 
\begin{equation*}
    h_{ij}= \sqrt{\dfrac{(n-K_n+1) \hat{b}_{ij}}{K_n \log p}}  \quad \text{and} \quad  h_{\eta,ii} = \sqrt{\dfrac{n \hat{b}_{\eta,ii}}{2 \log p}},
\end{equation*}
where the asymptotic variance estimators $\hat{b}_{ij}$ and $\hat{b}_{\eta,ii}$ can be obtained by (4.9) and (5.2) in  \citet{fan2018robust}.
We also employed the jump adjusted pre-averaging realized volatility matrix (PRVM) estimator \citep{ait2016increased, christensen2010pre, jacod2009microstructure} as follows:
\begin{equation*}
	\hat{\Gamma}_{ij}^{\text{PRV}}=\frac{1}{\phi K_n} \sum^{n-K_n}_{k=1}\left\{Z_i\left(\tau_{k}\right)Z_j\left(\tau_{k}\right)-\frac{1}{2}\hat{Y}_{i,j}\left(\tau_{k}\right) \right\}\mathbf{1}\left\{ \left| Z_{i}\(\tau_{k}\) \right| \leq v_{i,n}\right\}\mathbf{1}\left\{ \left| Z_{j}\(\tau_{k}\) \right| \leq v_{j,n}\right\},
\end{equation*}
where
\begin{eqnarray*}
	&&\hat{Y}_{i,j}\left(\tau_{k}\right) =\sum^{K_n}_{l=1}\Bigg[\left\{g\left(\frac{l}{K_n}\right) - g\left(\frac{l-1}{K_n}\right) \right\}^2 \cr
	&& \qquad \qquad \qquad  \quad \times  \left(Y_i(\tau_{k+l})-Y_i(\tau_{k+l-1})\right)\left(Y_j(\tau_{k+l})-Y_j(\tau_{k+l-1})\right)\Bigg],
\end{eqnarray*}
$v_{i,n}$ is defined the same as in the case of the RPRV estimator, $K_n = \lfloor n^{1/2} \rfloor$, and $g\left(x\right) = x \wedge \left(1 - x\right)$.
Thus, we calculated the input volatility matrix using the adaptive robust pre-averaging realized volatility matrix (ARPM), universal robust pre-averaging realized volatility matrix (URPM),  robust pre-averaging realized volatility matrix (RPRVM),  and jump adjusted pre-averaging realized volatility matrix (PRVM) estimators.
We used the tuning parameters discussed in Section \ref{SEC-5.3}.
We note that the PRVM estimator cannot account for the heavy-tail and that the URPM  and  RPRVM estimators cannot explain the heterogeneity of different degrees of the heaviness of tail distributions.

To make the estimates positive semi-definite, we projected the volatility matrix estimators onto the positive semi-definite cone in the spectral norm.
To calculate the POET estimators, we used the hard thresholding scheme  and  selected  the thresholding level  by minimizing the corresponding Frobenius norm.
Then, we adjusted the POET estimators using the identity matrix to avoid the singularity for the finite sample as described in \eqref{eq-5.9}.
The average estimation errors under the Frobenius norm, relative Frobenius norm, $\| \cdot \|_{\bGamma}$, $\ell_2$-norm (spectral norm), and maximum norm  were computed based on 1000 simulations. 
The average numbers of synchronized time points with the refresh time scheme were   300.5, 599.8, 1199.7 for $n^{all}=1000, 2000, 4000$, respectively.

\begin{table}[!ht]
\caption{The mean squared errors (MSEs) of  estimators for  $\alpha_{ij}$ given $n^{all}=1000, 2000, 4000$.}\label{Table-1}
\centering
\scalebox{0.91}{
\begin{tabular}{l l c c c  c c c c c c c}
\hline
\multicolumn{2} {c}{} & \multicolumn{5} {c}{MSE}   \\ \cline{3-7}
 Tail type $\backslash$ $n^{all}$	&&		1000 	&&	2000	&&	4000 \\ \hline
 Heterogeneous 	&&			0.056 	&& 0.039 	&&	0.028			\\
				
  Homogeneous 	&& 0.119 	&&	0.054	&&	0.019	 \\\hline
\end{tabular}
}
 \end{table}

Table \ref{Table-1} reports the mean squared errors (MSEs) of estimators for $\alpha_{ij}$ against the sample size $n^{all}$ for two heavy-tail processes.
For the heterogeneous heavy-tail process, 2$\alpha_i$'s were generated from the  unif(3, 4),  and 2$\alpha_i$=5 for the homogeneous heavy-tail process.
We calculated $\alpha_{ij}$ using \eqref{eq-5.14} in Section \ref{SEC-5.2}.
From Table \ref{Table-1}, we can find that for the heterogeneous and homogeneous heavy-tail processes, the MSE decreases as the sample size $n^{all}$ increases.
We also find that the MSEs for the heterogenous case are much smaller than the
homogeneous case when $n^{all}=  1000, 2000$ but larger when $n^{all} =4000$.
This is because for the homogeneous case, the true $\alpha$ is 2.5, while, according to the proposed procedure and RPRVM, the maximum possible $\alpha$ is 2. 
 Thus, we calculated the MSEs with $\alpha=2$. 
 When $n^{all}$ is small, there were relatively many estimates smaller than 2.
 In contrast, when $n^{all}$ is large ($n^{all}=4000$),  most of estimates were greater than 2.
 Finally, we note that, for the sub-Gaussian process, more than 99 percent of $\alpha_{ij}$ was estimated to be 2 for $n^{all}=4000$  (this is regarded as correctly estimated, due to the sub-Gaussianity of the truncated average).
These results indicate that  the proposed tail index estimator works well.

\begin{figure}[!ht]
\centering
\includegraphics[width = 0.75\textwidth]{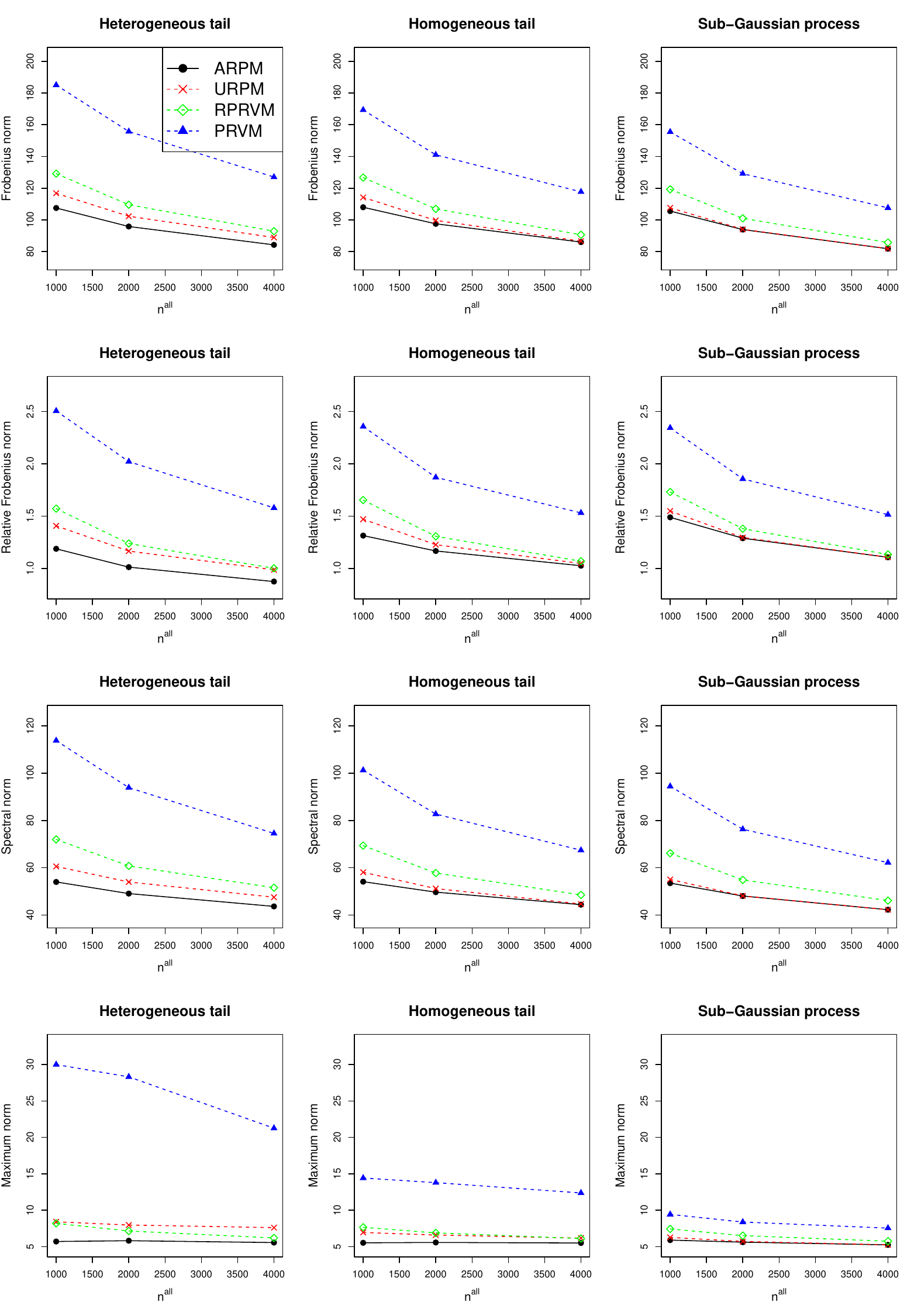}
\caption{The Frobenius, relative Frobenius, spectral, and max norm error plots (corresponding to four rows) of the POET estimators with the ARPM (black dot), URPM (red cross), RPRVM (green diamond), and PRVM (blue triangle) estimators for $p=200$ and $n^{all}=1000, 2000, 4000$.}\label{Fig-2}
\end{figure}

\begin{figure}[!ht]
\centering
\includegraphics[width = 0.78\textwidth]{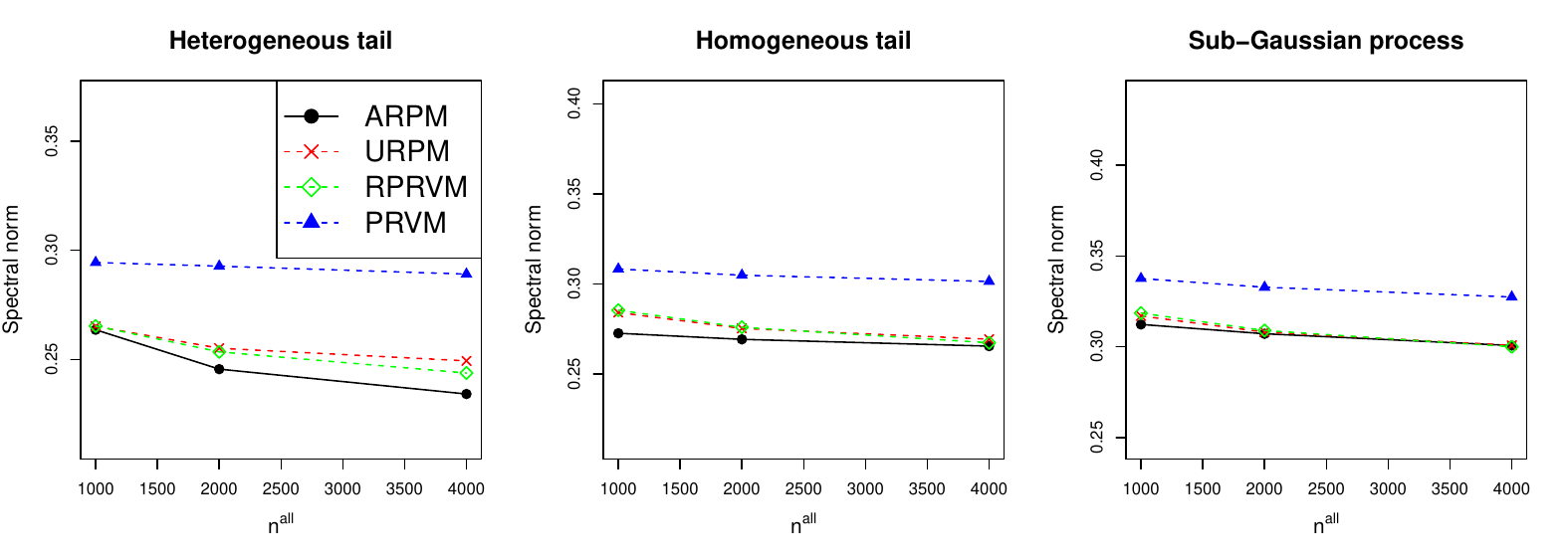}
\caption{The spectral norm error plots of the inverse POET estimators with ARPM (black dot), URPM (red cross), RPRVM (green diamond), and PRVM (blue triangle) estimators for $p=200$ and $n^{all}=1000, 2000, 4000$.}\label{Fig-3}
\end{figure}

Figure \ref{Fig-2} plots the Frobenius, relative Frobenius, spectral, and max norm errors against the sample size $n^{all}$ for the POET estimators from the ARPM, URPM, RPRVM, and PRVM estimators. 
Figure \ref{Fig-3} depicts the spectral norm errors against the sample size $n^{all}$ for the inverse POET estimators with the ARPM, URPM, RPRVM, and PRVM estimators.
As expected, the ARPM estimator outperforms other estimators for the heterogeneous heavy-tail process.
For the homogeneous heavy-tail and sub-Gaussian processes, the ARPM, URPM, and RPRVM estimators perform similarly and outperform the PRVM estimator.
One possible explanation of the poor performance of the PRVM estimator in the Gaussian noise case is that the true return process contains  heavy distributions over time; hence, robust methods outperform.
In sum, the ARPM estimator is robust to the heterogeneity of the heaviness of tails and adapts to the homogeneity of the heaviness of tails.

\subsection{An empirical study}\label{SEC-6.2}

\begin{figure}[!ht]
\centering
\includegraphics[width = 0.8 \textwidth]{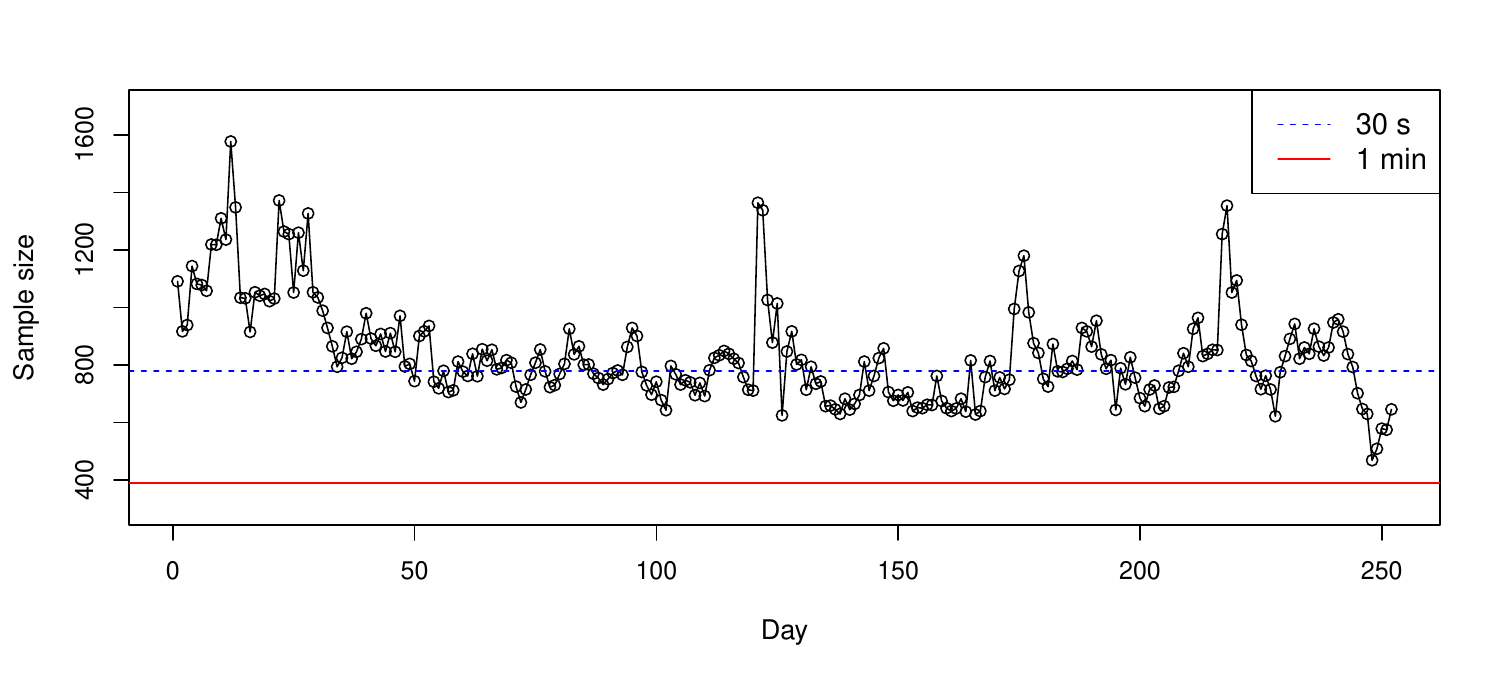}
\caption{The number of daily synchronized samples from the refresh time scheme for 200 assets over 251 days in 2016.
The  blue dash  and red solid  lines indicate the numbers of possible observations for 30-sec and 1-min log-returns in each trading day, which are 780  and 390, respectively.}\label{Fig-4}
\end{figure}

In this section, we applied the proposed ARP estimator to the high-frequency trading data of 200 assets, collected from January 2015 to December 2016 (501 trading days).
The 200  largest-volume stocks were selected from the S\&P 500, and the data was obtained from the Wharton Data Service (WRDS) system.
Days with half trading hours were excluded.
Figure \ref{Fig-4} plots the daily synchronized sample sizes from the refresh time scheme for the 200 assets.
As seen in Figure \ref{Fig-4}, sampling frequency higher than 1 minute can lead to the nonexistence of the observation between some consecutive sample points.
Hence, we employed 1-min log-return data with the previous tick scheme to mitigate the potential heterogeneity from observed time intervals and the irregular observation time errors  \citep{mykland2019algebra}.

\begin{figure}[!ht]
\centering
\includegraphics[width = 0.64\textwidth]{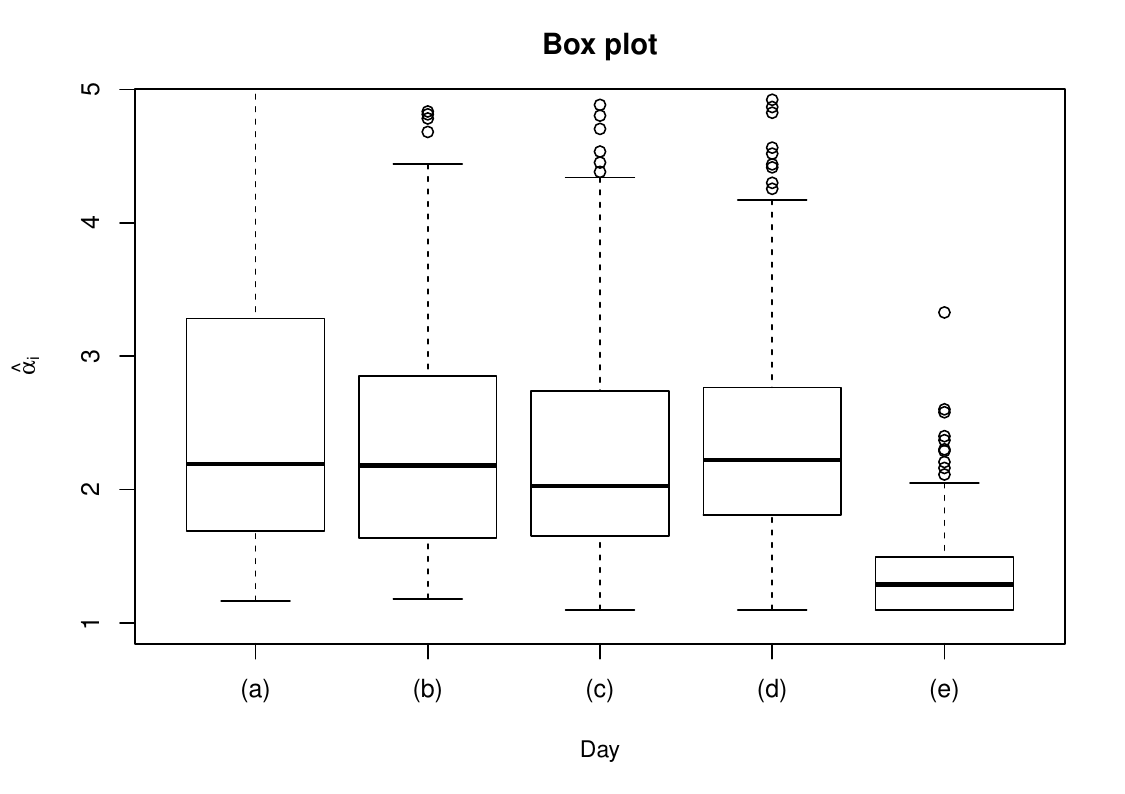}
\caption{The box plots of the distributions of the daily estimated tail indices $\hat{\alpha}_i$ for the 200 most liquid stocks among the S\&P 500 index in 2015--2016.
Day (a) has the largest IQR, and days (b)--(e) have the 75th, 50th, 25th, and 0th (minimum) percentile of the IQR among 501 trading days in 2015--2016, respectively. }\label{Fig-5}
\end{figure}

To measure the heterogeneous heavy-tailedness over time, we estimated the tail indices using $\hat{\alpha}_i$ proposed in \eqref{eq-5.16}.  
Figure \ref{Fig-5} shows the box plots of daily estimated tail indices $\hat{\alpha}_{i}$  of 200 assets for each of five selected days in 2015--2016: from the day with the largest IQR to the day with the smallest IQR among 501 days.
It provides stark evidence that the tail indices of observed log-returns are heterogeneous over time, which matches the daily kurtoses result in Figure \ref{Fig-1}.
This supports the heterogeneous heavy-tail assumption.

\begin{figure}[!ht]
\centering
\includegraphics[width = 0.8\textwidth]{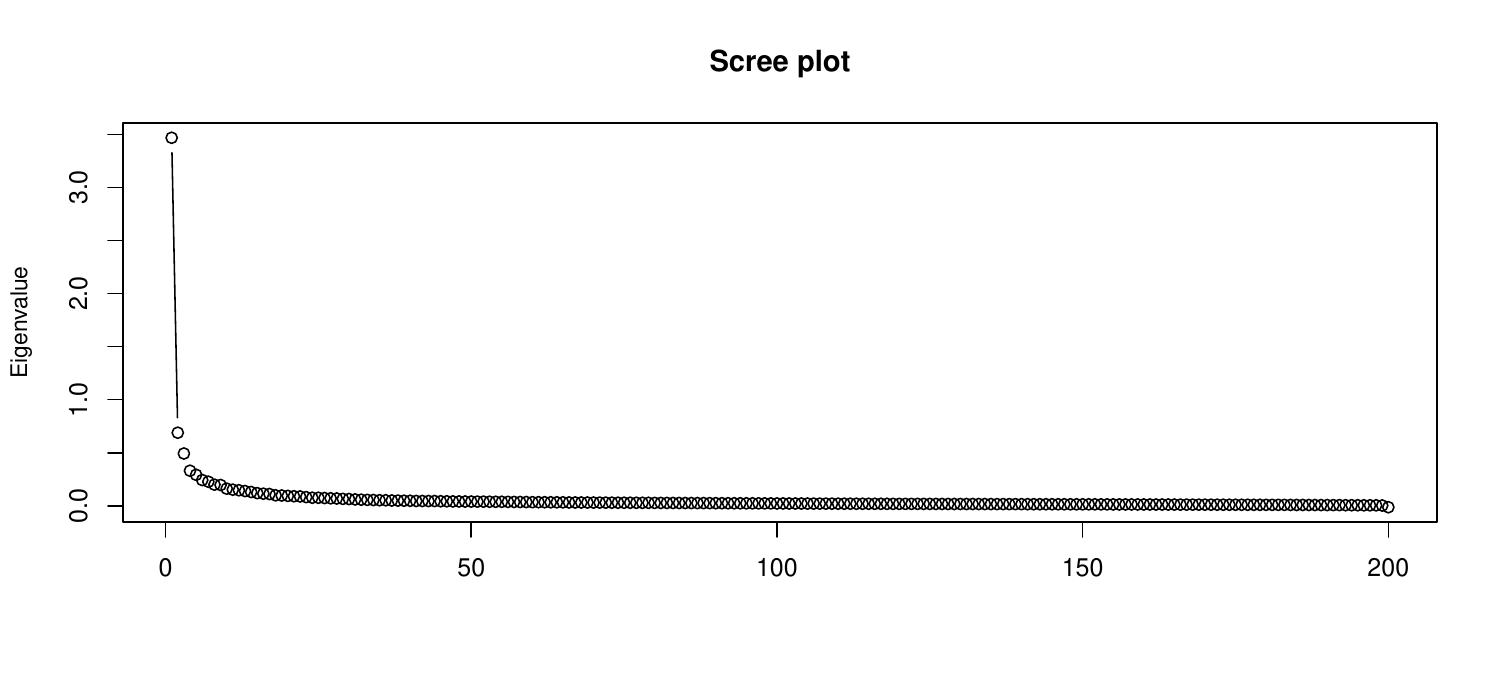}
\vspace{-0.65cm}
\caption{The scree plot of eigenvalues of the sum of 501 ARPM estimates.  }\label{Fig-6}
\end{figure}

To apply POET estimation procedures, we first needed to determine the rank $r$.
We calculated 501 daily integrated volatility matrices using the ARPM estimation procedure
with c=0.5, which is used in the simulation study.
Figure \ref{Fig-6} shows the scree plot, drawn using the eigenvalues from the sum of 501 ARPM estimates.
As seen in Figure \ref{Fig-6}, the possible values of the rank $r$ are 1, 2, and 3; hence, we conducted the empirical study for  $r=1,2,3$.
However, we reported the results only for $r=3$ since the empirical results are similar for $r=1,2,3$ and choosing $r=3$ gives the best overall performance.
We note that the errors caused by the underestimation of the number of factors is more severe than that caused by the overestimation \citep{fan2013large}.

To estimate the sparse volatility matrix $\bSigma$, we used the Global Industry Classification Standard (GICS) \citep{fan2016incorporating}.
Specifically, the covariance matrix for idiosyncratic components for the different sectors are set to zero, and those for the same sector are maintained.  
This corresponds to hard-thresholding using the sector information.
To make the estimates positive semi-definite, we projected the POET estimators onto the positive semi-definite cone in the spectral norm.  Then, we adjusted the POET estimators with the identity matrix to prevent the singularity as described in \eqref{eq-5.9}.

To choose the tuning parameter $c$ in the empirical study, we first defined 
\begin{equation*}
 \Lambda^{\text{ARPM}}(c) = \frac{1}{s}\sum_{d=1}^{s}\| \hat{\bGamma}_d^{\text{ARPM}}(c) - \hat{\bGamma}_{d+1}^{\text{PRVM}}\|^2_F, 
\end{equation*}
where $\hat{\bGamma}_d^{\text{ARPM}}(c)$ and $\hat{\bGamma}_d^{\text{PRVM}}$ are the POET estimator from the ARPM estimator with tuning parameter $c$ and POET estimator from the PRVM estimator for the $d$-th day, respectively.
Then, we selected $c$ by minimizing $\Lambda^{\text{ARPM}}(c)$ over $c \in (0, 1).$
We chose the in-sample period as day 1 to day 150 and obtained $c=0.35$ for the ARP estimator.
We note that stationarity is a reasonable assumption on volatility in financial time series, which motivates the above tuning parameter selection procedure.
Similarly, for the URPM estimator, we defined
\begin{equation*}
 \Lambda^{\text{URPM}}(c, \alpha) =  \frac{1}{s}\sum_{d=1}^{s}\| \hat{\bGamma}_d^{\text{URPM}}(c, \alpha) - \hat{\bGamma}_{d+1}^{\text{PRVM}}\|^2_F, 
\end{equation*}
where $\hat{\bGamma}_d^{\text{URPM}}(c, \alpha)$ is the POET estimator from the URPM estimator with tuning parameter $c$ and universal tail index $\alpha$   for the $d$-th day.
We selected $c$ and $\alpha$ by minimizing $ \Lambda^{\text{URPM}}(c, \alpha)$ over $c \in (0, 1)$ and $\alpha \in [1.1, 2]$.
The selected parameters are $c=0.75$  and $\alpha=2$.

To check the performance of the proposed ARP estimation procedure, we first investigated the mean squared prediction error (MSPE) for the POET estimators defined by
\begin{equation}\label{eq-6.1}
	 \text{MSPE}(\hat{\bGamma}) = \frac{1}{d_2 - d_1 +1}\sum_{d=d_1}^{d_2}\| \hat{\bGamma}_d^l - \hat{\bGamma}_{d+1}^{\text{PRVM}}\|^2_F,
\end{equation}
where $(d_2 - d_1 +1)$ is the number of days in the out-of-sample period, 
\begin{equation*}
\hat{\bGamma}_d^l=\dfrac{1}{l}\sum_{m=1}^{l}\hat{\bGamma}_{d-m},
\end{equation*} 
$\hat{\bGamma}_d$ can be POET estimators from the ARPM, URPM, RPRVM, and PRVM estimators for the $d$-th day, and $l$ is the averaging length.
We used $l=1,5,10$, and set the out-of-sample period as day 161 to day 501.
Then, we split the out-of-sample period into two parts. 
The two periods are denoted by period 1 (day 161 to day 331) and period2 (day 332 to day 501).
We note that since the PRVM estimator is not robust, the MSPE in \eqref{eq-6.1} may not be a perfect  measure. 
 Table \ref{Table-2} reports the MSPE results for the POET estimators from the inputs of the ARPM, URPM, RPRVM, and PRVM estimators for $l=1,5,10$.
We find that the ARPM estimator outperforms other estimators for all periods and averaging lengths.
This may be  because  the proposed ARPM estimator can help deal with heterogeneous heavy-tailed distributions and incorporating the heterogeneous heavy-tailedness helps account for the volatility dynamics.

\begin{table}[!ht]
\caption{The MSPEs of  the POET estimators from the ARPM, URPM, RPRVM, and PRVM estimators for the three averaging lengths  (whole period: day 161 to day 501, period 1: day 161 to day 331, period 2: day 332 to day 501).
}\label{Table-2}
\centering
\scalebox{0.88}{
\begin{tabular}{l l c c c c c c c c c c c c c  c c c c c c c c}
\hline
\multicolumn{4} {c}{} & \multicolumn{7} {c}{MSPE $\times 10^4$}   \\ \cline{5-11}
\multicolumn{4} {c}{} & \multicolumn{7} {c}{Estimator}   \\ \cline{5-11}
 &&	 	&&		ARPM	&&	URPM  && RPRVM && PRVM \\ \hline
1 day &&  whole period	&&   6.501 	&&  7.106 	&&	  6.664		&&	  9.748 			\\

 && period 1 		&&		11.632 	&&   12.807 	&&	 11.955	&&	 17.533 			\\

 && period 2 		&&		1.339 	 && 1.371 	&&   1.343		&&	 1.916 			\\

&&	&&	  	&&				&&		&& 					\\

5 day &&   whole period	 &&  6.555	&&		7.083 	&&	  6.604	&&	  7.044			\\

 && period 1 		&&		11.778 	&&	 12.736	&&	 11.851		&&	 12.659		\\

 && period 2		&&		1.300 	&&	1.397 	&&	 1.326		&& 1.395	\\

&&	&&	  	&&				&&		&& 					\\

10 day &&   whole period 	&&		6.875 	&&	 7.360 	&&	  6.939		&&	 7.080		\\

 && period 1 		&&	   12.428 	&&	  13.281	&&	  12.522		&&   12.811		\\

 && period 2			&&  1.290 	&&	 1.404 	&&	 1.324 	&&	 1.316	\\ \hline

\end{tabular}
}
 \end{table}

To check the out-of-sample performance, we applied the ARPM, URPM, RPRVM, and PRVM estimators to the following minimum variance portfolio allocation problem:
\begin{eqnarray} \label{eq6-2}
\min_{\omega} \omega^ \top\hat{\bGamma}_d^l  \omega,\quad \text{ subject to } \omega^\top\mathbf{J}=1 \text{ and } \left\| \omega\right\|_{1} \leq c_{0},
\end{eqnarray}
where $\mathbf{J}=\left(1, \ldots ,1\right)^\top \in \mathbb{R}^{p}$, the gross exposure constraint $c_{0}$ was changed from 1 to 6, 
\begin{equation*}
\hat{\bGamma}_d^l=\dfrac{1}{l}\sum_{m=1}^{l}\hat{\bGamma}_{d-m},
\end{equation*} 
and $\hat{\bGamma}_d$ could be the POET estimators from the ARPM, URPM, RPRVM, and PRVM estimators for the $d$-th day.
To calculate the out-of-sample risks, we constructed the portfolios at the beginning of each trading day using the stock weights \eqref{eq6-2}, calculated using the data from the previous $l$ days.
We then held this for one day and calculated the realized volatility using the 10-min portfolio log-returns.
The average of their square root was used for out-of-sample risk.
We used $l=1,5,10$, and tested the performances for the whole period (day 161 to day 501), period 1 (day 161 to day 331), and period 2 (day 332 to day 501).

\begin{figure}[!ht]
\centering
\includegraphics[width = 0.86\textwidth]{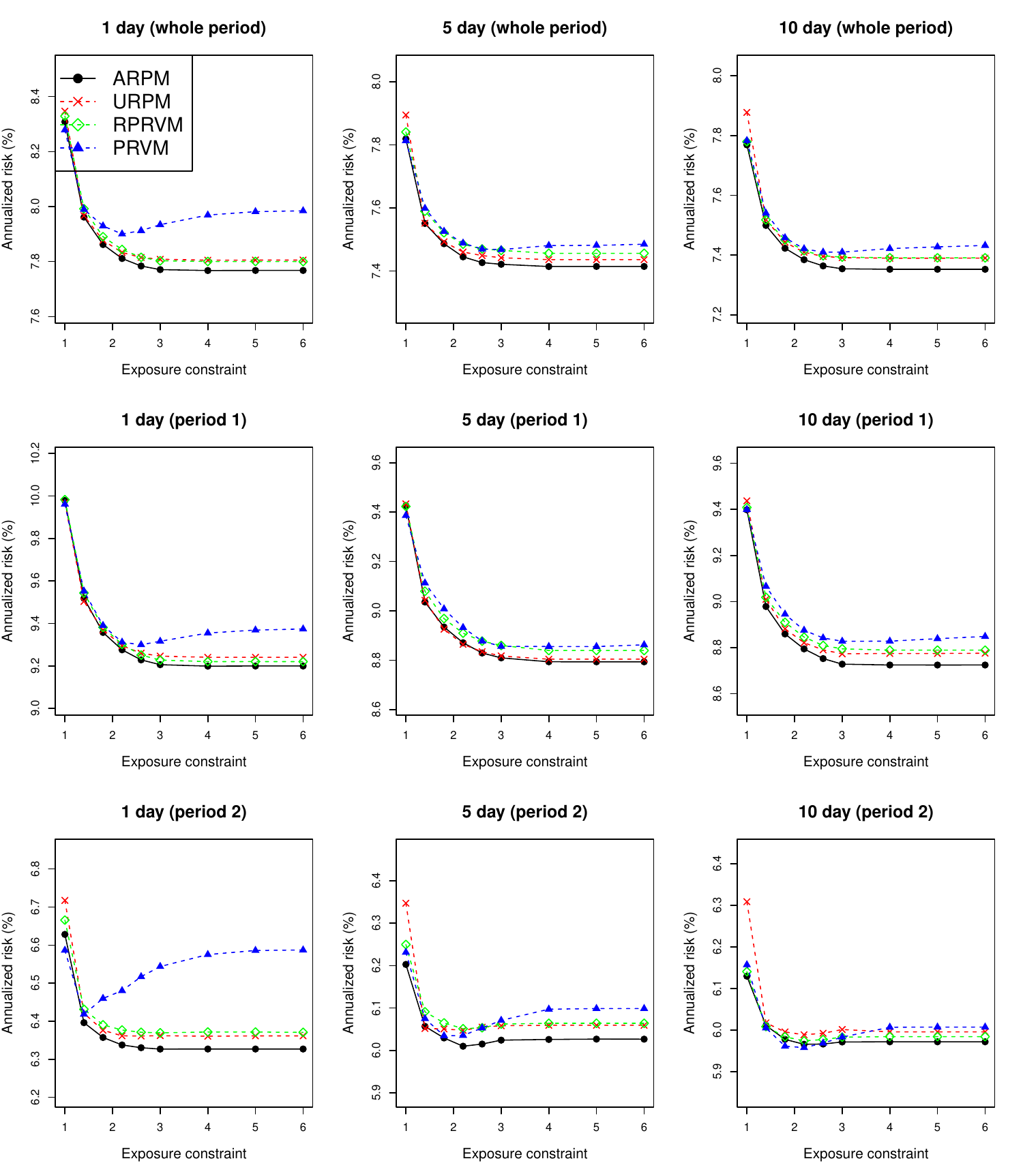}
\caption{The out-of-sample risks of the optimal portfolios constructed  using the POET estimators from the ARPM (black dot), URPM (red cross), RPRVM (green diamond), and PRVM (blue triangle) estimators for the three averaging lengths  (whole period: day 161 to day 501, period 1: day 161 to day 331, period 2: day 332 to day 501).}\label{Fig-7}
\end{figure}
Figure \ref{Fig-7} depicts the out-of-sample risks of the portfolios constructed by the POET estimators from the ARPM, URPM, RPRVM, and PRVM estimators for $l=1,5,10$.
We can find that the averaged volatility matrix predictors $(l = 5, 10$) perform slightly better than the one-day ($l=1$) volatility matrix predictors. 
One possible explanation is that the market is volatile and the one-day volatility is too volatile to be used as a predictor.
When comparing the estimation procedures, the ARPM estimator shows a stable result and performs the best overall.
This result lends further support to our claim that the heavy-tailed distributions of observed log-returns are heterogeneous, as shown in Figure \ref{Fig-1} and Figure \ref{Fig-5}, and that the proposed ARP estimation procedure can account for the heterogeneity of the degrees of heaviness of tail distributions.

\section{Conclusion}\label{SEC-7}

In this paper, we develop the adaptive robust pre-averaging realized volatility (ARP) estimation method to handle the heterogeneous heavy-tailed distributions of stock-returns.
To account for the heterogeneity of the heavy-tailedness from microstructural noises and price jumps, the ARP estimator truncates quadratic pre-averaged random variables according to daily tail indices.
We demonstrate that the proposed ARP estimator achieves sub-Weibull tail concentration with the optimal convergence rate by showing that its upper bound matches its lower bound.
To estimate large integrated volatility matrices, the ARP estimator is further regularized using the POET procedure, and the asymptotic properties of the POET estimator from the ARP estimator are also investigated.
In the empirical study, for the purpose of portfolio allocation, the POET estimator based on the ARP estimator performs best  overall.
These findings suggest that when it comes to estimating the integrated volatility matrices, the proposed ARP estimation procedure helps handle the heterogeneous tail distributions of observed log-returns.

The non-synchronization could be another source of the heavy-tailedness, and the heterogeneity of time intervals can cause some heterogeneous  variation.
However, in this paper, we do not focus on this issue and mainly consider the noise and jump as the source of the heavy-tailedness.
It would be worthwhile to study the observation time point as it relates to the heavy-tailedness.
Furthermore, there are other possible sources of the heavy-tailedness, and it is important to know what actually causes heavy-tailedness.
On the other hand, in this paper, we consider the simple i.i.d. structure for microstructural noises.
However, successfully adjusting the bias term coming from microstructural noises depends on their structure. 
In particular, it is difficult to estimate the cross-sectional dependent structure  with the non-synchronization observations. 
Thus, it is a worthwhile challenge to develop a robust estimator of the cross-sectional dependent structure of microstructural noises.
 We leave these interesting questions for a future study.

\begin{singlespace}
\bibliography{myReferences}
\end{singlespace}

\newpage

\appendix
\section{Appendix}\label{Appendix}

\subsection{Proof of Theorem \ref{Theorem1}}

For simplicity, we denote $K_n$ by $K$.
Let
\begin{eqnarray*}
\bar{X}_i(\tau_k) &=& \sqrt {\frac {n-K}{\phi K}} \sum_{l=0}^{K-1} g \( \frac{l}{K}\) \{ {X}^{c}_i(\tau_{i,k+l+1}) -{X}^{c}_i(\tau_{i,k+l}) \}\cr
&=&\sqrt {\frac {n-K}{\phi K}} \sum_{l=0}^{K-1} g \( \frac{l}{K}\) \int_{\tau_{i,k+l}}^{\tau_{i,k+l+1}}\mu_{i}(t)dt \cr
&&+ \sqrt {\frac {n-K}{\phi K}} \sum_{l=0}^{K-1} g \( \frac{l}{K}\) \int_{\tau_{i,k+l}}^{\tau_{i,k+l+1}}e_{i}^{\top}\bsigma^{\top}(t)d\bW_t \cr
&=& \bar{X}_i^{\mu}(\tau_k)+\bar{X}_i^{\sigma}(\tau_k),
\end{eqnarray*}
where $e_i$ is the $p \times 1$ vector with all 0s except for a 1 in the $i$-th coordinate and
\begin{eqnarray*}
	 \bar{\epsilon}_i(\tau_k) &=& \sqrt {\frac {n-K}{\phi K}} \sum_{l=0}^{K-1} g \( \frac{l}{K}\) \{ \epsilon_i(\tau_{i,k+l+1}) -\epsilon_i(\tau_{i,k+l})\} \cr
	&=&\sqrt {\frac {n-K}{\phi K}} \sum_{l=0}^{K-1} \left \{ g \( \frac{l}{K}\) - g \( \frac{l+1}{K}\) \right \} \epsilon_i(\tau_{i,k+l+1}).
\end{eqnarray*}
Then, we have
\begin{eqnarray} \label{eq-A.1}
	{Q}^{c}_{ij}(\tau_k)=\[\bar{X}_i(\tau_k)+\bar{\epsilon}_i(\tau_k) \]\[\bar{X}_j(\tau_k)+\bar{\epsilon}_j(\tau_k)\].
\end{eqnarray}
Also,  let
\begin{equation*}
{Q}^{(2)}_{ij}(\tau_k) = \[\bar{X}_i(\tau_k)+\bar{\epsilon}_i(\tau_k) + B_i(\tau_{k})\]\[ \bar{X}_j(\tau_k)+\bar{\epsilon}_j(\tau_k) + B_j(\tau_{k})\],
\end{equation*}
\begin{eqnarray*}
&& {Q}^{(2)}_{\rho,ij}(\tau_k) = \frac{1}{2 }  \{ Y_i^{c}(\tau_{i,k+1})- Y_i^{c}(\tau_{i,k}) + J_{i3} (\tau_{i,k+l+1}) - J_{i3} (\tau_{i,k+l}) \} \cr 
&& \quad \quad \quad \quad \quad \quad \times \{ Y_j^{c}(\tau_{j,k+1})- Y_j^{c}(\tau_{j,k}) + J_{j3} (\tau_{i,k+l+1}) - J_{j3} (\tau_{i,k+l})\},
\end{eqnarray*}
where
\begin{equation*}
 B_i(\tau_{k}) = \sqrt {\frac {n-K}{\phi K}}  \sum _{l=0}^{K-1} g\left(\frac{l}{K}\right)\left  \{ J_{i3} (\tau_{i,k+l+1}) - J_{i3} (\tau_{i,k+l}) \right \},
\end{equation*}  
\begin{equation*}
J_{i3}(t) =  \int_{0}^{t} \int_{|x| \leq n^{(1-\alpha_{i})/(4\alpha_{i}-2\alpha_{i} \pi)}} x[N_i(dx,ds) - I_i(dx)ds ].
\end{equation*} 

\begin{proposition}\label{Proposition3}
Under models \eqref{eq-2.1} and \eqref{eq-2.3}, (a) and (b) hold for all $1\leq i,j \leq p$ and sufficiently large $n$.
\begin{enumerate}
\item [(a)] Under Assumption \ref{Assumption1}, there exist positive constants $U_{ij}\left( \tau_{k}\right)$ whose values are free of $n$ and $p$ such that
\begin{equation*}
\E \left\{  \left |{Q}^{(2)}_{ij}(\tau_k)    \right  | ^{\alpha_{ij}} \Big | \FF_{\tau_k} \right\} \leq U_{ij}\left( \tau_{k}\right) \text{ a.s.}
\end{equation*}
for all $1\leq k \leq n-K_{n}$.
\item [(b)] Under Assumption \ref{Assumption1}(a)--(b) and Assumption \ref{Assumption2}, there exist positive constants $U_{\rho,ij}\left( \tau_{k}\right)$ whose values are free of $n$ and $p$ such that
\begin{equation*}
\E \left\{  \left |{Q}^{(2)}_{\rho,ij}(\tau_k)    \right  | ^{\alpha_{ij}} \Big | \FF_{\tau_{k-1}} \right\} \leq U_{\rho,ij}\left( \tau_{k}\right)  \text{ a.s.}
\end{equation*}
for all $1\leq k \leq n-1.$
\end{enumerate}
\end{proposition}
\textbf{Proof of Proposition \ref{Proposition3}.}
First, consider (a). By the H\"older's inequality,
\begin{eqnarray*}
&& \E \Big[  \left |{Q}^{(2)}_{ij}(\tau_k)    \right  | ^{2\alpha _{i}\alpha _{j}/(\alpha _{i}+\alpha _{j})} \Big | \FF_{\tau_{k}} \Big] \cr
&& \leq  \left\{E \Big[  \left |{Q}^{(2)}_{ii}(\tau_k)    \right  | ^{\alpha_{i}} \Big | \FF_{\tau_{k}} \Big] \right\}^{\alpha _{j}/(\alpha _{i}+\alpha _{j})}
\left\{E \Big[  \left |Q_{jj}^{(2)}(\tau_k)    \right  | ^{\alpha_{j}} \Big | \FF_{\tau_{k}} \Big] \right\}^{\alpha _{i}/(\alpha _{i}+\alpha _{j})}
 \text{ a.s.}
\end{eqnarray*}
Therefore, it suffices to show that
\begin{eqnarray*}
\E \Big[  \left |{Q}^{(2)}_{ii}(\tau_k)    \right  | ^{\alpha_{i}} \Big | \FF_{\tau_{k}} \Big] \leq  C  \text{ a.s.}
\end{eqnarray*}
By Assumption \ref{Assumption1}(c) and \eqref{eq-A.1}, we have
\begin{eqnarray*}
	\nu_Q&\geq& \E \Big[  \left |Q_{ii}^{c}(\tau_k)    \right  | ^{\alpha_{i}} \Big] \cr
	&=& \E \[  \left |\bar{X}_i(\tau_k)+\bar{\epsilon}_i(\tau_k)    \right| ^{2\alpha_{i}}  \] \cr
	&=& \E\left\{\E \Big[  \left |\bar{X}_i(\tau_k)+\bar{\epsilon}_i(\tau_k)    \right  | ^{2\alpha_{i}}  \Big| \bar{\epsilon}_i(\tau_k) \Big]\right\} \cr
	&\geq& \E\left\{\left |\E \Big[  \bar{X}_i(\tau_k)+\bar{\epsilon}_i(\tau_k)  | \bar{\epsilon}_i(\tau_k) \Big] \right  |^{2\alpha_{i}}\right\} \cr
	&\geq&  \E\left\{\Big| \left|\bar{\epsilon}_i(\tau_k) \right| -\left|\E\[\bar{X}_i(\tau_k) \]\right| \Big|^{2\alpha_{i}}\right\}.
\end{eqnarray*}
Also, by the fact that $\left|\bar{X}_i^{\mu}(\tau_k)\right| \leq Cn^{-1/4}\text{ a.s.}$, we have
\begin{equation*}
\left|\E\[\bar{X}_i(\tau_k) \]\right| = \left|\E\[ \bar{X}_i^{\mu}(\tau_k) \]\right| \leq Cn^{-1/4} \text{ a.s.}
\end{equation*}
Hence, using the H\"older's inequality, we can show
\begin{equation*}
\E \[\left|\bar{\epsilon}_i(\tau_k) \right|^{2\alpha_i} \] \leq C.
\end{equation*}
Then, by the Lipschitz continuity of $g(\cdot)$, we have
\begin{eqnarray} \label{eq-A.2}
	 \E \Big[  \left |{Q}^{c}_{ii}(\tau_k)    \right  | ^{\alpha_{i}} \Big | \FF_{\tau_{k}} \Big] &=& \E \[  \left |\bar{X}_i^{\mu}(\tau_k)+\bar{X}_i^{\sigma}(\tau_k)+\bar{\epsilon}_i(\tau_k)    \right  | ^{2\alpha_{i}} \Big | \FF_{\tau_{k}} \] \cr
	&\leq& C+ C\E \[  \left |\bar{X}_i^{\sigma}(\tau_k)   \right  | ^{2\alpha_{i}} \Big | \FF_{\tau_{k}} \]    \cr
	&\leq& C+ C\left| \frac {n-K}{\phi K}\sum ^{K-1}_{l=0}g ^{2}\left( \frac {l}{K}\right)\frac {\nu_{\gamma}}{n}\right| ^{\alpha_{i} }  \cr
	&\leq& C \text{ a.s.},
\end{eqnarray}
where the first inequality is due to the H\"older's inequality and the second inequality is from the Burkholder-Davis-Gundy inequality.
For the jump part, by Assumption \ref{Assumption1}(e), we have
\begin{eqnarray*}
 \E \[\left| J_{i3} (\tau_{i,k+l+1}) - J_{i3} (\tau_{i,k+l}) \right|^{2\alpha_i} \] &\leq& Cn^{-1} \int_{0}^{n^{(1-\alpha_{i})/(4\alpha_{i}-2\alpha_{i} \pi)}} |x|^{2\alpha_i} I_i(dx)      \cr
 &\leq&  Cn^{-1+(2\alpha_i -\pi)(1-\alpha_{i})/(4\alpha_{i}-2\alpha_{i} \pi)} \cr
 &\leq& Cn^{-(\alpha_i +1 )/2}
\end{eqnarray*}
and
\begin{eqnarray*}
\E \[\left| J_{i3} (\tau_{i,k+l+1}) - J_{i3} (\tau_{i,k+l}) \right|^{2} \] &\leq& Cn^{-1} \int_{0}^{n^{(1-\alpha_{i})/(4\alpha_{i}-2\alpha_{i} \pi)}} x^{2} I_i(dx)      \cr
 &\leq&  Cn^{-1+(2 -\pi)(1-\alpha_{i})/(4\alpha_{i}-2\alpha_{i} \pi)} \cr
 &\leq& Cn^{-1}.
\end{eqnarray*}
Then, by the Lipschitz continuity of $g(\cdot)$ and Rosenthal's inequality, we have
\begin{eqnarray*}
&& \E \[\left|  B_i (\tau_k)\right|^{2\alpha_i} \Big | \FF_{\tau_{k}} \] \cr
&& \leq Cn^{\alpha_i/2}\left\{ \sum _{l=0}^{K_n-1} \E \[\left| J_{i3} (\tau_{i,k+l+1}) - J_{i3} (\tau_{i,k+l}) \right|^{2\alpha_i} \] +\left(\sum _{l=0}^{K_n-1} \E \[\left| J_{i3} (\tau_{i,k+l+1}) - J_{i3} (\tau_{i,k+l}) \right|^{2}\] \right)^{\alpha_i} \right\} \cr
&& \leq C.
\end{eqnarray*}
Therefore, (a) is proved by H\"older's inequality, and we can show (b) similar to the proof of (a).
\endpf

\textbf{Proof of Theorem \ref{Theorem1}.}
Without loss of generality, we assume that $n = K (L+1)$ for some $L  \in \mathbb{ N}$.
We have
\begin{eqnarray}\label{eq-A.3}
| \hat{T}^{\alpha}_{ij,\theta} - T_{ij} | &\leq& \left| \dfrac {1}{(n-K)\theta_{ij}}\sum_{k=1}^{n-K}\psi_{\alpha_{ij}}\left\{\theta _{ij}Q_{ij}^{(2)}\left(\tau_{k}\right)\right\}-T_{ij} \right|	  \cr
&& + \left| \dfrac {1}{(n-K)\theta_{ij}}\sum_{k=1}^{n-K}\Big[\psi_{\alpha_{ij}}\left\{\theta _{ij}Q_{ij}\left(\tau_{k}\right)\right\}- \psi_{\alpha_{ij}}\left\{\theta _{ij}Q_{ij}^{(2)}\left(\tau_{k}\right)\right\}\Big] \right| \cr
&=& (\uppercase\expandafter{\romannumeral1}) +(\uppercase\expandafter{\romannumeral2}).
\end{eqnarray}
First, consider $(\uppercase\expandafter{\romannumeral1})$.
Let
\begin{equation*}
	 A_{ij}(\tau_k)=\E \Big[ Q_{ij}^{(2)}(\tau_k) \Big | \FF_{\tau_k} \Big].
\end{equation*}
Then for any $s > 0$, we obtain that
\begin{eqnarray} \label{eq-A.4}
	&&\Pr \left \{ \dfrac {1}{(n-K)\theta_{ij}}\sum_{k=1}^{n-K}\psi_{\alpha_{ij}}\left\{\theta _{ij}Q_{ij}^{(2)}\left(\tau_{k}\right)\right\}-\frac{1}{ n-K }  \sum_{k=1}^{n-K}  A_{ij} (\tau_k) \geq \frac{K}{n-K}s  \right \}  \cr
	&&\leq \exp \left\{-\theta _{ij}s\right\} \E  \[ \prod ^{K-1}_{m=0}\prod ^{L-1}_{k=0}\exp \left\{ \frac {1}{K}\Big[ \psi_{\alpha_{ij}} \left\{ \theta _{ij}Q_{ij}^{(2)}\left(\tau_{K k+ m+1}\right) \right\} -\theta _{ij}A_{ij}\left( \tau_{K k+ m+1}\right) \Big]\right\}   \]  \cr
	&&\leq \exp \left\{-\theta _{ij}s\right\} \prod ^{K-1}_{m=0}\E  \[\prod ^{L-1}_{k=0}\exp \Big[ \psi_{\alpha_{ij}} \left\{ \theta _{ij}Q_{ij}^{(2)}\left( \tau_{K k+ m+1}\right) \right\} -\theta _{ij}A_{ij}\left( \tau_{K k+ m+1}\right) \Big]   \] ^{1/K} \cr
	&&= \exp \left\{-\theta _{ij}s\right\}  \prod_{m=0}^{K-1}  \E  \Biggl [  \prod_{k=0}^{L-2}  \exp \Big[ \psi_{\alpha_{ij}} \left\{ \theta _{ij}Q_{ij}^{(2)}\left(\tau_{K k+ m+1}\right) \right\} -\theta _{ij}A_{ij}\left(\tau_{K k+ m+1}\right) \Big]   \cr
	&& \qquad\quad \times  \E \Biggl [\exp \Big[ \psi_{\alpha_{ij}} \left\{ \theta _{ij}Q_{ij}^{(2)}\left(\tau_{K (L-1)+ m+1}\right) \right\} -\theta _{ij}A_{ij}\left(\tau_{K (L-1)+ m+1}\right) \Big]  \Bigg | \FF_{\tau_{K (L-1)+ m+1}}\Biggr ]   \Biggr ] ^{1/K}  \cr
	&&\leq \exp \left\{-\theta _{ij}s\right\}  \prod_{m=0}^{K-1}  \E  \Biggl [  \prod_{k=0}^{L-2}  \exp \Big[ \psi_{\alpha_{ij}} \left\{ \theta _{ij}Q_{ij}^{(2)}\left(\tau_{K k+ m+1}\right) \right\} -\theta _{ij}A_{ij}\left(\tau_{K k+ m+1}\right) \Big]\Biggr ]^{1/K}   \cr
	&& \qquad\qquad\qquad\quad \times  \exp \left\{K^{-1}\sum ^{K-1}_{m=0}c_{\alpha_{ij}}U_{ij}(\tau_{K(L-1)+m+1})\theta^{\alpha_{ij}}_{ij}\right\} \cr
	&&\leq \exp \left\{-\theta _{ij}s+\frac{n-K}{K}c_{\alpha_{ij}}S_{ij} \theta _{ij}^{\alpha_{ij}} \right\}  ,
\end{eqnarray}
where the first and second inequalities are due to the Markov inequality and H\"older's inequality, respectively, and the third and fourth inequalities can be obtained by \eqref{eq-A.5}.
Since we can get $- \log (1- x +c_{\alpha_{ij}}|x|^{\alpha_{ij}}) \leq \psi _{\alpha_{ij} }\left( x\right) \leq \log (1+ x + c_{\alpha_{ij}}|x|^{\alpha_{ij}})$ from Lemma A.2 \citep{minsker2018sub},
we have
\begin{eqnarray} \label{eq-A.5}
	 && \E  \Biggl[      \exp \Big[ \psi_{\alpha_{ij}} \left\{ \theta _{ij}Q_{ij}^{(2)}\left(\tau_{k}\right) \right\} -\theta _{ij}A_{ij}\left(\tau_{k}\right) \Big] \Bigg | \FF_{\tau_k} \Biggr]  \cr
	 && \leq \E  \Biggl[   \exp \Big[ \log \left\{ 1+\theta _{ij}Q_{ij}^{(2)}\left(\tau_{k}\right)+c_{\alpha_{ij}}\left|\theta _{ij}Q_{ij}^{(2)}\left(\tau_{k}\right)\right|^{\alpha_{ij}}\right\}-\theta _{ij}A_{ij}\left(\tau_{k}\right) \Big] \Bigg  | \FF_{\tau_k}  \Biggr]   \cr
	 && = \exp \Big[ \log \left\{ 1+\theta _{ij}A_{ij}\left(\tau_{k}\right)+c_{\alpha_{ij}}\theta _{ij}^{\alpha_{ij}}\E \[  \left |Q_{ij}^{(2)}(\tau_k)    \right  | ^{\alpha_{ij}}  \ | \FF_{\tau_{k}} \]\right\}-\theta _{ij}A_{ij}\left(\tau_{k}\right) \Big]     \cr
	 && \leq  \exp \Big[c_{\alpha_{ij}}\theta _{ij}^{\alpha_{ij}}\E \[  \left |Q_{ij}^{(2)}(\tau_k)    \right  | ^{\alpha_{ij}} \ | \FF_{\tau_{k}}   \] \Big]   \leq  \exp \Big[c_{\alpha_{ij}} U_{ij}(\tau_{k}) \theta _{ij}^{\alpha_{ij}} \Big] \text{ a.s.},
\end{eqnarray}
where the second inequality is due to the fact that $\log \left( 1+x\right) \leq x$ for any $x>-1$, and the last inequality is  from Proposition \ref{Proposition3}.
Choose
\begin{equation*}
\theta_{ij} =    \left(\frac {K\log y ^{-1}}{\left( \alpha_{ij} -1\right) c_{\alpha_{ij} }S_{ij}(n-K)} \right)^{1/\alpha_{ij}}, \quad
s= \left( \frac {\alpha_{ij} ^{\alpha_{ij} }c_{\alpha_{ij} }S_{ij}\left(n-K\right)\left( \log y ^{-1}\right)^{\alpha_{ij}-1} }{\left( \alpha_{ij} -1\right) ^{\alpha_{ij} -1}K}\right)^{1/\alpha_{ij}},
\end{equation*}
where $c\log n \leq \log y^{-1} \leq 2\sqrt{n}$.
Then, we have
\begin{eqnarray*}
	&& \Pr \Bigg[ \dfrac {1}{(n-K)\theta_{ij}}\sum_{k=1}^{n-K}\psi_{\alpha_{ij}}\left\{\theta _{ij}Q_{ij}^{(2)}\left(\tau_{k}\right)\right\}-\frac{1}{ n-K }  \sum_{k=0}^{n-K}  A_{ij} (\tau_k)  \cr
	&& \qquad \qquad   \qquad \qquad \geq\left( \frac {\alpha_{ij} ^{\alpha_{ij} }c_{\alpha_{ij} }S_{ij}K^{\alpha_{ij}-1}\left( \log y ^{-1}\right)^{\alpha_{ij}-1} }{\left( \alpha_{ij} -1\right) ^{\alpha_{ij} -1}\left(n-K\right)^{\alpha_{ij}-1}}\right)^{1/\alpha_{ij}} \Bigg]  \leq y.
\end{eqnarray*}
Similarly, we can show
\begin{eqnarray}\label{eq-A.6}
	&& \Pr \Bigg[ \left|\dfrac {1}{(n-K)\theta_{ij}}\sum_{k=1}^{n-K}\psi_{\alpha_{ij}}\left\{\theta _{ij}Q_{ij}^{(2)}\left(\tau_{k}\right)\right\}-\frac{1}{ n-K }  \sum_{k=0}^{n-K}  A_{ij} (\tau_k)\right|  \cr
	&& \qquad \qquad \qquad\qquad\qquad  \qquad \quad \quad\leq C  \left( n^{-1/2}\log y ^{-1}\right) ^{(\alpha_{ij} -1)/\alpha_{ij} } \Bigg] \geq 1-2y.
\end{eqnarray}

Now, we need to establish the relationship between $  \sum_{k=0}^{n-K}  A_{ij} (\tau_k)/(n-K)$ and $T_{ij}$.
Since $X$, $\epsilon$, and $J_3$ are mutually independent, we have
\begin{eqnarray*}
 A_{ij} (\tau_k) &=& \E\Big[ \bar{X}_i^{\mu}(\tau_k)\bar{X}_j^{\mu}(\tau_k)\Big | \FF_{\tau_k} \Big]+ \E\Big[\bar{X}_i^{\mu}(\tau_k)\bar{X}_j^{\sigma}(\tau_k) +  \bar{X}_i^{\sigma}(\tau_k)\bar{X}_j^{\mu}(\tau_k)  \Big | \FF_{\tau_k} \Big] \cr
 &&+ \E\Big[\bar{X}_i^{\sigma}(\tau_k)\bar{X}_j^{\sigma}(\tau_k)\Big | \FF_{\tau_k} \Big]+  \E\Big[\bar{\epsilon}_i(\tau_k)\bar{\epsilon}_j(\tau_k)\Big | \FF_{\tau_k} \Big] +  \E\Big[B_i(\tau_k) B_j(\tau_k) \Big]  \cr
 &=& (a)+(b)+(c)+(d)+(e).
\end{eqnarray*}
By the fact that $\left|\bar{X}_i^{\mu}(\tau_k)\right| \leq Cn^{-1/4}\text{ a.s.}$, we have
\begin{equation}\label{eq-A.7}
\left|(a)\right| \leq Cn^{-1/2}\text{ a.s.}
\end{equation}
Using the  Burkholder-Davis-Gundy inequality, we can show
\begin{equation}\label{eq-A.8}
\left|(b)\right| \leq Cn^{-1/4}\left(\sqrt{ \E\Big[\left\{\bar{X}_i^{\sigma}(\tau_k)\right\}^2\Big | \FF_{\tau_k} \Big]} +\sqrt{ \E\Big[\left\{\bar{X}_j^{\sigma}(\tau_k)\right\}^2\Big | \FF_{\tau_k} \Big]}\right) \leq Cn^{-1/4}\text{ a.s.}
\end{equation}
Consider  $(c)$.
Let
\begin{equation*}
\bar{X}_i^{\sigma}(\tau_k)=\sqrt {\frac {n-K}{\phi K}} \sum_{l=0}^{K-1}H_{i,k,l},
\end{equation*}
where
\begin{eqnarray*}
H_{i,k,l}&=& g\left(\frac{l}{K}\right) \int_{\tau_{k+l}}^{\tau_{i,k+l+1}}e_{i}^{\top}\bsigma^{\top}(t)d\bW_t
+g\left(\frac{l+1}{K}\right)\int_{\tau_{i,k+l+1}}^{\tau_{k+l+1}}e_{i}^{\top}\bsigma^{\top}(t)d\bW_t \cr
	&=&g\left(\frac{l}{K}\right) \int_{\tau_{k+l}}^{\tau_{k+l+1}}e_{i}^{\top}\bsigma^{\top}(t)d\bW_t
+ \left \{ g\left(\frac{l+1}{K}\right) - g\left(\frac{l}{K}\right)  \right \} \int_{\tau_{i,k+l+1}}^{\tau_{k+l+1}}e_{i}^{\top}\bsigma^{\top}(t)d\bW_t .
\end{eqnarray*}
Then, we have
\begin{equation*}
(c)=\frac {n-K}{\phi K}\sum_{l=0}^{K-1}\E\Big[H_{i,k,l}H_{j,k,l} \Big | \FF_{\tau_k} \Big] \text{ a.s.}
\end{equation*}
By the It\^o's isometry and the boundedness of $\gamma_{ij}(t)$, we can get for all $0\leq l \leq K-1$,
\begin{eqnarray*}
&&\E\Big[H_{i,k,l}H_{j,k,l} \Big | \FF_{\tau_k} \Big] -\left \{ g\left(\frac{l}{K}\right)  \right \} ^2 \E \left \{\int_{\tau_{k+l}} ^{\tau_{k+l+1}}  \gamma_{ij}(t)dt \middle | \FF_{\tau_k} \right \}   \cr
&&\leq Cn^{-1} \[   \left\{g \( \frac{l+1}{K} \)-  g \( \frac{l}{K}\)\right\}^{2}  +  \left |  g\left(\frac{l}{K}\right)  \left\{g \( \frac{l+1}{K} \)-  g \( \frac{l}{K}\)\right\}  \right| \] \cr
&& \leq Cn^{-3/2} \text{ a.s.},
\end{eqnarray*}
where the last inequality is by the piecewise Lipschitz derivative condition for $g(\cdot)$.
Thus, we have
\begin{equation}\label{eq-A.9}
\left|(c)- \frac{n-K}{ \phi K} \sum_{l=0}^{K-1} \left \{ g\left(\frac{l}{K}\right)  \right \} ^2 \E \left \{\int_{\tau_{k+l}} ^{\tau_{k+l+1}}  \gamma_{ij}(t)dt \middle | \FF_{\tau_k} \right \}\right| \leq Cn^{-1/2}\text{ a.s.}
\end{equation}
For $(d)$,  we have
\begin{equation}\label{eq-A.10}
(d)= \frac{n-K}{ \phi K}\sum_{l=0}^{K-1} \left \{ g \( \frac{l}{K} \) - g \( \frac{l+1}{K}\) \right \} ^2 \mathbf{1}(\tau_{i,k+l+1}=\tau_{j,k+l+1})  \eta_{ij} \text{ a.s.}
\end{equation}
Finally, consider $(e)$. Note that
\begin{eqnarray*}
\E \[\left| J_{i3} (\tau_{i,k+l+1}) - J_{i3} (\tau_{i,k+l}) \right|^{2} \] &\leq& Cn^{-1} \int_{0}^{n^{(1-\alpha_{i})/(4\alpha_{i}-2\alpha_{i} \pi)}} x^{2} I_i(dx)      \cr
 &\leq&  Cn^{-1+(2 -\pi)(1-\alpha_{i})/(4\alpha_{i}-2\alpha_{i} \pi)}.
\end{eqnarray*}
Hence, by the Lipschitz continuity of $g(\cdot)$, we have
\begin{eqnarray*}
\E\Big[B_i^2(\tau_k)\Big] &\leq&   Cn^{(1-\alpha_i)/(2\alpha_i)} \cr
&\leq&   Cn^{(1-\alpha_{ii})/(2\alpha_{ii})},
\end{eqnarray*}
which implies
\begin{equation}\label{eq-A.11}
\left|(e)\right| \leq C  n^{(1-\alpha_{ij})/(2\alpha_{ij}) }.
\end{equation}
Combining \eqref{eq-A.7}--\eqref{eq-A.11}, we have
\begin{equation} \label{eq-A.12}
\left|\frac{1}{ n-K }  \sum_{k=1}^{n-K}  A_{ij} (\tau_k)-A_{ij}^* \right| \leq C  n^{(1-\alpha_{ij})/(2\alpha_{ij}) } \text{ a.s.},
\end{equation}
where
\begin{equation*}
	  A_{ij}^* = \frac{1}{ \phi K}  \sum_{l=0}^{K-1} \left \{ g\left(\frac{l}{K}\right)  \right \} ^2 \sum_{k=1} ^{n-K}  \E \left \{\int_{\tau_{k+l}} ^{\tau_{k+l+1}}  \gamma_{ij}(t)dt \middle | \FF_{\tau_k} \right \} +\rho_{ij}.
\end{equation*}

Now, we investigate the relationship between $A_{ij}^*$ and $T_{ij}$.
Note that $\gamma_{ij}(t)$ is bounded and $\sum_{k=1} ^{n-K} \[  \E \left \{\int_{\tau_{k+l}} ^{\tau_{k+l+1}}  \gamma_{ij}(t)dt \middle | \FF_{\tau_k} \right \}   -\int_{\tau_{k+l}} ^{\tau_{k+l+1}}  \gamma_{ij}(t)dt \]$ is the sum of $l+1$ martingales.
Hence, using the Azuma-Hoeffding inequality for each martingale, we can show for all $0\leq l \leq K-1$,
\begin{equation*}
	\Pr \( \left | \sum_{k=1} ^{n-K} \[  \E \left \{\int_{\tau_{k+l}} ^{\tau_{k+l+1}}  \gamma_{ij}(t)dt \middle | \FF_{\tau_k} \right \}   -\int_{\tau_{k+l}} ^{\tau_{k+l+1}}  \gamma_{ij}(t)dt \] \right | \geq C \left( n^{-1/2}\log y ^{-1}\right) ^{1/2 }  \)  \leq Ky.
\end{equation*}
Also, simple algebraic manipulations show
\begin{eqnarray*}
	&&|A_{ij}^* - T_{ij}| \cr
	&&= \left | \frac{1}{ \phi K}  \sum_{l=0}^{K-1} \left \{ g\left(\frac{l}{K}\right)  \right \} ^2 \left( \[\sum_{k=1} ^{n-K}   \E \left \{\int_{\tau_{k+l}} ^{\tau_{k+l+1}}  \gamma_{ij}(t)dt \middle | \FF_{\tau_k} \right \} \]  -\int_{0} ^{1}  \gamma_{ij}(t)dt \right) \right |   \cr
	&&\leq \left | \frac{1}{ \phi K}  \sum_{l=0}^{K-1} \left \{ g\left(\frac{l}{K}\right)  \right \} ^2 \sum_{k=1} ^{n-K} \[  \E \left \{\int_{\tau_{k+l}} ^{\tau_{k+l+1}}  \gamma_{ij}(t)dt \middle | \FF_{\tau_k} \right \}   -\int_{\tau_{k+l}} ^{\tau_{k+l+1}}  \gamma_{ij}(t)dt \] \right |  + 2 \frac{K} {n} \nu_{\gamma}.
\end{eqnarray*}
Therefore, we have
\begin{equation} \label{eq-A.13}
	\Pr \left \{|A_{ij} ^* -T_{ij} | \leq  C \left( n^{-1/2}\log y ^{-1}\right) ^{1/2 } + 2 \frac{K} {n} \nu_{\gamma} \right \} \geq 1-K^{2}y .
\end{equation}
Combining \eqref{eq-A.6}, \eqref{eq-A.12}, and \eqref{eq-A.13}, we have
\begin{equation}\label{eq-A.14}
		\Pr  \left \{   (\uppercase\expandafter{\romannumeral1}) \leq  C  \left( n^{-1/2}\log y ^{-1}\right) ^{(\alpha_{ij} -1)/\alpha_{ij} }
		  \right \}  \geq  1 - 2K^{2} y.
\end{equation}

Consider $(\uppercase\expandafter{\romannumeral2})$.
Since $I_i(|x|>n^{(1-\alpha_{i})/(4\alpha_{i}-2\alpha_{i} \pi)}) < \infty$, for each asset, we have, with probability at least $1-y$, the number of jumps with size bigger than $n^{(1-\alpha_{i})/(4\alpha_{i}-2\alpha_{i} \pi)}$ is bounded by $C\log y^{-1}$.
Hence, using the fact that $\psi_{\alpha}(x)$ is a bounded function, we have
\begin{equation*}
		\Pr  \left\{ (\uppercase\expandafter{\romannumeral2}) \leq   C\left(\frac{K\log y^{-1}}{n\theta_{ij}}\right)
		  \right\}  \geq  1-2y,
\end{equation*}
which implies
\begin{equation}\label{eq-A.15}
		\Pr  \left\{ (\uppercase\expandafter{\romannumeral2}) \leq C  \left( n^{-1/2}\log y ^{-1}\right) ^{(\alpha_{ij} -1)/\alpha_{ij}}
		  \right\}  \geq  1-2y.
\end{equation}
Collecting \eqref{eq-A.3}, \eqref{eq-A.14}, and \eqref{eq-A.15},
we obtain that with probability at least $1-3K^{2} y$,
\begin{equation*}
		| \hat{T}^{\alpha}_{ij,\theta} - T_{ij} | \leq C  \left( n^{-1/2}\log y ^{-1}\right) ^{(\alpha_{ij} -1)/\alpha_{ij}},
\end{equation*}
and then substituting $\delta/\(3K^{2}\)$ for $y$ completes the proof.
\endpf

\subsection{Proof of Theorem \ref{Theorem2}}
\textbf{Proof of Theorem \ref{Theorem2}.}
Let $n_i=n_j=n$ and $t_{i,k}=t_{j,k}=\tau_{i,k}=\tau_{j,k}=\tau_k=k/n$ for $1 \leq k \leq n$.
To derive a lower bound, we construct two quadratic pre-averaged random variables $Q_{1,ij}\left( \tau_{k}\right)$ and $Q_{2,ij}\left( \tau_{k}\right)$ as follows.
Let  $d \bX_{1}(t) = d \bX_{2}(t) = \bsigma^\top(t)d \bW_t$ for any appropriate $\bsigma(t)$, which implies  $\bar{X}_{1,h}(\tau_k)=\bar{X}_{2,h}(\tau_k)$  for  $1\leq h \leq p$ and $1\leq k \leq n-K$.
Also, let $2\epsilon_{1,h}(t_{h,k})=\epsilon_{2,h}(t_{h,k})$  for  $1\leq h \leq p$ and $0\leq k \leq n$, where the distributions of $\epsilon_{1,h}(t_{h,k})$,  $1\leq h \leq p$ are defined as follows:
\begin{eqnarray*}
	\epsilon_{1,h}(t_{h,k})=
 \begin{cases}
  K^{(\alpha_{h}+1)/2\alpha_{h}}\left( \log (1/2\delta) \right) ^{-1/2\alpha_{h} }  & \text{ with probability } d   \\
 0 & \text{ with probability }1-2d  \\
 -K^{(\alpha_{h}+1)/2\alpha_{h}}\left( \log (1/2\delta) \right) ^{-1/2\alpha_{h} } & \text{ with probability } d, \\
\end{cases}
\end{eqnarray*}
where $d=C^{2}_{K}\log (1/2\delta)/4K^{2} $. 
For each $0\leq k \leq n$,  let $\Pr \left\{ \epsilon_{1,h}(t_{h,k})>0 \text{ for all }1\leq h \leq p \right\}$ $=\Pr \left\{ \epsilon_{1,h}(t_{h,k})<0 \text{ for all }1\leq h \leq p \right\}=d$ and  $\Pr \left\{ \epsilon_{1,h}(t_{h,k})=0  \text{ for all }1\leq h \leq p \right\}=1-2d$.
Then, using the fact that $1-x\geq \exp(-x/(1-x))$ for any $0\leq x \leq 1/2$, we can show
\begin{eqnarray}\label{eq-A.16}
\prod_{k=1}^{n}\Pr \left\{ \epsilon_{1,i}(\tau_k)=\epsilon_{1,j}(\tau_k)=\epsilon_{2,i}(\tau_k)=\epsilon_{2,j}(\tau_k) = 0\right\}
=\left( 1-\frac {1}{2n} \log \frac {1}{2\delta}\right) ^{n}\geq 2\delta.
\end{eqnarray}
Here, we need to check whether the construction satisfies Assumption \ref{Assumption1}(c).
It suffices to show
\begin{eqnarray}\label{eq-A.17}
\E \[  \left |Q_{1,ii}^{(2)}\left( \tau_{1}\right)   \right  | ^{\alpha_{i}} \] \leq  C.
\end{eqnarray}
Note that for all $1\leq k \leq n$,
\begin{eqnarray*}
 \E \[  \left |\epsilon_{1,i}\left( \tau_{k}\right)   \right  | ^{2\alpha_{i}} \]= \frac {C^{2}_{K}}{2}K ^{\alpha_{i}-1 } \quad \text{and} \quad \E \[  \epsilon_{1,i} ^{2}\left( \tau_{k}\right)    \]= \frac {C^{2}_{K}}{2}\left( K^{-1}\log \frac {1}{2\delta }\right) ^{(\alpha_{i}-1)/\alpha_{i} }.
\end{eqnarray*}
Hence, by the Lipschitz continuity of $g(\cdot)$, we have
\begin{eqnarray}\label{eq-A.18}
\E \[  \left |\bar{\epsilon}_{1,i}(\tau_1)   \right  | ^{2\alpha_{i}} \]
&=& \left( \sqrt {\frac {n-K}{\phi K}}\frac {1}{K}\right) ^{2\alpha_{i}}
    \E \Bigg[\left |  \sum_{l=0}^{K-1} K\left \{ g \( \frac{l}{K}\) - g \( \frac{l+1}{K}\) \right \} \epsilon_{1,i}(\tau_{l+2})  \right  | ^{2\alpha_{i}}\Bigg]  \cr
&\leq& CK^{-\alpha_{i}}\left\{\sum_{l=0}^{K-1}\E \[  \left |\epsilon_{1,i}\left( \tau_{l+2}\right)   \right  | ^{2\alpha_{i}} \]
        +\left(\sum_{l=0}^{K-1}\E \[  \epsilon_{1,i}^{2}\left( \tau_{l+2}\right)    \]\right)^{\alpha_{i}} \right\} \cr
&\leq& C+CK^{-\alpha_{i}}\left\{K\left( K^{-1}\log \frac {1}{2\delta}\right) ^{\left( \alpha_{i}-1\right) /\alpha_{i} }\right\}^{\alpha_{i}} \leq C,
\end{eqnarray}
where the first inequality is due to the Rosenthal's inequality.
Then similar to the proof of  Proposition \ref{Proposition3}, we can show \eqref{eq-A.17}.

Now, since
\begin{eqnarray*}
 \left |T_{1,ij}-T_{2,ij} \right| =\frac {3n\zeta }{\phi K}\E \[  \epsilon_{1,i}\left( \tau_{1}\right)\epsilon_{1,j}\left( \tau_{1}\right)   \],
\end{eqnarray*}
we have for any $\hat{T}_{ij}\left( Q_{ij}\left( \tau_{k}\right) ,\delta \right)$,
\begin{eqnarray}\label{eq-A.19}
&&\max\Bigg[\Pr  \left\{   \left| \hat{T}_{ij}\left( Q_{1,ij}\left( \tau_{k}\right) ,\delta \right) - T_{1,ij} \right| \geq   \frac {n\zeta}{\phi K} \E \[   \epsilon_{1,i}\left( \tau_{1}\right)\epsilon_{1,j}\left( \tau_{1}\right)   \] \right\}, \cr
&& \qquad \quad \Pr  \left\{   \left| \hat{T}_{ij}\left( Q_{2,ij}\left( \tau_{k}\right) ,\delta \right) - T_{2,ij} \right| \geq   \frac {n\zeta}{\phi K} \E \[  \epsilon_{1,i}\left( \tau_{1}\right)\epsilon_{1,j}\left( \tau_{1}\right) \]\right\} \Bigg]\cr
&& \geq \frac{1}{2} \Pr  \Bigg[  \left| \hat{T}_{ij}\left( Q_{1,ij}\left( \tau_{k}\right) ,\delta \right) - T_{1,ij} \right| \geq   \frac {n\zeta}{\phi K} \E \[   \epsilon_{1,i}\left( \tau_{1}\right)\epsilon_{1,j}\left(\tau_{1}\right)   \] \cr
&&\qquad  \qquad \text{ or }  \left| \hat{T}_{ij}\left( Q_{2,ij}\left( \tau_{k}\right) ,\delta \right) - T_{2,ij} \right| \geq   \frac {n\zeta}{\phi K} \E \[   \epsilon_{1,i}\left( \tau_{1}\right)\epsilon_{1,j}\left(\tau_{1}\right)\] \Bigg]\cr
&&\geq \frac{1}{2} \Pr  \left\{ \hat{T}_{ij}\left( Q_{1,ij}\left(\tau_{k}\right) ,\delta \right) =\hat{T}_{ij}\left( Q_{2,ij}\left(\tau_{k}\right) ,\delta \right)    \right\}    \cr
&& \geq \frac{1}{2}\prod_{k=1}^{n}\Pr \left\{ \epsilon_{1,i}(\tau_k)=\epsilon_{1,j}(\tau_k)=\epsilon_{2,i}(\tau_k)=\epsilon_{2,j}(\tau_k) = 0\right\} \geq \delta,
\end{eqnarray}
where the last inequality is from \eqref{eq-A.16}.
Combining \eqref{eq-A.19} and the fact that $\E \[   \epsilon_{1,i}\left( \tau_{1}\right)\epsilon_{1,j}\left( \tau_{1}\right)   \]=C  \left( K^{-1}\log (1/2\delta) \right)^{(\alpha_{ij} -1)/\alpha_{ij} }$, we have for sufficiently large $n$,
\begin{eqnarray}\label{eq-A.20}
&&\max\Bigg[\Pr  \left\{   \left| \hat{T}_{ij}\left( Q_{1,ij}\left( \tau_{k}\right) ,\delta \right) - T_{1,ij} \right| \geq    C  \left( n^{-1/2}\log \frac {1}{2\delta} \right)
 ^{(\alpha_{ij} -1)/\alpha_{ij} }  \right\}, \cr
&& \qquad \quad  \Pr  \left\{   \left| \hat{T}_{ij}\left( Q_{2,ij}\left( \tau_{k}\right) ,\delta \right) - T_{2,ij} \right| \geq    C  \left( n^{-1/2}\log \frac {1}{2\delta} \right)
^{(\alpha_{ij} -1)/\alpha_{ij} } \right\} \Bigg] \cr
&&\geq \delta,
\end{eqnarray}
which completes the proof.
\endpf

\subsection{Proof of Theorem \ref{Theorem3}}

\begin{proposition}\label{Proposition4}
Under Assumption \ref{Assumption1}(a),(b),(d) and Assumption \ref{Assumption2},
Assumption \ref{Assumption1}(c) is satisfied.
\end{proposition}

\textbf{Proof of Proposition \ref{Proposition4}.}
Similar to  the proof of Proposition \ref{Proposition3}, we can show
\begin{equation*}
\E \[  \left |\epsilon_{i}\left( \tau_{i,k}\right)   \right  | ^{2\alpha_{i}} \]\leq C.
\end{equation*}
Then, we have
\begin{eqnarray*}
	 \E \Big[  \left |{Q}^{c}_{ii}(\tau_k)    \right  | ^{\alpha_{i}} \Big] &=& \E \[  \left |\bar{X}_i^{\mu}(\tau_k)+\bar{X}_i^{\sigma}(\tau_k)+\bar{\epsilon}_i(\tau_k)    \right  | ^{2\alpha_{i}} \] \cr
	&\leq& C+ C\E \[  \left |\bar{X}_i^{\sigma}(\tau_k)   \right  | ^{2\alpha_{i}}  \] + C\E \[  \left |\bar{\epsilon}_i(\tau_k)   \right  | ^{2\alpha_{i}}  \]    \cr
	&\leq& C+ C\left| \frac {n-K}{\phi K}\sum ^{K-1}_{l=0}g ^{2}\left( \frac {l}{k}\right)\frac {\nu_{\gamma}}{n}\right| ^{\alpha_{i} }\cr
	&&+C\E\[ \left | \frac {n-K}{\phi K}  \sum_{l=0}^{K-1} \left \{ g \( \frac{l}{K}\) - g \( \frac{l+1}{K}\) \right \}^{2} \epsilon_i^{2}(\tau_{i,k+l+1})  \right  | ^{\alpha_{i}}\] \cr
	&\leq& C+CK^{-\alpha_{i}}\E \[ \left |  \sum_{l=0}^{K-1}  \epsilon_i^{2}(\tau_{i,k+l+1})  \right  | ^{\alpha_{i}} \] \cr
	&\leq& C+CK^{-1} \sum_{l=0}^{K-1}\E\[ \left |  \epsilon_i(\tau_{i,k+l+1})  \right  | ^{2\alpha_{i}} \] \leq C \text{ a.s.},
\end{eqnarray*}
where the first and fourth inequalities are due to the H\"older's inequality, and the second inequality is from the Burkholder-Davis-Gundy inequality.
\endpf

\textbf{Proof of Theorem \ref{Theorem3}.}
Without loss of generality, we assume that $n = 2L+1 $ for some $L  \in \mathbb{ N}$.
We have
\begin{eqnarray}\label{eq-A.21}
| \hat{\rho}^{\alpha}_{ij,\theta} - \rho_{ij} | &\leq& \left| \dfrac {\zeta}{\phi K\theta_{\rho,ij}}\sum_{k=1}^{n-1}\psi_{\alpha_{ij}}\left\{\theta _{\rho,ij}Q_{\rho,ij}^{(2)}\left(\tau_{k}\right)\right\}-\rho_{ij} \right|	  \cr
&& + \left| \dfrac {\zeta}{\phi K\theta_{\rho,ij}}\sum_{k=1}^{n-1}\Big[\psi_{\alpha_{ij}}\left\{\theta _{\rho,ij}Q_{\rho,ij}\left(\tau_{k}\right)\right\}- \psi_{\alpha_{ij}}\left\{\theta _{\rho,ij}Q_{\rho,ij}^{(2)}\left(\tau_{k}\right)\right\}\Big] \right| \cr
&=& (\uppercase\expandafter{\romannumeral1}) +(\uppercase\expandafter{\romannumeral2}).
\end{eqnarray}
First, consider $(\uppercase\expandafter{\romannumeral1})$.
Let
$$(\zeta/\phi K\theta_{\rho,ij})\sum_{k=1}^{n-1}\psi_{\alpha_{ij}}\left\{\theta _{\rho,ij}Q_{\rho,ij}^{(2)}\left(\tau_{k}\right)\right\} = \hat{\rho}^{\alpha }_{1,ij, \theta} + \hat{\rho}^{\alpha }_{2,ij, \theta},
$$ where
\begin{eqnarray*}
	&&\hat{\rho}^{\alpha }_{1,ij, \theta} =  \dfrac {\zeta}{\phi K\theta_{\rho,ij}}\sum_{k=1}^{L} \psi_{\alpha_{ij}}\left\{\theta _{\rho,ij}Q_{\rho,ij}^{(2)}\left( \tau_{2k-1}\right)\right\}, \cr
	&&\hat{\rho}^{\alpha }_{2,ij, \theta} =  \dfrac {\zeta}{\phi K\theta_{\rho,ij}}\sum_{k=1}^{L} \psi_{\alpha_{ij}}\left\{\theta _{\rho,ij}Q_{\rho,ij}^{(2)}\left( \tau_{2k}\right)\right\}.
\end{eqnarray*}
Also, define
\begin{equation*}
	  A_{\rho, ij}(\tau_k)= \E \Big[ Q_{\rho,ij}^{(2)}\left( \tau_{k}\right) \Big | \FF_{\tau_{k-1}} \Big].
\end{equation*}
Then, we can show for any $s>0$,
\begin{eqnarray*}
	&&\Pr \left \{ \hat{\rho}^{\alpha }_{1,ij, \theta}-\frac{\zeta}{ \phi K }  \sum_{k=1}^{L}  A_{\rho, ij}(\tau_{2k-1}) \geq \frac{\zeta s}{\phi K}  \right \}  \cr
	&&\leq \exp \left\{-\theta_{\rho,ij}s\right\}\E  \[\exp \left\{ \sum ^{L}_{k=1}\Big[ \psi_{\alpha_{ij}} \left\{ \theta_{\rho,ij}Q_{\rho,ij}^{(2)}\left( \tau_{2k-1}\right)\right\} -\theta_{\rho,ij}A_{\rho, ij}(\tau_{2k-1}) \Big]\right\}   \]  \cr
	&&= \exp \left\{-\theta_{\rho,ij}s\right\}   \E  \Biggl [    \exp \left\{ \sum ^{L-1}_{k=1}\Big[ \psi_{\alpha_{ij}} \left\{ \theta_{\rho,ij}Q_{\rho,ij}^{(2)}\left( \tau_{2k-1}\right)\right\} -\theta_{\rho,ij}A_{\rho, ij}(\tau_{2k-1}) \Big]\right\}   \cr
	&& \qquad \qquad\qquad \qquad \times  \E \Biggl [ \exp \Big[ \psi_{\alpha_{ij}} \left\{ \theta_{\rho,ij}Q_{\rho,ij}^{(2)}\left( \tau_{2L-1}\right)\right\} -\theta_{\rho,ij}A_{\rho, ij}(\tau_{2L-1}) \Big]  \Bigg | \FF_{\tau_{2L-2}}\Biggr ]   \Biggr ]   \cr
	&&\leq \exp \left\{-\theta_{\rho,ij}s\right\}   \E  \Biggl [    \exp \left\{ \sum ^{L-1}_{k=1}\Big[ \psi_{\alpha_{ij}} \left\{ \theta_{\rho,ij}Q_{\rho,ij}^{(2)}\left( \tau_{2k-1}\right)\right\} -\theta_{\rho,ij}A_{\rho, ij}(\tau_{2k-1}) \Big]\right\}\Biggr]   \cr
	&& \qquad\qquad\qquad \qquad \times  \exp \left\{c_{\alpha_{ij}}U_{\rho,ij}(\tau_{2L-2})\theta^{\alpha_{ij}}_{\rho,ij}\right\} \cr
	&&\leq \exp \left\{-\theta_{\rho,ij}s+ (n-1)c_{\alpha_{ij}}S_{\rho,ij} \theta^{\alpha_{ij}}_{\rho,ij}   \right\}.
\end{eqnarray*}
Choose
\begin{equation*}
\theta_{\rho,ij} =    \left(\frac {\log y ^{-1}}{\left( \alpha_{ij} -1\right) c_{\alpha_{ij} }S_{\rho,ij}(n-1)} \right)^{1/\alpha_{ij}}, \quad
s= \left( \frac {\alpha_{ij} ^{\alpha_{ij} }c_{\alpha_{ij} }S_{\rho,ij}(n-1)\left( \log y ^{-1}\right)^{\alpha_{ij}-1} }{\left( \alpha_{ij} -1\right) ^{\alpha_{ij} -1}}\right)^{1/\alpha_{ij}},
\end{equation*}
where $c\log n \leq \log y^{-1} \leq 2\sqrt{n}$.
Then, we have
\begin{eqnarray*}
&&\Pr \left \{ \hat{\rho}^{\alpha }_{1,ij, \theta}-\frac{\zeta}{ \phi K }  \sum_{k=1}^{L}  A_{\rho, ij}(\tau_{2k-1}) \geq C  \left( n^{-1}\log y ^{-1}\right) ^{(\alpha_{ij} -1)/\alpha_{ij} }  \right \} \leq y.
\end{eqnarray*}
Similarly, we can show
\begin{eqnarray}\label{eq-A.22}
&&\Pr \Bigg[ \left|\dfrac {\zeta}{\phi K\theta_{\rho,ij}}\sum_{k=1}^{n-1}\psi_{\alpha_{ij}}\left\{\theta _{\rho,ij}Q_{\rho,ij}^{(2)}\left(\tau_{k}\right)\right\}-\frac{\zeta}{ \phi K }  \sum_{k=1}^{n-1}  A_{\rho, ij}(\tau_{k})\right| \cr
&&  \qquad \qquad \qquad   \qquad \qquad \qquad \quad \geq C  \left( n^{-1}\log y ^{-1}\right) ^{(\alpha_{ij} -1)/\alpha_{ij} }  \Bigg] \leq 4y.
\end{eqnarray}
Now, we need to establish the relationship between $ \zeta\sum_{k=1}^{n-1}A_{\rho, ij}(\tau_{k})/(\phi K)$ and $\rho_{ij}$.
Since $X$, $\epsilon$, and $J_3$ are mutually independent, similar to the proof of Theorem \ref{Theorem1}, we can show
\begin{equation}\label{eq-A.23}
\left|\rho_{ij}-\frac{\zeta}{ \phi K }  \sum_{k=1}^{n-1}  A_{\rho, ij}(\tau_{k})\right| \leq Cn^{-1}\text{ a.s.}
\end{equation}
Combining \eqref{eq-A.22} and \eqref{eq-A.23}, we have
\begin{equation}\label{eq-A.24}
		\Pr  \left \{   (\uppercase\expandafter{\romannumeral1}) \leq  C  \left( n^{-1}\log y ^{-1}\right) ^{(\alpha_{ij} -1)/\alpha_{ij} }
		  \right \}  \geq  1 - 4y.
\end{equation}

Consider $(\uppercase\expandafter{\romannumeral2})$. 
Note that $I_i(|x|>n^{(1-\alpha_{i})/(4\alpha_{i}-2\alpha_{i} \pi)}) < \infty$. 
Hence, for each asset, we have, with probability at least $1-y$, the number of jumps with size bigger than $n^{(1-\alpha_{i})/(4\alpha_{i}-2\alpha_{i} \pi)}$  is bounded by $C\log y^{-1}$.
Then, using the fact that $\psi_{\alpha}(x)$ is a bounded function, we have
\begin{equation*}
		\Pr  \left\{ (\uppercase\expandafter{\romannumeral2}) \leq   C\left(\frac{\log y^{-1}}{n\theta_{\rho,ij}}\right)
		  \right\}  \geq  1-2y,
\end{equation*}
which implies
\begin{equation}\label{eq-A.25}
		\Pr  \left\{ (\uppercase\expandafter{\romannumeral2}) \leq C  \left( n^{-1}\log y ^{-1}\right) ^{(\alpha_{ij} -1)/\alpha_{ij}}
		  \right\}  \geq  1-2y.
\end{equation}
Collecting \eqref{eq-A.21}, \eqref{eq-A.24}, and \eqref{eq-A.25}, we obtain that with probability at least $1-6y$,
\begin{equation*}
		| \hat{\rho}^{\alpha}_{ij,\theta} - \rho_{ij} | \leq C  \left( n^{-1}\log y ^{-1}\right) ^{(\alpha_{ij} -1)/\alpha_{ij}},
\end{equation*}
and then substituting $\delta/6$ for $y$ completes the proof of \eqref{eq-4.5}.
Also, \eqref{eq-4.6} is proved by \eqref{eq-4.4} and Theorem \ref{Theorem1}.
\endpf

\subsection{Proof of Theorem \ref{Theorem4}}

\textbf{Proof of Theorem \ref{Theorem4}.}
Let $n_i=n_j=n$ and $t_{i,k}=t_{j,k}=\tau_{i,k}=\tau_{j,k}=\tau_k=k/n$ for $1 \leq k \leq n$.
Similar to the proof of Theorem \ref{Theorem2}, we construct two quadratic log-return random variables $Q_{1,\rho,ij}\left( \tau_{k}\right)$ and $Q_{2,\rho,ij}\left( \tau_{k}\right)$ as follows.
Let $d \bX_{1}(t) = d \bX_{2}(t) = \bsigma^\top(t)d \bW_t$ for any appropriate $\bsigma(t)$, which implies $X_{1,h}(\tau_{k+1})-X_{1,h}(\tau_{k})=X_{2,h}(\tau_{k+1})-X_{2,h}(\tau_{k})$ for  $1\leq h \leq p$ and $1\leq k \leq n-1$.
Also, let $2\epsilon_{1,h}(t_{h,k})=\epsilon_{2,h}(t_{h,k})$ for $1\leq h \leq p$  and $0\leq k \leq n$, where the distributions of $\epsilon_{1,h}(t_{h,k})$, $1\leq h \leq p$  are defined as follows:
\begin{eqnarray*}
	\epsilon_{1,h}(t_{h,k})=
 \begin{cases}
  n^{1/2\alpha_{h}}\left( \log (1/2\delta) \right) ^{-1/2\alpha_{h} }  & \text{ with probability } d   \\
 0 & \text{ with probability }1-2d  \\
 -n^{1/2\alpha_{h}}\left( \log (1/2\delta) \right) ^{-1/2\alpha_{h} } & \text{ with probability } d, \\
\end{cases}
\end{eqnarray*}
where $d=\log (1/2\delta)/8n $. 
For each $0\leq k \leq n$,  assume that $\Pr \left\{ \epsilon_{1,h}(t_{h,k})>0 \text{ for all }1\leq h \leq p \right\}$ $=\Pr \left\{ \epsilon_{1,h}(t_{h,k})<0 \text{ for all }1\leq h \leq p \right\}=d$ and  $\Pr \left\{ \epsilon_{1,h}(t_{h,k})=0  \text{ for all }1\leq h \leq p \right\}=1-2d$.
Then, using the fact that $1-x\geq \exp(-x/(1-x))$ for any $0\leq x \leq 1/2$, we can show
\begin{eqnarray}\label{eq-A.26}
\prod_{k=1}^{n} \Pr \left\{ \epsilon_{1,i}(\tau_k)=\epsilon_{1,j}(\tau_k)=\epsilon_{2,i}(\tau_k)=\epsilon_{2,j}(\tau_k) = 0\right\}
=\left( 1-\frac {1}{4n}\log\frac{1}{2\delta}\right) ^{n}\geq 2\delta.
\end{eqnarray}
Here, we need to check whether the construction satisfies Assumption \ref{Assumption2}.
It suffices to show
\begin{eqnarray}\label{eq-A.27}
\E \[  \left |\epsilon_{1,i}\left( \tau_{1}\right)   \right  | ^{2\alpha_{i}} \]\leq C.
\end{eqnarray}
Note that
\begin{eqnarray*}
 \E \[  \left |\epsilon_{1,i}\left( \tau_{1}\right)   \right  | ^{2\alpha_{i}} \]= \frac {1}{4}
\quad \text{and} \quad  \E \[  \epsilon_{1,i}\left( \tau_{1}\right)\epsilon_{1,j}\left( \tau_{1}\right)   \]= \frac {1}{4}\left( n^{-1}\log \frac {1}{2\delta }\right) ^{(\alpha_{ij}-1)/\alpha_{ij} }.
\end{eqnarray*}
Hence, \eqref{eq-A.27} is satisfied, and since
\begin{eqnarray*}
 \left |\rho_{1,ij}-\rho_{2,ij} \right| =\frac {3n\zeta }{\phi K}\E \[  \epsilon_{1,i}\left( \tau_{1}\right)\epsilon_{1,j}\left( \tau_{1}\right)   \],
\end{eqnarray*}
we have for any $\hat{\rho}_{ij}\left( Q_{\rho,ij}\left( \tau_{k}\right) ,\delta \right)$,
\begin{eqnarray}\label{eq-A.28}
&&\max\Bigg[\Pr  \left\{   \left| \hat{\rho}_{ij}\left( Q_{1,\rho,ij}\left( \tau_{k}\right) ,\delta \right) - \rho_{1,ij} \right| \geq   \frac{n\zeta }{4\phi K}\left( n^{-1}\log \frac {1}{2\delta }\right) ^{(\alpha_{ij}-1)/\alpha_{ij} } \right\}, \cr
&& \qquad \quad \Pr  \left\{   \left|  \hat{\rho}_{ij}\left( Q_{2,\rho,ij}\left( \tau_{k}\right) ,\delta \right) - \rho_{2,ij}\right| \geq \frac{n\zeta }{4\phi K}\left( n^{-1}\log \frac {1}{2\delta }\right) ^{(\alpha_{ij}-1)/\alpha_{ij} } \right\} \Bigg]\cr
&& \geq \frac{1}{2} \Pr  \Bigg[  \left|  \hat{\rho}_{ij}\left( Q_{1,\rho,ij}\left( \tau_{k}\right) ,\delta \right) - \rho_{1,ij} \right| \geq   \frac{n\zeta }{4\phi K}\left( n^{-1}\log \frac {1}{2\delta }\right) ^{(\alpha_{ij}-1)/\alpha_{ij} } \cr
&& \qquad \qquad \text{ or } \left|  \hat{\rho}_{ij}\left( Q_{2,\rho,ij}\left( \tau_{k}\right) ,\delta \right) - \rho_{2,ij} \right| \geq  \frac{n\zeta }{4\phi K}\left( n^{-1}\log \frac {1}{2\delta }\right) ^{(\alpha_{ij}-1)/\alpha_{ij} } \Bigg]\cr
&&\geq \frac{1}{2} \Pr  \left\{  \hat{\rho}_{ij}\left( Q_{1,\rho,ij}\left( \tau_{k}\right) ,\delta \right) = \hat{\rho}_{ij}\left( Q_{2,\rho,ij}\left( \tau_{k}\right) ,\delta \right)   \right\}    \cr
&&\geq \frac{1}{2}\prod_{k=0}^{n} \Pr \left\{ \epsilon_{1,i}(\tau_k)=\epsilon_{1,j}(\tau_k)=\epsilon_{2,i}(\tau_k)=\epsilon_{2,j}(\tau_k) = 0\right\} \geq \delta,
\end{eqnarray}
where the last inequality is from \eqref{eq-A.26}.
\endpf

\subsection{Proof of Theorem \ref{Theorem5}}
\begin{lemma}\label{Lemma1}
Let $Z^{(1)} \geq \ldots \geq Z^{(n)}$ be the order statistics of $n$ independent Pareto random variables with a scale parameter $x_m$ and a shape parameter $\alpha$, which have the following cumulative distribution function:
 \begin{eqnarray*}
	F_Z(x)=
 \begin{cases}
  1-\(x_m/x \)^{\alpha}  & x \geq x_m  \\
 0 & x < x_m,  \\
\end{cases}
\end{eqnarray*}
where $\alpha>2$.
Then, for any given positive constant $a$, we have, with probability at least $1-p^{-a}$,
\begin{equation}\label{eq-A.29}
|\log(Z^{(j)}) - \log(x_m) - \log(n/j)/\alpha| \leq C j^{-1/2} \log p \quad \text{for all}  \, \,  \, 1 \leq j \leq n.
\end{equation}
\end{lemma}

\textbf{Proof of Lemma \ref{Lemma1}.}
By the R\'enyi's representation of exponential order statistics,  $\log(Z^{(j)}/x_m)$ can be represented as $\sum_{l=j}^{n}E_l$, where $E_l$ is independent exponential random variable with the parameter $\alpha l$.
Using the Markov inequality, for each $1 \leq j \leq n$, we have
\begin{eqnarray}\label{eq-A.30}
	&&\Pr \left \{  \sum_{l=j}^{n} \[ E_l - 1/(\alpha l) \] \geq  (a+1+1/\kappa) j^{-1/2} \log p \right \}  \cr
	&&\leq p^{-a-1-1/\kappa}  \prod_{l=j}^{n} \E  \[\exp \left\{ j^{1/2}[E_l - 1/(\alpha l)] \right\}   \]  \cr
	&&\leq p^{-a-1-1/\kappa}  \prod_{l=j}^{n} \dfrac{\alpha l}{\alpha l - j^{1/2}} \exp \left\{-j^{1/2}/(\alpha l) \right\} \cr
	&&\leq p^{-a-1-1/\kappa}  \prod_{l=j}^{n} \exp \left\{j^{1/2}/(\alpha l-j^{1/2})\right\} \exp \left\{-j^{1/2}/(\alpha l) \right\} \cr
	&&\leq p^{-a-1-1/\kappa}  \exp \left\{ \sum_{l=j}^{n} 2j/(\alpha^2 l^2) \right\} \cr
	&&\leq Cp^{-a-1-1/\kappa},
\end{eqnarray}
where the second inequality is from the moment generating function of the exponential distribution and the third inequality is due to the fact that $1+x \leq \exp(x)$ for $x \in \mathbb{R}$. 
Similarly, we can show
\begin{eqnarray}\label{eq-A.31}
	&&\Pr \left \{  \sum_{l=j}^{n} \[ - E_l + 1/(\alpha l) \] \geq  (a+1+1/\kappa) j^{-1/2} \log p \right \}  \cr
	&&\leq Cp^{-a-1-1/\kappa}.
\end{eqnarray}
Also, by the fact that
\begin{equation*}
 \int_{j}^{n} \dfrac{1}{x}dx \leq \sum^{n}_{l=j}\dfrac{1}{l} \leq \int_{j-1}^{n} \dfrac{1}{x}dx \quad \text{for} \quad 2 \leq j \leq n,
\end{equation*}
we have
\begin{equation}\label{eq-A.32}
 \left| \sum^{n}_{l=j}\dfrac{1}{l} - \log(n/j) \right| \leq Cj^{-1} \quad \text{for} \quad 1 \leq j \leq n.
\end{equation}
Combining \eqref{eq-A.30}--\eqref{eq-A.32}, we can obtain \eqref{eq-A.29}. 
\endpf

\textbf{Proof of Theorem \ref{Theorem5}.}
Without loss of generality, we assume that $n = 2L+1 $ and  $u_n = 2u_n^{'}$ for some $L, u_n^{'}  \in \mathbb{ N}$.
Define
\begin{equation*}
G_{i,1}^{(j)} = j \log \(\dfrac{U_{i,1}^{(j)}}{U_{i,1}^{(j+1)}} \) \quad \text{and} \quad G_{i,2}^{(j)} = j \log \(\dfrac{U_{i,2}^{(j)}}{U_{i,2}^{(j+1)}} \)  
\quad \text{for} \quad j=1, \ldots, L-1,
\end{equation*}
where $U_{i,1}^{(1)} \geq \ldots \geq U_{i,1}^{(L)}$ are the order statistics of the sample $| \epsilon_{i}\left(\tau_{i,2k}\right)-\epsilon_{i}\left(\tau_{i,2k-1}\right)|$, $k=1, \ldots, L$, and $U_{i,2}^{(1)} \geq \ldots \geq U_{i,2}^{(L)}$ are the order statistics of the sample $| \epsilon_{i}\left(\tau_{i,2k+1}\right)-\epsilon_{i}\left(\tau_{i,2k}\right)|$, $k=1, \ldots, L$.  
Note that for any $c \geq 0$, $\Pr \left\{ U_{i,1}^{(u_n+1)}  \leq c \right\}$ and $\Pr \left\{ U_{i,2}^{(u_n+1)}  \leq c \right\}$ are determined by $F_i(c)$.
Hence, by Assumption \ref{Assumption4}(a) and Lemma \ref{Lemma1}, we have 
\begin{equation*}
\Pr \left\{ \min_{1 \leq i \leq p} \min_{1 \leq l \leq 2} U_{i,l}^{(u_n+1)}\geq  C \(n/u_n\)^{1/(2\alpha_i)} \right\} \geq 1-p^{-1-a}. 
\end{equation*}
Therefore, without loss of generality, for $1 \leq j \leq u_n +1$, we assume that $U_{i,1}^{(j)}$ and $U_{i,2}^{(j)}$ are the $j$-th largest order statistics of $L$ independent Pareto random variables with a scale parameter $c_i$ and a shape parameter $2\alpha_i$.

Consider $\dfrac{1}{u_n^{'}}\sum_{j=1}^{u_n^{'}} G_{i,1}^{(j)}$. 
By the R\'enyi's representation of exponential order statistics, we can replace $G_{i,1}^{(j)}$ with $E_{j,1}$ for $1 \leq j \leq u_n $, where $E_{j,1}$'s are independent exponential random variables with a parameter $2\alpha_i$.
Thus, similar to the proofs of \eqref{eq-A.30}, we can show
\begin{equation}\label{eq-A.33}
\Pr \left\{ \left| \dfrac{1}{u_n^{'}}\sum_{j=1}^{u_n^{'}} G_{i,1}^{(j)} - \dfrac{1}{2\alpha_i}  \right| \leq  Cu_n^{-1/2}\log p \right\} \geq 1-p^{-2-a}, 
\end{equation}
and similarly, we can show
\begin{equation}\label{eq-A.34}
\Pr \left\{ \left| \dfrac{1}{u_n^{'}}\sum_{j=1}^{u_n^{'}} G_{i,2}^{(j)} - \dfrac{1}{2\alpha_i}  \right| \leq  Cu_n^{-1/2}\log p \right\} \geq 1-p^{-2-a}. 
\end{equation}
Define 
\begin{equation*}
G_{i}^{(j)} = j \log \(\dfrac{U_{i}^{(j)}}{U_{i}^{(j+1)}} \)  \quad \text{for} \quad j=1, \ldots, n-2,
\end{equation*}
where $U_{i}^{(1)} \geq \ldots \geq U_{i}^{(n-1)}$ are the order statistics of the sample $| \epsilon_{i}\left(\tau_{i,k}\right)-\epsilon_{i}\left(\tau_{i,k-1}\right)|$, $k=2, \ldots, n$.
Since   
\begin{equation}\label{eq-A.35}
\min(U_{i,1}^{(j)}, U_{i,2}^{(j)}) \leq  U_{i}^{(2j-1)}, U_{i}^{(2j)} \leq \max(U_{i,1}^{(j)}, U_{i,2}^{(j)}),
\end{equation}
we have
\begin{eqnarray*}
&& \left| \dfrac{1}{u_n}\sum_{j=1}^{u_n} G_{i}^{(j)} - \dfrac{1}{u_n}\sum_{j=1}^{u_n^{'}} G_{i,1}^{(j)} - \dfrac{1}{u_n}\sum_{j=1}^{u_n^{'}} G_{i,2}^{(j)} \right|  \cr
&& \leq \dfrac{1}{u_n}\sum_{j=1}^{u_n^{'}}\left|\log\left(U_{i}^{(2j-1)}\right) +\log\left(U_{i}^{(2j)}\right) - \log\left(U_{i,1}^{(j)}\right) - \log\left(U_{i,2}^{(j)}\right)\right| \cr
&& \quad + \left|\log\left(U_{i}^{(u_n +1)}\right) - \dfrac{1}{2}\log\left(U_{i,1}^{(u_n^{'} +1)}\right) - \dfrac{1}{2}\log\left(U_{i,2}^{(u_n^{'} +1)}\right)\right| \cr
&& \leq  \dfrac{1}{u_n}\sum_{j=1}^{u_n^{'}}  \left| \log\left(U_{i,1}^{(j)}\right) - \log \left(U_{i,2}^{(j)}\right)  \right| +  \dfrac{1}{2}\left| \log\left(U_{i,1}^{(u_n^{'} +1)}\right) - \log \left(U_{i,2}^{(u_n^{'} +1)}\right)  \right|.
\end{eqnarray*}
Therefore, by Lemma \ref{Lemma1}, we can show
\begin{equation}\label{eq-A.36}
\Pr \left\{ \left| \dfrac{1}{u_n}\sum_{j=1}^{u_n} G_{i}^{(j)} - \dfrac{1}{u_n}\sum_{j=1}^{u_n^{'}} G_{i,1}^{(j)} - \dfrac{1}{u_n}\sum_{j=1}^{u_n^{'}} G_{i,2}^{(j)} \right| \leq  Cu_n^{-1/2}\log p \right\} \geq 1-p^{-2-a}.
\end{equation}
Combining \eqref{eq-A.33}, \eqref{eq-A.34}, and \eqref{eq-A.36}, we have
\begin{equation}\label{eq-A.37}
\Pr \left\{ \left| \dfrac{1}{u_n}\sum_{j=1}^{u_n} G_{i}^{(j)} - \dfrac{1}{2\alpha_i}  \right| \leq  Cu_n^{-1/2}\log p \right\} \geq 1-3p^{-2-a}.
\end{equation}

Define
\begin{equation*}
H_{c,i}^{(j)} = j \log \(\dfrac{D_{c,i}^{(j)}}{D_{c,i}^{(j+1)}} \) \quad \text{for} \quad j=1, \ldots, n-2,
\end{equation*}
where $D_{c,i}^{(1)} \geq \ldots \geq D_{c,i}^{(n-1)}$ are the order statistics of  $\{| Y_{i}^c\left(\tau_{i,2}\right)-Y_{i}^c\left(\tau_{i,1}\right)|, \ldots, | Y_{i}^c\left(\tau_{i,n}\right)-Y_{i}^c\left(\tau_{i,n-1}\right)|\}$. 
By Assumption \ref{Assumption1}(a)--(b), we have
 \begin{equation*}
\Pr \left\{ \sup_{1 \leq k \leq n-1}\left| X_i^{c}(\tau_{i,k+1})- X_i^{c}(\tau_{i,k})  \right| \leq  Cn^{-1/2} \sqrt{\log p} \right\} \geq 1-p^{-2-a},
\end{equation*}
which implies
 \begin{equation}\label{eq-A.38}
\Pr \left\{ \max_{1 \leq j \leq n-1}\left| U_{i}^{(j)} - D_{c,i}^{(j)} \right| \leq  Cn^{-1/2} \sqrt{\log p} \right\} \geq 1-p^{-2-a}.
\end{equation}
Also, by Lemma \ref{Lemma1} and \eqref{eq-A.35}, we have
 \begin{equation}\label{eq-A.39}
\Pr \left\{ U_{i}^{(u_n+1)} \geq  C \(n/u_n\)^{1/(2\alpha_i)} \right\} \geq 1-p^{-2-a}.
\end{equation}
Then, using the fact that $x/(1+x) \leq \log(1+x) \leq x$ for all $x>-1$, we can show, with probability at least $1-2p^{-2-a}$,
 \begin{eqnarray}\label{eq-A.40}
&& \left| \dfrac{1}{u_n}\sum_{j=1}^{u_n} \left( G_{i}^{(j)} - H_{c,i}^{(j)}  \right) \right| \cr
&& \leq \dfrac{1}{u_n}\sum_{j=1}^{u_n} \left| \log\left(U_i^{(j)}\right) - \log\left(D_{c,i}^{(j)}\right) \right| +  \left| \log\left(U_i^{(u_n +1)}\right) - \log\left(D_{c,i}^{(u_n +1)}\right) \right| \cr
&& \leq \dfrac{C}{u_n} \sum_{j=1}^{u_n} \left|  \dfrac{U_i^{(j)} - D_{c,i}^{(j)}}{U_i^{(j)}} \right| + C \left|  \dfrac{U_i^{(u_n+1)} - D_{c,i}^{(u_n+1)}}{U_i^{(u_n+1)}} \right| \cr
&& \leq Cn^{-1/2} \sqrt{\log p}.
\end{eqnarray}
Combining \eqref{eq-A.37} and \eqref{eq-A.40}, we have
\begin{equation}\label{eq-A.41}
\Pr \left\{ \left| \dfrac{1}{u_n}\sum_{j=1}^{u_n} H_{c,i}^{(j)} - \dfrac{1}{2\alpha_i}  \right| \leq  Cu_n^{-1/2}\log p \right\} \geq 1-5p^{-2-a}.
\end{equation}

Consider $\dfrac{1}{u_n} \left|\sum_{j=1}^{u_n} H_i^{(j)} \1 \( H_i^{(j)} \leq \omega \)  \right|$.
We have 
\begin{eqnarray}\label{eq-A.42}
&& \dfrac{1}{u_n} \left|\sum_{j=1}^{u_n} \left[H_{c,i}^{(j)} - H_i^{(j)} \1 \( H_i^{(j)} \leq \omega \)\right]  \right| \cr
&& \leq \dfrac{1}{u_n} \left|\sum_{j=1}^{u_n} \left[H_{c,i}^{(j)} - H_{c,i}^{(j)} \1 \( H_{c,i}^{(j)} \leq \omega \)\right]  \right| \cr 
&& \quad + \dfrac{1}{u_n} \left|\sum_{j=1}^{u_n} \left[H_{c,i}^{(j)} \1 \( H_{c,i}^{(j)} \leq \omega \) - H_i^{(j)} \1 \( H_i^{(j)} \leq \omega \)\right]  \right| \cr
&& = (\uppercase\expandafter{\romannumeral1}) +(\uppercase\expandafter{\romannumeral2}).
\end{eqnarray}
For $(\uppercase\expandafter{\romannumeral1})$, since $G_{i,1}^{(j)}$ and $G_{i,2}^{(j)}$ can be replaced with exponential random variables with the parameter $2\alpha_i$, we have
\begin{equation*}
\Pr \left\{ \max_{1 \leq j \leq u_n} \max_{1 \leq l \leq 2}  G_{i,l}^{(j)}  \leq  C \log p  \right\} \geq 1-p^{-2-a},
\end{equation*}
which implies
\begin{equation}\label{eq-A.43}
\Pr \left\{ \max_{1 \leq j \leq u_n}  G_{i}^{(j)}  \leq  C \log p  \right\} \geq 1-p^{-2-a}.
\end{equation}
Now, we investigate $\max\limits_{1 \leq j \leq u_n} \left| H_{c,i}^{(j)} - G_{i}^{(j)}\right|$.
We have
\begin{eqnarray*}
\left|  H_{c,i}^{(j)} - G_{i}^{(j)} \right| &=&  j \left|  \log \left(\dfrac{U_{i}^{(j)} D_{c,i}^{(j+1)}}{U_{i}^{(j+1)} D_{c,i}^{(j)}}\right) \right| \cr
&=&  j \left|  \log \left(1 + \dfrac{U_{i}^{(j)}\left[D_{c,i}^{(j+1)} - U_{i}^{(j+1)}\right] +  U_{i}^{(j+1)}\left[U_{i}^{(j)} - D_{c,i}^{(j)}\right] }{U_{i}^{(j+1)} D_{c,i}^{(j)}}\right) \right|.
\end{eqnarray*}
By \eqref{eq-A.38} and \eqref{eq-A.39}, we have, with probability at least $1-2p^{-2-a}$,
\begin{equation*}
 \max_{1 \leq j \leq u_n}  \left|    \dfrac{U_{i}^{(j)}\left[D_{c,i}^{(j+1)} - U_{i}^{(j+1)}\right] +  U_{i}^{(j+1)}\left[U_{i}^{(j)} - D_{c,i}^{(j)}\right] }{U_{i}^{(j+1)} D_{c,i}^{(j)}} \right| \leq C n^{-(\alpha_i+1)/2\alpha_i} u_n^{1/(2\alpha_i)} \sqrt{\log p}.
\end{equation*}
Thus, using the fact that $x/(1+x) \leq \log(1+x) \leq x$ for all $x>-1$, we can show, with probability at least $1-2p^{-2-a}$,
\begin{equation}\label{eq-A.44}
\max_{1 \leq j \leq u_n} \left| H_{c,i}^{(j)} - G_{i}^{(j)}\right| \leq C n^{-(\alpha_i+1)/2\alpha_i} u_n^{(2\alpha_i+1)/(2\alpha_i)} \sqrt{\log p}.
\end{equation}
Combining \eqref{eq-A.43} and \eqref{eq-A.44}, we have
\begin{equation*}
\Pr \left\{ \max_{1 \leq j \leq u_n}  H_{c,i}^{(j)}  \leq  C \log p  \right\} \geq 1-3p^{-2-a}.
\end{equation*}
Therefore, we have
\begin{equation}\label{eq-A.45}
\Pr \left\{(\uppercase\expandafter{\romannumeral1}) = 0 \right\} \geq 1-3p^{-2-a}.
\end{equation}
For $(\uppercase\expandafter{\romannumeral2})$, since $\pi=0$, with probability at least $1-p^{-2-a}$, the number of jumps is bounded by $C \log p$.
Hence, we have
\begin{equation}\label{eq-A.46}
\Pr \left\{(\uppercase\expandafter{\romannumeral2}) \leq C \(\log p \)^2/u_n \right\} \geq 1-p^{-2-a}.
\end{equation}
By \eqref{eq-A.41}, \eqref{eq-A.42}, \eqref{eq-A.45}, and \eqref{eq-A.46}, we have
\begin{equation}\label{eq-A.47}
\Pr \left\{ \left| \dfrac{1}{u_n}\sum_{j=1}^{u_n}  H_i^{(j)} \1 \( H_i^{(j)} \leq \omega \) - \dfrac{1}{2\alpha_i}  \right| \leq  Cu_n^{-1/2}\log p \right\} \geq 1-9p^{-2-a},
\end{equation}
which completes the proof.
\endpf

\textbf{Proof of Proposition \ref{Proposition2}.}
By Theorem \ref{Theorem1}, Theorem \ref{Theorem3}, and \eqref{eq-5.13}, it is enough to show
\begin{equation}\label{eq-A.48}
 \left( n^{-1/2}\log p \right) ^{(\hat{\alpha}_{ij} -1)/\hat{\alpha}_{ij}} \leq C \left( n^{-1/2}\log p \right) ^{({\alpha}_{ij} -1)/{\alpha}_{ij}}
\end{equation}
under the event
\begin{equation*}
 \left\{ \max_{1 \leq i \leq p} \hat{\alpha}_i - \alpha_i <0 \quad \text{and} \quad  \max_{1 \leq i \leq p} | \hat{\alpha}_i - \alpha_i | \leq  C  n^{-\xi}\log p \right\}.
\end{equation*}
We first investigate $\hat{\alpha}_{ij} - {\alpha}_{ij}$.
We have
\begin{equation*}
\left|  \dfrac {\hat{\alpha}_{i}\hat{\alpha}_{j}}{\hat{\alpha}_{i}+\hat{\alpha}_{j}} - \dfrac {{\alpha}_{i}{\alpha}_{j}}{{\alpha}_{i}+{\alpha}_{j}}\right| = \left| \dfrac{\hat{\alpha}_{i}{\alpha}_{i}\left(\hat{\alpha}_{j}-\alpha_j\right) + \hat{\alpha}_{j}{\alpha}_{j}\left(\hat{\alpha}_{i}-\alpha_i\right)   }{\left(\hat{\alpha}_{i}+\hat{\alpha}_{j}\right)\left({\alpha}_{i}+{\alpha}_{j}\right)} \right| \leq C  n^{-\xi}\log p,
\end{equation*}
which implies
\begin{equation*}
\left| \hat{\alpha}_{ij} - {\alpha}_{ij} \right| \leq C  n^{-\xi}\log p
\end{equation*}
and
\begin{equation*}
\left| (\hat{\alpha}_{ij} -1)/\hat{\alpha}_{ij} - ({\alpha}_{ij} -1)/{\alpha}_{ij}  \right| \leq C  n^{-\xi}\log p.
\end{equation*}
Thus, we have
\begin{eqnarray*}
 && \log\left(n^{1/2}/ \log p \right) \left| (\hat{\alpha}_{ij} -1)/\hat{\alpha}_{ij} - ({\alpha}_{ij} -1)/{\alpha}_{ij}  \right| \cr
 && \leq  C n^{w-\xi} \log n \leq C,
\end{eqnarray*}
which completes the proof.
\endpf

\end{spacing}
\end{document}